\documentclass[11pt]{amsart}

\usepackage{geometry}
\usepackage{fancyhdr}

\usepackage[french,english]{babel}

\usepackage{amsmath,amsfonts,amssymb,amsthm}
\usepackage{bm}
\usepackage[all]{xy}
\usepackage{tikz-cd}

\usepackage{hyperref}

\theoremstyle{plain}
\newtheorem{thm}{Theorem}[section]
\newtheorem{pro}[thm]{Proposition}
\newtheorem{lem}[thm]{Lemma}
\newtheorem{cor}[thm]{Corollary}

\theoremstyle{definition}
\newtheorem{de}[thm]{Definition}
\newtheorem{expl}[thm]{Example}
\newtheorem{rem}[thm]{Remark}

\theoremstyle{plain}
\newtheorem*{conj}{Conjecture}
\newtheorem*{assu}{Assumption}
\newtheorem*{thm*}{Theorem}

\theoremstyle{definition}
\newtheorem*{rem*}{Remark}

\newcommand{\thmref}[1]{Theorem~\ref{#1}}
\newcommand{\propref}[1]{Proposition~\ref{#1}}
\newcommand{\lemref}[1]{Lemma~\ref{#1}}
\newcommand{\corref}[1]{Corollary~\ref{#1}}

\newcommand{\defref}[1]{Definition~\ref{#1}}
\newcommand{\explref}[1]{Example~\ref{#1}}
\newcommand{\remref}[1]{Remark~\ref{#1}}

\newcommand{\RNum}[1]{\uppercase\expandafter{\romannumeral #1\relax}}

\def\Z{\mathbb{Z}}
\def\N{\mathbb{N}}
\def\Q{\mathbb{Q}}
\def\R{\mathbb{R}}
\def\C{\mathbb{C}}

\def\G{\mathbb{G}}
\def\A{\mathbb{A}}
\def\T{\mathbb{T}}

\def\e{\mathbf{e}}

\def\bG{\mathbf{G}}

\def\bX{\mathbf{X}}

\def\LG{{}^L G}
\def\LH{{}^L H}

\def\a{{\mathfrak a}}

\def\g{{\mathfrak g}}
\def\h{{\mathfrak h}}

\def\m{{\mathfrak m}}

\def\t{{\mathfrak t}}
\def\z{{\mathfrak z}}

\newcommand{\cA}{\mathcal{A}}

\newcommand{\cC}{\mathcal{C}}
\newcommand{\cD}{\mathcal{D}}

\newcommand{\cF}{\mathcal{F}}
\newcommand{\cG}{\mathcal{G}}
\newcommand{\cH}{\mathcal{H}}

\newcommand{\cO}{\mathcal{O}}

\newcommand{\cQ}{\mathcal{Q}}
\newcommand{\cR}{\mathcal{R}}
\newcommand{\cS}{\mathcal{S}}

\newcommand{\loccit}{\textit{loc. cit. }}

\newcommand{\dslash}{{/\mkern-6mu/}}

\newcommand\Cat{{\mathbf {Cat}}}
\newcommand\Top{{\mathbf {Top}}}
\newcommand\Gp{{\mathbf {Gp}}}

\newcommand\Gpd{{\mathbf {Gpd}}}
\newcommand\Set{{\mathbf {Sets}}}
\newcommand\sSet{{s\mathbf {Sets}}}

\newcommand\Alg{{\mathbf {Alg}}}

\newcommand\CR{{\mathbf {CR}}}
\newcommand\sCR{{s\mathbf {CR}}}
\newcommand\sCROk{{{}_{\cO} \backslash \sCR /_{k} }}
\newcommand\sCRO{{{}_{\cO} \backslash \sCR  }}
\newcommand\CNL{{ \mathbf {CNL} }}
\newcommand\Art{{ \mathbf {Art} }}
\newcommand\sArt{{{}_{\cO} \backslash s\mathbf {Art} /_{k}}}
\newcommand\sArtB{{{}_{\cO} \backslash s\mathbf {Art} /_{B}}}
\newcommand\Mod{{\mathbf {Mod}}}
\newcommand\sMod{{s\mathbf {Mod}}}
\newcommand{\Ch}{{\mathbf {Ch}}}
\newcommand\sHom{{\mathbf {sHom}}}
\newcommand\shom{{\mathbf {shom}}}
\newcommand{\Ho}{{\mathrm {Ho}}}
\newcommand\hocolim{{\mathrm{hocolim}}}
\newcommand\holim{{\mathrm{holim}}}
\newcommand{\hofib}{{\mathrm {hofib}}}
\newcommand{\DK}{{\mathrm {DK}}}
\newcommand{\Ex}{{\mathrm {Ex}}}
\newcommand{\Bi}{{\mathrm {Bi}}}
\newcommand{\BG}{{\mathcal {B}G}}
\newcommand{\BPG}{{\mathcal {B}PG}}

\newcommand{\uotimes}{{ \underline{\otimes} }}
\newcommand{\lotimes}{{\stackrel{L}{\otimes}}}

\newcommand\SL{{\mathrm {SL}}}

\newcommand\GL{{\mathrm {GL}}}

\newcommand\Sp{{\mathrm {Sp}}}
\newcommand\GSp{{\mathrm {GSp}}}
\newcommand\OO{{\mathrm {O}}}

\newcommand\GSO{{\mathrm {GSO}}}
\newcommand\GO{{\mathrm {GO}}}

\newcommand\U{{\mathrm {U}}}

\newcommand\GSpin{{\mathrm {GSpin}}}

\newcommand{\depth}{{\mathrm{depth}}}
\newcommand{\pd}{{\mathrm{pd}}}

\newcommand{\Ann}{{\mathrm{Ann}}}

\newcommand{\Tor}{{\mathrm {Tor}}}
\newcommand{\Ext}{{\mathrm {Ext}}}
\newcommand\rank{{\mathrm {rank} }}

\newcommand{\Hom}{{\mathrm {Hom}}}
\newcommand{\End}{{\mathrm {End}}}

\newcommand{\iso}{{\> \stackrel{\sim}{\rightarrow} \> }} 
\newcommand\im{{\mathrm {im}}}

\newcommand\divides{{\> \small| \>}}
\newcommand\notdivides{{\> \not\small| \>}}

\newcommand\Frob{{\mathrm {Frob}}}
\newcommand{\Der}{{\mathrm {Der}}}

\newcommand{\Gal}{{\mathrm {Gal}}}

\newcommand\Res{{\mathrm {Res}}}

\newcommand\Lie{{\mathrm {Lie}}}
\newcommand\Def{{\mathrm {Def}}}

\newcommand\crys{{\mathrm {crys}}}
\newcommand\proj{{\mathrm {proj}}}
\newcommand\inj{{\mathrm {inj}}}
\newcommand\rec{{\mathrm {rec}}}

\newcommand\nord{{\mathrm {n.o}}}
\newcommand\ord{{\mathrm {ord}}}
\newcommand\FL{{\mathrm {FL}}}
\newcommand\Resm{{\mathrm {Res}_\m}}
\newcommand\Galm{{\mathrm {Gal}_\m}}
\newcommand\Vanm{{\mathrm {Van}_\m}}

\newcommand\ad{{\mathrm {ad}}}

\newcommand\id{{\mathrm {id}}}

\newcommand\op{{\mathrm {op}}}

\newcommand\ur{{\mathrm {ur}}}

\newcommand\temp{{\mathrm {temp}}}

\DeclareFontFamily{U}{wncy}{}
\DeclareFontShape{U}{wncy}{m}{n}{<->wncyr10}{}
\DeclareSymbolFont{mcy}{U}{wncy}{m}{n}
\DeclareMathSymbol{\Sha}{\mathord}{mcy}{"58} 

\DeclareMathOperator{\diag}{diag}

\begin{document}
	
	\title[Derived deformation rings allowing congruences]{Derived deformation rings allowing congruences}
	\author{Yichang CAI}
	\thanks{The author is partially supported by the grants PerCoLaTor ANR-14-CE25 and CoLoss AAPG2019.}
	\address{LAGA, Institut Galil\'ee, Université Sorbonne Paris Nord, 99 av. J.-B. Cl\'ement, 93430 Villetaneuse, France.}
	\subjclass{}
	\keywords{}
	\email{cai@math.univ-paris13.fr}
	
	\begin{abstract}
		We generalize a result of Galatius and Venkatesh (\cite[Theorem 14.1]{GV18}) which relates the graded module of cohomology of locally symmetric spaces to the graded homotopy ring of the derived Galois deformation rings, by removing certain assumptions, and in particular by allowing congruences inside the localized Hecke algebra.
	\end{abstract}
	
	\maketitle
	
	\tableofcontents

\section*{Introduction}

The cohomology of locally symmetric spaces associated to reductive algebraic groups defined over number fields is a central object in modern Number Theory.  As a complex vector space endowed with an action of the Hecke algebra, it generalizes the space of modular forms for general groups; on the other hand, this Hecke module admits natural underlying integral structures (for instance over rings of $p$-adic integers).  Given a cohomological cuspidal automorphic representation $\pi$, the $\pi$-isotypical component of the cohomology under the Hecke action may occur in several degrees. 
In the Shimura variety case, this phenomenon can be avoided by restricting to tempererd representations, but in general, it cannot be avoided.  This phenomenon has been explained over $\C$ by Borel and Wallach by calculations of $(\g,K)$-cohomology.  More recently, a motivic interpretation of this phenomenon has been investigated by A. Venkatesh.  Over the $p$-adic integers, the fundamental works are those by Calegari and Geraghty \cite{CG18} and by Galatius and Venkatesh \cite{GV18}.  Following these works, we will study the relation between the graded structure of cohomology of locally symmetric spaces and the graded homotopy ring of the derived Galois deformation rings under assumptions similar but lighter than those of \cite{GV18}.  

Let $F$ be a fixed number field.  We use the bold $\bG$ to denote a connected reductive linear algebraic group over $F$ and write $G = {}^L \bG = \widehat{\bG} \rtimes \Gal(\bar{F}/F)$ for convenience, as the paper will mainly focus on the $L$-group.  Let $\bG_f = \bG(\A_F^{\infty})$ and $\bG_{\infty} = \bG(F\otimes_{\Q} \R)$.  Let $X_\bG = \bG_{\infty}/K_{\infty}$ 
be the symmetric space associated to $\bG$, where $K_{\infty} = C_{\infty} \cdot A(\R)$, $C_{\infty}$ is a maximal compact subgroup of the real Lie group $\bG_{\infty}$ and $A$ is a maximal $\Q$-split torus of the center of 
$\Res_{\Q}^F \bG$.  Let $q_0$ and $\ell_0$ be integers associated to $\bG$ such that 
\begin{displaymath}
\left\{\begin{array}{l}{  2q_0+\ell_0 = \dim X_\bG=d;  } \\ { \ell_0 = \rank \, \bG_{\infty} - \rank \, K_{\infty}.}\end{array}\right.
\end{displaymath}
For an open compact subgroup $U \subseteq \bG_f$, the locally symmetric space of $\bG$ with level structure $U$ is defined to be $X_\bG^U = \bG(F) \backslash (X_\bG \times \bG_f /U)$.  

Let $p>2$ be an odd prime number.  Let $K$ be a large enough $p$-adic number field containing all embeddings of $F$ into $\bar{\Q}_p $, let $\cO$ be its ring of integers, $k$ be its residue field and $\varpi$ be a uniformizing parameter.  
For a dominant weight $\lambda = (\lambda_{\tau,i})_{ \tau\colon F \hookrightarrow K,1\leq i \leq n }$ for $G$, we write $V_{\lambda} = \otimes_{\tau\colon F \hookrightarrow K} V_{\lambda_{\tau}}$ 
for the irreducible algebraic representation of $G$ of highest weight $\lambda$, and write $\widetilde{V}_{\lambda}(R)$ for the associated sheaf for an $\cO$-algebra $R$.  

The complex cohomology of $X_\bG^U$ can be studied via the theory of $(\g, K)$-cohomology.  In particular, the tempered part $H^{\ast}_{\temp}(X_\bG^U, \widetilde{V}_{\lambda}(\C))$ (an embedding $K\hookrightarrow \C$ is fixed) is concentrated in the interval $[q_0,q_0+\ell_0]$ and we have 
$$\dim H^{q_0+i}_{\temp}(X_\bG^U, \widetilde{V}_{\lambda}(\C)) = \binom{\ell_0}{i} \cdot \dim H^{q_0}_{\temp}(X_\bG^U, \widetilde{V}_{\lambda}(\C)). $$
In fact, in \cite[Section 3]{PV16}, the authors constructed an action of $\wedge^{\ast} \a_\bG^{\ast}$ on $H^{\ast}_{\temp}(X_\bG^U, \widetilde{V}_{\lambda}(\C))$, where $\a_\bG^{\ast}$ 
is the dual of the Lie algebra of the split part of a fundamental Cartan algebra, such that $H^{d-\ast}_{\temp}(X_\bG^U, \widetilde{V}_{\lambda}(\C))$ is freely generated in degree $q_0$ over $\wedge^{\ast} \a_\bG^{\ast}$.  

It's natural to consider the analogous question for integral coefficients.  Under some assumptions, the Calegari-Geraghty method (see \cite{CG18}) implies that, $H^{\ast}(X_\bG^U, \widetilde{V}_{\lambda}(\cO))_{\m}$, where $\m$ is a non-Eisenstein maximal ideal of the associated Hecke algebra, is a free graded module over a graded commutative ring which arises naturally in the Taylor-Wiles method.  However, this graded commutative ring is not canonically defined, and the idea of \cite{GV18} is that the better object is the derived generalization of the Galois deformation ring. We will now explain in more details the backgrounds and results of this paper.

We suppose that $p$ is very good for $G$ in the sense of \cite[Page 10]{BHKT19} and $\zeta_p \notin F$.  Let $S_p$ be the set of places of $F$ dividing $p$ and $S_{\infty}$ be the set of archimedean places of $F$.  Let $S \supseteq S_p$ be a finite set of finite places of $F$.  We write $\bG_S = \prod_{v\in S}\bG(F_v)$ and $\bG^S$ for the image of the natural projection $\bG_f \to \prod_{v\notin S}\bG(F_v)$.  
Let's fix a faithful representation $\bG \to \GL_N$ and define $\underline{\bG}$ to be the schematic closure of $\bG$ in $\GL_{N,\cO_F}$.  Suppose $U = U_S \times U^S = (\prod_{v\in S} U_v) \times (\prod_{v\notin S} U_v)$
with $U_v \subseteq \underline{\bG}(\cO_v)$ for every finite place $v$ and each $U_v$ ($v\notin S \backslash S_p$) hyperspecial maximal; the spherical Hecke algebra $\cH(\bG^S, U^S)$ acts on $H^{\ast}(X_\bG^U, \widetilde{V}_{\lambda}(\cO))$.  Note that the image of this action, which we denote by $h$, is a finite commutative $\cO$-algebra.  We say that a maximal ideal $\m$ is non-Eisenstein if any ($h \otimes_{\cO} k$)-isotypical component appearing in $H^{\ast}(X_\bG^U, \widetilde{V}_{\lambda}(\cO))_{\m} \otimes_{\cO} k$ doesn't come from 
$H^{\ast}(X_\bG^U, \widetilde{V}_{\lambda}(k)) /H^{\ast}_{!}(X_\bG^U, \widetilde{V}_{\lambda}(k))$.  Let $\m$ be a non-Eisenstein maximal ideal of $h$ and let $\T=h_\m$.
Let $\pi$ be a cuspidal automorphic representation occuring in 
$H^{\ast}(X_\bG^U, \widetilde{V}_{\lambda}(\cO))_{\m}$.  

Let $\Gamma_v = \Gal(\bar{F}_v/F_v)$ and let $\Gamma_S$ be the Galois group of the maximal $S$-ramified extension of $F$.  We make the following assumption: 

\begin{assu}[$\Resm$]
	There exists an absolutely irreducible (see \cite[Definition 3.5]{BHKT19}) Galois representation $\bar\rho \colon \Gamma_S \to G(k)$ associated to $\pi$ (see \cite[Section 5]{BG10}) such that 
	\begin{enumerate}
		\item for $v\notin S$, the $\widehat{\bG}(k)$-conjugacy class of $\bar\rho(\Frob_v)$ is given by the Satake parameter of $\pi_v$ modulo $\m$;
		\item $\bar\rho\vert_{\Gamma_v}$ is minimal for $v\in S\backslash S_p$;
		\item $\bar\rho\vert_{\Gamma_v}$ is simultaneously either ordinary, or Fontaine-Laffaille with Hodge–Tate weights differing by at most $p-2$ for $v\in S_p$.  In the ordinary case, $\bar\rho\vert_{\Gamma_v}$ is furthermore assumed to be regular and dual regular (see \cite[Propostion 6.2 and Propostion 6.3]{Til96}).
	\end{enumerate}
	We require further that $\bar\rho$ is odd (see \defref{odd}) and has an enormous image (see \defref{REI}).
\end{assu}

Let $\cS$ be the global deformation problem for $\bar\rho \colon \Gamma_S \to G(k)$, which is either minimal ordinary or minimal Fontaine-Laffaille.  Then the deformation functor for $\bar\rho$ of type $\cS$ (denoted $\cD_{\cS}$) is represented by a complete Noetherian local $\cO$-algebra $R_{\cS}$.  

The method of \cite{CG18} relies significantly on the following conjectures:

\begin{conj}[$\Vanm$]
	The cohomology group $H^i(X_\bG^U, \widetilde{V}_{\lambda}(\cO))_{\m}$ vanishes unless $i\in [q_0,q_0+\ell_0]$.  
\end{conj}

\begin{conj}[$\Galm$]
	There is a Galois representation $\rho_{\m}\colon \Gamma_S \to G(\T)$ lifting $\bar\rho$ such that 
	\begin{enumerate}
		\item $\rho_{\m}\vert_{\Gamma_v}$ is minimal for $v\in S\backslash S_p$;
		\item $\rho_{\m}\vert_{\Gamma_v}$ is simultaneously either ordinary, or Fontaine-Laffaille for every $v\in S_p$;
		\item $\rho_{\m}\vert_{\Gamma_v}$ satisifies local-global compatibility for any Taylor-Wiles prime $v$. 
	\end{enumerate}
	(see \cite[Assumption 2]{GV18} and \cite[Conjecture 6.27]{KT17}).  This means that there is a natural morphism $R_{\cS} \to \T$ in $\CNL_{\cO}$ and similarly for the "Taylor-Wiles thickenings" $R_{Q}$ and $\T_{Q}$ of the rings $R_{\cS}$ and $\T$. 
\end{conj}

\begin{rem*}
	For $\GL_N$ over CM fields, \cite{ACC+18} gives strong evidences for the existence of $R_{\cS} \to \T$, as these authors do prove it after quotient by a nilpotent ideal of $\T$. 
\end{rem*}

Then Calegari and Geraghty (\cite{CG18}) constructed $R_{\infty} = \cO[[ X_1, \dots, X_g ]]$ and $S_{\infty} = \cO[[X_1, \dots, X_{g+\ell_0}]]$ ($g$ is a constant) with an $\cO$-algebra morphism $S_{\infty} \to R_{\infty}$, as well as a complex $C_{\infty}^{\ast}$ of finite free $S_{\infty}$-modules concentrated in degrees $[q_0,q_0+\ell_0]$ and an $S_{\infty}$-algebra morphism $R_{\infty} \to \End_{S_{\infty}}(H^{\ast}(C_{\infty}^{\ast}))$, such that $H^{\ast}(C_{\infty}^{\ast} \otimes_{S_{\infty}} \cO) \cong H^{\ast}(X_\bG^U, \widetilde{V}_{\lambda}(\cO))_{\m}$ and the following result holds: 

\begin{thm*}[Calegari-Geraghty]
	Let the notations be as above.  Assume $(\Resm)$, $(\Galm)$ and $(\Vanm)$.  Then  
	\begin{enumerate}
		\item $H^i(C_{\infty}^{\ast}) = 0$ for $i\neq q_0+\ell_0$ and $H^{q_0+\ell_0}(C_{\infty}^{\ast})$ is free over $R_{\infty}$.
		\item There is an isomorphism $H^{q_0+\ell_0-i}(X_\bG^U, \widetilde{V}_{\lambda}(\cO))_{\m} \cong \Tor_i^{S_{\infty}}(H^{q_0+\ell_0}(C_{\infty}^{\ast}),\cO),$ and $\Tor_{\ast}^{S_{\infty}}(H^{q_0+\ell_0}(C_{\infty}^{\ast}),\cO)$ is a natural graded $\Tor_{\ast}^{S_{\infty}}(R_{\infty},\cO)$-module freely generated by $\Tor_0^{S_{\infty}}(H^{q_0+\ell_0}(C_{\infty}^{\ast}),\cO)$.
		\item $R_{\cS} \to \T$ is an isomorphism.
	\end{enumerate}
\end{thm*}

Note the isomorphism $\Tor_{\ast}^{S_\infty}(R_\infty, \cO) \cong \pi_{\ast} (R_{\infty} \uotimes_{S_{\infty}} \cO)$ as graded commutative $\cO$-algebras, where $- \uotimes_{S_{\infty}} \cO$ can be thought of as a model for calculating the total left derived functor of the degreewise-extended tensor on simplicial rings (see Section \ref{simplicial commutative ring}), so it's natural to think of $\Tor_{\ast}^{S_{\infty}}(R_{\infty},\cO)$ as the homotopy of some derived generalization of $R_{\cS}$.  It is observed in \cite{GV18} that this is indeed the case.  

For a complete and cocomplete category $\cC$, the simplicial category $s\cC$ is defined to be the category of contravariant functors from $\bm{\Delta}$ to $\cC$, where $\bm{\Delta}$ is the cosimplicial indexing category 
(the objects are totally ordered sets $[n]=\{0,\ldots,n\}$ and morphisms are non-decreasing maps).  When $\cC$ is the category of sets, modules or algebras, the category $s\cC$ is naturally a simplicial model category.  
In particular in these categories  
\begin{enumerate}
	\item we can define homotopy groups and a weak equivalence relation, such that $f\colon A\to B$ is a weak equivalence if and only if $f$ induces isomorphisms on all homotopy groups;
	\item there is an enriched hom $\sHom(A,B) \in \sSet$, with the property $\sHom(A,B)_0 \cong \Hom(A,B)$.  
\end{enumerate}

Note $\cD_{\cS}$ restricts to a functor from the category of artinian local $\cO$-algebras $\Art_{\cO}$ to the category of sets $\Set$.  Following \cite{GV18}, $\cD_{\cS}$ can be extended to a functor $s\cD_{\cS}$ from the category of simplicial artinian local $\cO$-algebras $\sArt$ to the category of simplicial sets $\sSet$.  By saying extended, we mean $\cD_{\cS}(A) \cong \pi_0 s\cD_{\cS}(A)$ when $A$ is a classical artinian local $\cO$-algebra (on the right hand side $A$ is regarded as a constant object in $\sArt$).

It is proved that the functor $s\cD_{\cS}$ is pro-representable.  More precisely, there exists a projective system $\cR_{\cS} = (\cR_n)_{n\in \N}$ with each $\cR_n \in \sArt$ being cofibrant, 
such that $s\cD_{\cS}(A)$ is weakly equivalent to $\varinjlim\limits_n\sHom_{\sCROk}(\cR_n, A)$ for each $A\in \sArt$.  Note $\cR_{\cS}$ is unique only in the homotopy category, nonetheless $\pi_{\ast} \cR_{\cS}$ is well-defined.  Indeed, by regarding $\pi_{\ast} \cR_{\cS}$ as the projective limit, it is naturally a graded commutative $\cO$-algebra, and at degree $0$ we have $\pi_0 \cR_{\cS} \cong R_{\cS}$.  We can now state our main result (see \thmref{main theorem}), which is a generalization of \cite[Theorem 14.1]{GV18}:

\begin{thm*}
	With the above notations, there is an isomorphism of graded commutative $\cO$-algebras $\pi_{\ast}\cR_{\cS} \cong \Tor_{\ast}^{S_{\infty}}(R_{\infty},\cO)$ (where $\pi_{\ast} \cR_{\cS}$ is defined as the projective limit).  Moreover, $H^{\ast}(X_\bG^U, \widetilde{V}_{\lambda}(\cO))_{\m}$ is a graded $\pi_{\ast}\cR_{\cS}$-module freely generated by $H^{q_0+\ell_0}(X_\bG^U, \widetilde{V}_{\lambda}(\cO))_{\m}$.  
\end{thm*}

We mention the differences with \cite[Theorem 14.1]{GV18}: 
\begin{enumerate}
	\item In \cite[Theorem 14.1]{GV18} the group $G$ is assumed to have a trivial center.  In the general case, as already pointed out in \cite{GV18}, one has to modify the derived (local and global) universal deformation functors to take the center into account.  
	\item More importantly, one has to redefine the derived local deformation problems, for in \cite[Section 9]{GV18} it is assumed that the classical local (unframed) deformation functors are represented by formally smooth rings, which is not the case for us.  
	\item In \cite{GV18}, only the case $R_{\cS} = \T = \cO$ is considered since the application in \cite[Section 15]{GV18} uses the surjectivity of the homomorphism $S_{\infty} \twoheadrightarrow R_{\infty}$ (see \cite[Remark 1.1]{GV18}).  This surjectivity is obtained by imposing strong restrictions on the local deformation conditions (\cite[Section 10]{GV18}) which we don't have.  Here, we have to recalculate the Poitou-Tate Euler characteristics in order to verify \cite[Theorem 11.1]{GV18} in our more general setting.  See also \cite{TU21}, where some partial results are proved without the surjection $S_{\infty} \twoheadrightarrow R_{\infty}$.
\end{enumerate}

Here is the outline of the paper.  In Section 1, we will recall the Calegari-Geraghty method which describes the graded structure of the integral cohomology after non-Eisenstein localizations.  In Section 2, we will prepare the necessary backgrounds on simplicial theory to study functors from simplicial Artinian $\cO$-algebras to simplicial sets.  In Section 3, we will extend the classical deformation functors to simplicial categories and study the homotopy of their pro-representing rings.  The main result (\thmref{main theorem}) is a generalization of \cite[Theorem 14.1]{GV18}, where in particular the congruences inside the localized Hecke algebra are allowed.  In Section 4, we will discuss the examples of general linear groups and orthogonal similitude groups, and we will try to compare the derived deformation rings and the cohomology of locally symmetric spaces under certain Langlands transfers.

\section{Calegari-Geraghty method}\label{Calegari-Geraghty}

In this section we will present the Calegari-Geraghty method with an emphasis on the graded structures of $H^{\ast}(X_\bG^U, \widetilde{V}_{\lambda}(\cO))_{\m}$ and $\Tor_{\ast}^{S_{\infty}}(R_{\infty},\cO)$.  As mentioned in the introduction, the isomorphism $\Tor_{\ast}^{S_\infty}(R_\infty, \cO) \cong \pi_{\ast} (R_{\infty} \uotimes_{S_{\infty}} \cO)$ provides a motivation for considering deformations to simplicial rings.

We suppose the center $Z$ of $G$ is smooth over $\cO$.  Let $\g_k=\Lie(G/\cO) \otimes_{\cO} k$ (resp. $\z_k=\Lie(Z/\cO)\otimes_{\cO} k$).

Let $\Gamma_S$ be the Galois group of the maximal $S$-ramified extension of $F$ and let $\bar\rho \colon \Gamma_S \to G(k)$ be a fixed absolutely irreducible continuous Galois representation; $\bar\rho$ will eventually be the representation described in $(\Resm)$.  Since $\bar\rho$ is absolutely irreducible, we have $H^0(\Gamma_S, \g_k) = \z_k$ (see \cite[Lemma 5.1]{BHKT19}).

\subsection{Galois deformation theory}

We begin by recalling some deformation theory for $\bar\rho$.  Let $\CNL_{\cO}$ be the category of complete Notherian local $\cO$-algebras with residue field $k$.  The universal framed deformation functor $\Def_S^{\square}\colon \CNL_{\cO} \to \Set$ for $\bar\rho$ is defined by associating $A\in\CNL_{\cO}$ to the set of continuous liftings $\rho\colon \Gamma_S\to G(A)$ which make the following diagram commute:
\begin{displaymath}
\xymatrix{
	\Gamma_S \ar[r]^{ \rho }\ar[dr]^{\bar\rho} & G(A) \ar[d]\\
	& G(k).
}
\end{displaymath}
Moreover, the universal deformation functor $\Def_S \colon \CNL_{\cO} \to \Set$ is defined by associating $A\in\CNL_{\cO}$ to the set of $\ker(G(A) \to G(k))$-conjugacy classes of $\Def_S^{\square}(A)$.  As an application of Schlessinger's criterion (see \cite[Theorem 2.11]{Sch68}), the functors $\Def^{\square}_S$ and $\Def_S$ are representable (for the latter we require the condition $H^0(\Gamma, \g_k) = \z_k$, see \cite[Theorem 3.3]{Til96}).

For each place $v\in S$, we define similarly the universal framed deformation functor $\Def_v^{\square}$ and the universal deformation functor $\Def_v$ for $\bar\rho\vert_{\Gamma_v}$ (recall $\Gamma_v = \Gal(\bar{F}_v/F_v))$.  Again Schlessinger's criterion implies that $\Def_v^{\square}$ is representatble, say by $R_v^{\square} \in \CNL_{\cO}$.  However the functor $\Def_v$ is generally not representable, for $H^0(\Gamma_v, \g_k)= \z_k$ is usually not true. 

\begin{de}
	Let $v$ be a finite place of $F$.  A local deformation problem for $\bar\rho\vert_{\Gamma_v}$ is a subfunctor $\cD_v$ of $\Def_v^{\square}$ satisfying the following conditions:
	\begin{enumerate}
		\item $\cD_v$ is represented by a quotient $R_v \in \CNL_{\cO}$ of $R_v^{\square}$.
		\item For any $A\in \CNL_{\cO}$, $\rho \in \cD_v(A)$ and $a \in \widehat{G}(A)$, we have $a \rho a^{-1} \in \cD_v(A)$.
	\end{enumerate}
\end{de}

Let $k[\epsilon] = k[t]/(t^2)$.  Then it's well-known that $\cD_v(k[\epsilon])$ can be identified with a subspace $\widetilde{L}_v \subseteq Z^1(\Gamma_v, \g_k)$, which is the preimage of a subspace $L_v \subseteq H^1(\Gamma_v, \g_k)$ under the projection $Z^1(\Gamma_v, \g_k) \to H^1(\Gamma_v, \g_k)$.  Note $R_v$ can be generated by $$\dim_k \widetilde{L}_v = \dim_k \g_k + (\dim_k L_v - \dim_k H^0(\Gamma_v,\g_k))$$ 
variables over $\cO$.  We say $\cD_v$ or $R_v$ is formally smooth if $R_v$ is a power series ring over $\cO$; note then the number of generators is $\dim_k \widetilde{L}_v$.

\begin{de}\label{global deformation problem}
	A global deformation problem is a tuple $\cS =(S, \{\cD_v\}_{v \in S}),$ where $\cD_v$ is a local deformation problem for $\bar\rho\vert_{\Gamma_v}$ for each $v\in S$.
\end{de}

\begin{de}
	We say a lifting $\rho \colon \Gamma_S \to G(A)$ ($A\in \CNL_{\cO}$) of $\bar\rho$ is of type $\cS$ if $\rho\vert_{\Gamma_v} \in \cD_v(A)$ for every $v\in S$.  Two liftings $\rho_1$,$\rho_2 \colon \Gamma_S \to G(A)$ of type $\cS$ are said to be equivalent if there exists $a\in \ker(G(A) \to G(k))$ such that $\rho_2 = a \rho_1 a^{-1}$.  An equivalent class of liftings of type $\cS$ is called a deformation of type $\cS$.  We denote by $\cD_{\cS} \colon \CNL_{\cO} \to \Set$ the functor which sends $A\in \CNL_{\cO}$ to the set of deformations to $G(A)$ of type $\cS$.
\end{de}

Under our condition $H^0(\Gamma_S, \g_k) = \z_k$, it's well-known that the functor $\cD_{\cS}$ is representable, say by $R_{\cS}\in \CNL_{\cO}$.  

We define $C^{\ast}_{\cS}(\Gamma_S, \g_k)$ by the cone construction: let $C^{\ast}_{\cS}(\Gamma_S, \g_k) = C^{\ast}[-1]$, where $C^{\ast}$ is the mapping cone of the natural morphism 
\begin{displaymath}
\xymatrix{
	0 \ar[r]& C^0(\Gamma_S, \g_k) \ar[r]\ar[d] & C^1(\Gamma_S, \g_k)  \ar[r]\ar[d] & C^2(\Gamma_S, \g_k) \ar[d]\ar[r] &\dots \\
	0 \ar[r] & 0 \ar[r] & \bigoplus_{v \in S} C^1(\Gamma_v, \g_k)/\widetilde{L}_v  \ar[r] & \bigoplus_{v \in S} C^2(\Gamma_v, \g_k) \ar[r] &\dots 
}
\end{displaymath}

Let $H_{\cS}^{\ast}(\Gamma_S, \g_k)$ be the cohomology of $C^{\ast}_{\cS}(\Gamma_S, \g_k)$.  Then we have the following exact sequence:
\begin{align*}
0 \to &H^0_{\cS}(\Gamma_S, \g_k) \to H^0(\Gamma_S, \g_k) \to 0 \\
\to &H^1_{\cS}(\Gamma_S, \g_k) \to H^1(\Gamma_S, \g_k) \to  \bigoplus_{v \in S} H^1(\Gamma_v, \g_k)/L_v \\
\to &H^2_{\cS}(\Gamma_S, \g_k) \to H^2(\Gamma_S, \g_k) \to \bigoplus_{v \in S} H^2(\Gamma_v, \g_k) \\
\to &H^3_{\cS}(\Gamma_S, \g_k) \to 0.
\end{align*}

For a finite $\cO$-module $M$ equipped with a Galois group action, we write $M^{\vee}=\Hom_{\cO}(M,K/\cO)$ and $M^{\ast}=\Hom_{\cO}(M,K/\cO(1))$.  In particular, if $M$ is a $k$-vector space, then $M^{\vee}\cong \Hom_k(M,k)$ and $M^{\ast}\cong \Hom_k(M,k(1))$.

Define $H^1_{\cS^{\perp}}(\Gamma_S, \g_k^{\ast}) = \ker( H^1(\Gamma_S, \g_k^{\ast})\to \bigoplus_{v \in S}  H^1(\Gamma_v, \g_k^{\ast})/L_v^{\perp} ) $, where $L_v^{\perp} \subseteq H^1(\Gamma_v, \g_k^{\ast})$ is the dual of $L_v \subseteq H^1(\Gamma_v,\g_k)$ under the local Tate duality.  As an application of the Poitou-Tate duality, we have $H^1_{\cS^{\perp}}(\Gamma_S, \g_k^{\ast})^{\vee} \cong H^2_{\cS}(\Gamma_S, \g_k)$ and $H^0(\Gamma_S, \g_k^{\ast})^{\vee} \cong H^3_{\cS}(\Gamma_S, \g_k)$ (see the proof of \cite[Proposition 6.2.24]{ACC+18}).

\begin{lem}\label{relative dimension}
	There is an $\cO$-algebra surjection $\cO[[X_1, \dots, X_g]] \twoheadrightarrow R_{\cS}$, with 
	\begin{align*}
	g = \dim_k H^1_{\cS}(\Gamma_S, \g_k) = &\dim_k H^1_{\cS^{\perp}}(\Gamma_S, \g_k^{\ast}) + \dim_k H^0(\Gamma_S, \g_k) - \dim_k H^0(\Gamma_S, \g_k^{\ast}) \\ &- \sum\limits_{v \divides \infty} \dim_k H^0(\Gamma_v,\g_k) + \sum\limits_{v\in S} (\dim_k L_v - \dim_kH^0(\Gamma_v,\g_k) ).
	\end{align*}
\end{lem}  

\begin{proof}
	See \cite[Proposition 6.2.24]{ACC+18}.
\end{proof}

\begin{rem}\label{relative dimension remark}
	Suppose $\zeta_p \notin F$, and suppose $H=\bar\rho(\Gal_{F(\zeta_p)})$ satisfies $\g_k^H = \z_k$ (this is part of the enormous image condition for $\bar\rho$), then it's easy to see $H^0(\Gamma_S, \g_k^{\ast})=0$.  
\end{rem}

\subsection{Taylor-Wiles primes}

Given $\cS$ and a finite set of places $Q$ disjoint from $S$, we write $\cS_Q=(S \cup Q, \{\cD_v\}_{v \in S \cup Q})$ where $\cD_v = \Def_v^{\square}$ for every $v \in Q$.  

\begin{de}\label{allowable Taylor-Wiles}
	\begin{enumerate}
		\item A place $v\notin S$ is called a Taylor-Wiles prime if $N(v) \equiv 1 \pmod p$ and $\bar\rho(\Frob_v)$ is conjugated to a strongly regular element of $T(k)$ ({\it i.e.}, an element $t\in T(k)$ whose centralizer in $G$ coincides with $T$).
		\item An allowable Taylor-Wiles datum of level $m$ is a set of Taylor-Wiles primes $Q=(v_1, \dots, v_r)$, together with a strongly regular element $t_{v_i} \in T(k)$ conjugate to $\bar\rho(\Frob_{v_i})$ for each $i \in \{1,\dots, r\}$, such that
		\begin{enumerate}
			\item $N(v_i) \equiv 1 \pmod{p^m}$, for every $i \in \{1,\dots, r\}$;
			\item $H^2_{\cS_Q}(\Gamma_{S\cup Q}, \g_k) = 0$.
		\end{enumerate}
	\end{enumerate}
\end{de}

\begin{rem}
	By the Poitou-Tate duality, condition (b) is equivalent to $H^1_{\cS^{\perp}}(\Gamma_S, \g_k^{\ast})=0$.
\end{rem}

The existence of Taylor-Wiles data relies on the enormous image assumption for $\bar\rho$ (see \cite[Definition 6.2.28]{ACC+18}):

\begin{de}\label{REI}
	Let $\g_k^{\prime}$ be the Lie algebra of the derived group $G^{\prime}$.  We say $\bar\rho\colon \Gamma_S \to G(k)$ has an enormous image, if $H = \bar\rho(\Gal_{F(\zeta_p)})$ satisfies the following:
	\begin{enumerate}
		\item $H$ has no non-trivial $p$-power order quotient.
		\item $H^0(H, \g_k^{\prime}) = H^1(H, \g_k^{\prime})=0$.
		\item For any simple $k[H]$-module $W\subseteq \g_k^{\prime}$, there is a regular semisimple $h\in H$ such that $W^h \neq 0$.
	\end{enumerate}
\end{de}

\begin{lem}
	Suppose $\bar\rho\colon \Gamma_S \to G(k)$ has an enormous image.  Let $r\geq \dim_k H^1_{\cS}(\Gamma_S, \g_k^{\ast})$ and $m\geq 1$.  Then there exists an allowable Taylor-Wiles datum $Q$ of level $m$ and cardinal $r$.
\end{lem}

\begin{proof}
	\cite[Lemma 6.2.31]{ACC+18} proved this for $\GL_N$, but the proof applies verbatim for general $G$.
\end{proof}

Now fix $r\geq \dim_k H^1_{\cS}(\Gamma_S, \g_k^{\ast})$.  Let $\cQ=(Q_m)_{m \geq 1}$ be a system of disjoint allowable Taylor-Wiles data, such that each $Q_m$ is of level $m$ and cardinal $r$.  For simplicity, we write $\Gamma_m = \Gamma_{S\cup Q_m}$, $\cD_m = \cD_{\cS_{Q_m}}$ for the deformation functor of type $\cS_{Q_m}$ and $R_m = R_{\cS_{Q_m}}$ for the representing ring of $\cD_m$.  Let $n = \rank \, G$.

\begin{lem}\label{relative dimension 1}
	For every $m\geq 1$, there is an $\cO$-algebra surjection $\cO[[X_1, \dots, X_g]] \twoheadrightarrow R_m$, with 
	\begin{align*}
	g = \dim_k H^1_{\cS_{Q_m}}(\Gamma_m, \g_k) = nr +  \dim_k H^1_{\cS}(\Gamma_S, \g_k) - \dim_k H^1_{\cS^{\perp}}(\Gamma_S, \g_k^{\ast}).
	\end{align*}
\end{lem}

\begin{proof}
	We apply \lemref{relative dimension} to the global deformation problem $\cS_{Q_m}$.  Then $H^1_{\cS_{Q_m}^{\perp}}(\Gamma_m, \g_k^{\ast}) = 0$ by the definition of $Q_m$, so 
	\begin{align*}
	g = \dim_k H^0(\Gamma_S, \g_k) -\sum\limits_{v \divides \infty} \dim_k H^0(\Gamma_v,\g_k) &+ \sum\limits_{v\in S} (\dim_k L_v - \dim_k H^0(\Gamma_v,\g_k) )  \\ &+ \sum\limits_{v\in Q_m} (\dim_k L_v - \dim_k H^0(\Gamma_v,\g_k) ).
	\end{align*}
	
	By \lemref{relative dimension}, we have 
	\begin{align*}
	\dim_k H^1_{\cS}(\Gamma_S, \g_k) - \dim_k H^1_{\cS^{\perp}}(\Gamma_S, \g_k^{\ast}) = &\dim_k H^0(\Gamma_S, \g_k) -\sum\limits_{v \divides \infty} \dim_k H^0(\Gamma_v,\g_k) \\ &+ \sum\limits_{v\in S} (\dim_k L_v - \dim_k H^0(\Gamma_v,\g_k) ).
	\end{align*}
	On the other hand, for $v\in Q_m$, we have $L_v = H^1(\Gamma_v,\g_k)$ and hence 
	$$\dim_k L_v - \dim_k H^0(\Gamma_v,\g_k) = \dim_k H^0(\Gamma_v,\g_k^{\ast}) = n$$ 
	(here the first equality follows from the local Euler characteristic formula and the second equality is because $N(v) \equiv 1 \pmod{p}$ and $\bar\rho(\Frob_v)$ is conjugated to a strongly regular element of $T(k)$).  So the conclusion follows.
\end{proof}

\begin{de}\label{odd}
	We say $\bar\rho$ is odd, if $\sum\limits_{v \divides \infty} \dim_k H^0(\Gamma_v,\g_k) = \ell_0 + [F:\Q] (\dim G -\dim B) + \dim_k H^0(\Gamma_S, \g_k)$.
\end{de}

\begin{rem}
	This definition seems rather deliberate.  For the locally symmetric space associated to $\Res_{\Q}^F \GL_N$ and $\bar\rho \colon \Gamma_S \to \GL_N(k)$, one has $\ell_0 = [\frac{N+1}{2}] r_1  + N r_2 -1$, where $r_1$ (resp. $r_2$) is the numbers of real (resp. complex) places of $F$, and hence 
	$$\ell_0 + [F:\Q] (\dim G -\dim B) + \dim_k H^0(\Gamma_S, \g_k) = [\frac{N+1}{2}] r_1  + N r_2 -1 + (r_1 +2 r_2) \frac{N^2-N}{2}+1.$$
	Therefore the oddness of $\bar\rho$ means precisely $\dim_k H^0(\Gamma_v,\g_k) = [\frac{N^2+1}{2}]$ for every real place $v$, or in other words, $H^0(\Gamma_v,\g_k)$ has the minimal possible dimension.
\end{rem}

Write $\rho_m \colon \Gamma_m \to G(R_m)$ be any representative of the universal deformation.  Then for each $v \in Q_m$, there exists a conjugation of $\rho_m\vert_{\Gamma_v}$ which takes values in $T(R_m)$ (see \cite[Remark 8.4]{GV18}).  By restricting to $\cO_v^{\ast}$ via the local Artin reciprocity, we get an $\cO$-algebra homomorphism $\cO[\Delta_v] \to R_m$ where $\Delta_v$ is the Sylow $p$-subgroup of $(k_v^{\ast})^n$.  Define $\Delta_{Q_m} =\prod_{v \in Q_m}\Delta_v$, then $R_m$ is naturally an $\cO[\Delta_{Q_m}]$-algebra and it's clear that $R_{\cS} \cong R_m \otimes_{\cO[[\Delta_{Q_m}]]} \cO$.

Let $S_{\infty} = \cO[[X_1, \dots, X_{nr}]]$, $J_m = \langle ((1+X_i)^{p^m}-1)_{1 \leq i \leq nr} \rangle$ and $S_m = S_{\infty}/J_m$ ($m\geq 1$).  Note that $J_1 \supseteq J_2 \supseteq \dots$ is a decreasing sequence and $\cap_{i \geq 1} J_i = 0$.  Since $\Delta_{Q_m}$ is a product of $nr$ cyclic groups, each of order at least $p^m$, the ring $S_m$ is a quotient of $\cO[\Delta_{Q_m}]$.  We introduce $\bar{S}_m = S_m /p^m$ and $\bar{R}_m = R_m \otimes_{\cO[[\Delta_{Q_m}]]} \bar{S}_m$.  Let $R_{\infty} = \cO[[X_1, \dots, X_g]]$ with $g = nr +  \dim_k H^1_{\cS}(\Gamma_S, \g_k) - \dim_k H^1_{\cS^{\perp}}(\Gamma_S, \g_k^{\ast})$.  Then \lemref{relative dimension 1} implies there is a surjection $R_{\infty} \twoheadrightarrow R_m$ for every $m\geq 1$.

\subsection{Calegari-Geraghty setting}\label{Calegari-Geraghty setting}

Let $U = U_S \times U^S = (\prod_{v\in S} U_v) \times (\prod_{v\notin S} U_v)$ be a neat open compact subgroup of $\bG_f$ such that $U_v \subseteq \underline{\bG}(\cO_v)$ and each $U_v$ ($v\notin S$) is hyperspecial maximal.  Let $h$ be the image of $\cH(\bG^S, U^S)$ acting on $H^{\ast}(X_{\bG}^U, \widetilde{V}_{\lambda}(\cO))$ and let $\m$ be a non-Eisenstein maximal ideal of $h$, and we write $\T = h_{\m}$.  Assume $(\Resm)$, $(\Galm)$ and $(\Vanm)$.

Set $\cS$ to be the global deformation problem for $\bar\rho \colon \Gamma_S \to G(k)$ described in $(\Resm)$ (more precisely, it is simultaneously either ordinary or Fontaine-Laffaille for $v\in S_p$, and minimal for $v \in S\backslash S_p$).  In Section \ref{local deformation problems} we will show 
$$\sum\limits_{v\in S} (\dim_k L_v - \dim_k H^0(\Gamma_v,\g_k) ) =  [F:\Q] (\dim G -\dim B),$$
so together with the oddness condition, one has 
$$\dim_k H^1_{\cS}(\Gamma_S, \g_k) - \dim_k H^1_{\cS^{\perp}}(\Gamma_S, \g_k^{\ast}) = -\ell_0,$$
and hence $\dim S_{\infty} - \dim R_{\infty} = \ell_0$. 

The conjecture $(\Vanm)$ implies that we can choose a minimal cochain complex of $\cO$-modules $C^{\ast}$ concentrated in degrees $[q_0,q_0+\ell_0]$ such that $H^{\ast}(C^{\ast}) \cong H^{\ast}(X_{\bG}^U, \widetilde{V}_{\lambda}(\cO))_{\m}$ (see \cite[Lemma 2.3]{KT17}).  For each allowable Taylor-Wiles datum $Q_m$, it is explained in \cite[Section 13.6]{GV18} that, under the local-global compatibilities at Taylor-Wiles primes, there exists a cochain complex $C_m^{\ast}$ of finite free $\bar{S}_m$-modules such that $C_m^{\ast} \otimes_{\bar{S}_m} \cO/p^m$ is quasi-isomorphic to $C^{\ast}/p^m$ and there is a natural action of $\bar{R}_m$ on $H^{\ast}(C_m^{\ast})$ which is compatible as $\bar{S}_m$-algebras, and compatible with the $R_{\cS}$-action on $H^{\ast}(C^{\ast})$ after descending to $C_m^{\ast} \otimes_{\bar{S}_m} \cO/p^m \simeq C^{\ast}/p^m$.  We can further require $C_m^{\ast}$ to be minimal so that it is also concentrated in degrees $[q_0,q_0+\ell_0]$, and the quasi-isomorphism $C_m^{\ast} \otimes_{\bar{S}_m} \cO/p^m \simeq C^{\ast}/p^m$ is then induced from an isomorphism of chain complexes.  

To summarize the data, we have: 
\begin{enumerate}
	\item a minimal complex of $\cO$-modules $C^{\ast}$ concentrated in degrees $[q_0,q_0+\ell_0]$;
	\item an $\cO$-algebra homomorphism $R_{\cS} \to \End_{\cO}(H^{\ast}(C^{\ast}))$;
	\item a minimal complex of $\bar{S}_m$-modules $C_m^{\ast}$ concentrated in degrees $[q_0,q_0+\ell_0]$, such that $C_m^{\ast} \otimes_{\bar{S}_m} \cO/p^m \cong C^{\ast}/p^m$ for each $m\geq 1$;
	\item a commutative diagram of $\bar{S}_m$-algebra homomorphisms for each $m\geq 1$: 
	\begin{displaymath}
	\xymatrix{
		\bar{R}_m \ar[r]\ar[d]^{-\otimes_{\bar{S}_m} \cO/p^m} & \End_{\cO}(H^{\ast}(C_m^{\ast})) \ar[d]^{-\otimes_{\bar{S}_m} \cO/p^m} \\
		R_{\cS}/p^m \ar[r] & \End_{\cO/p^m}(H^{\ast}(C^{\ast}/p^m));
	}
	\end{displaymath}
	\item a surjective $\cO$-algebra homomorphism $R_{\infty} \to \bar{R}_m$ for each $m\geq 1$.
\end{enumerate}

Now by the patching argument (see \cite[Proposition 3.1]{KT17}), we can find the following data:
\begin{enumerate}
	\item[(a)] a complex of finite free $S_{\infty}$-modules $C_{\infty}^{\ast}$ concentrated in degrees $[q_0,q_0+\ell_0]$ together with an isomophism $C_{\infty}^{\ast} \otimes_{S_{\infty}} \cO \cong C^{\ast}$;
	\item[(b)] an $\cO$-algebra homomorphism $S_{\infty} \to R_{\infty}$;
	\item[(c)] a commutative diagram of $S_{\infty}$-algebra homomorphisms:
	\begin{displaymath}
	\xymatrix{
		R_{\infty} \ar[r]\ar@{->>}[d]^{-\otimes_{S_{\infty}} \cO} & \End_{S_{\infty}}(H^{\ast}(C_{\infty}^{\ast})) \ar[d]^{-\otimes_{S_{\infty}} \cO } \\
		R_{\cS} \ar[r] & \End_{\cO}(H^{\ast}(C^{\ast})).
	}
	\end{displaymath}
\end{enumerate}

\begin{rem}
	An important point in the patching argument is that $\bar{R}_m \to \End_{\cO}(H^{\ast}(C_m^{\ast}))$ factors through $\bar{R}_m / \m_{\bar{R}_m}^{c(m)}$ for a constant $c(m)$ only depending on $m$, so essentially the datum $(\bar{R}_m / \m_{\bar{R}_m}^{c(m)},C_m^{\ast})$ admits only finite choices, and hence we can select a compatible system satisfying conditions (3)-(5) and pass to the inverse limit.
\end{rem}

The difference with the Taylor-Wiles method is the appearance of the positive $\ell_0$, both as $\dim S_{\infty} - \dim R_{\infty}$ and the length of the interval $[q_0, q_0+\ell_0]$.  Note $\dim_{S_{\infty}} H^{\ast}(C_{\infty}^{\ast}) = \dim_{R_{\infty}} H^{\ast}(C_{\infty}^{\ast}) \leq \dim R_{\infty} = \dim S_{\infty}-\ell_0$ (the first equality is because $R_{\infty}/\Ann_{R_{\infty}}(H^{\ast}(C_{\infty}^{\ast}))$ acts faithfully on the finite $S_{\infty}$-module $H^{\ast}(C_{\infty}^{\ast})$, so it is finite over $S_{\infty}$).  By the commutative algebra lemma \ref{essential commutative algebra lemma} (applying to $S=S_{\infty}$ and $D^{\ast}=C^{\ast}_{\infty}$), we know $H^i(C_{\infty}^{\ast})$ is non-zero only at degree $i=q_0+\ell_0$, and 
\begin{displaymath}
\left\{\begin{array}{l}{ \depth_{S_{\infty}} \, H^{q_0+\ell_0}(C_{\infty}^{\ast}) =\dim_{S_{\infty}} H^{q_0+\ell_0}(C_{\infty}^{\ast})= \dim S_{\infty}-\ell_0;} \\ {  \pd_{S_{\infty}} H^{q_0+\ell_0}(C_{\infty}^{\ast}) = \ell_0.}\end{array}\right.
\end{displaymath}
See also \cite[Theorem 2.1.1]{Han12}, which proves above results by a different approach.

\begin{lem}\label{essential commutative algebra lemma}
	Let $S$ be a Noetherian local ring.  Let $D^{\ast}$ be a complex of finite free $S$-modules concentrated in degrees $[q_m,q_s]$.  Let $\ell=q_s-q_m$.  Suppose $H^{\ast}(D^{\ast})\neq 0$, then $\dim_S H^{\ast}(D^{\ast})\geq  \depth \, S-\ell$.  If equality holds, then $H^i(D^{\ast})$ is non-zero only at degree $i=q_s$, and we have $\depth_S \, H^{q_s}(D^{\ast})=\depth \, S-\ell$, $\pd_S H^{q_s}(D^{\ast}) = \ell$.
\end{lem}

\begin{proof}
	Let $q$ be the smallest integer that $H^q(D^{\ast})\neq 0$, and set $K^q  = D^q/\im(D^{q-1})$.  
	
	Note that $	0 \to D^{q_m} \to \dots \to D^q$ is a projective resolution of $K^q$, so $\pd_SK^q \leq q-q_m$.  On the other hand, by Ischebeck's Lemma (see \cite[(15.E) Lemma 2]{Mat80}), $\Ext^i_S(H^q(D^{\ast}), K^q) = 0$ for $i< \depth_S \, K^q - \dim_S H^q(D^{\ast})$.  In particular, since $H^q(D^{\ast})$ is a non-zero submodule of $K^q$, we must have $\depth_S \, K^q \leq \dim_S H^q(D^{\ast})$.
	
	By the Auslander–Buchsbaum formula (see \cite[Tag 090V]{Sta}), we get the desired inequality:
	\begin{displaymath}
	\depth \, S = \depth_S \, K^q+\pd_SK^q \leq \dim_SH^q(D^{\ast}) + (q-q_m) \leq \dim_S H^q(D^{\ast}) + \ell.
	\end{displaymath}  
	
	If the two inequalities are actually equalities, then the second one gives $q=q_s$, which implies $K^q =H^{q_s}(D^{\ast}) $, and the first one then gives 
	\begin{displaymath}
	\left\{\begin{array}{l}{\depth_S \, H^{q_s}(D^{\ast})=\dim_S H^{q_s}(D^{\ast}) =\depth \, S-\ell;} \\ {  \pd_S H^{q_{s}}(D^{\ast}) = \ell.}\end{array}\right.
	\end{displaymath}
\end{proof}

\begin{cor}\label{Calegari-Geraghty corollary}
	\begin{enumerate}
		\item $H^{\ast}(C_{\infty}^{\ast}) = H^{q_0+\ell_0}(C_{\infty}^{\ast})$ is free over $R_{\infty}$.
		\item There is an isomorphism $H^{q_0+\ell_0-i}(C^{\ast}) \cong \Tor_i^{S_{\infty}}(H^{q_0+\ell_0}(C_{\infty}^{\ast}),\cO)$.
		\item $H^{q_0+\ell_0}(C^{\ast})$ is free over $R_{\cS}$ and $R_{\cS} \to \T$ is an isomorphism.
	\end{enumerate}
\end{cor}

\begin{proof}
	\begin{enumerate}
		\item Since the map $S_{\infty} \to R_{\infty}$ repects the module structures of $H^{q_0+\ell_0}(C_{\infty}^{\ast})$, it sends a regular sequence in $S_{\infty}$ for $H^{q_0+\ell_0}(C_{\infty}^{\ast})$ to a regular sequence in $R_{\infty}$ for $H^{q_0+\ell_0}(C_{\infty}^{\ast})$, so 
		$$\depth_{R_{\infty}} \, H^{q_0+\ell_0}(C_{\infty}^{\ast}) \geq \depth_{S_{\infty}} \, H^{q_0+\ell_0}(C_{\infty}^{\ast}) = \depth \, R_{\infty}.$$  
		Since $H^{q_0+\ell_0}(C_{\infty}^{\ast})$ is finitely generated over the regular local ring $R_{\infty}$, it's well-known that the projective dimension of $H^{q_0+\ell_0}(C_{\infty}^{\ast})$ over $R_{\infty}$ is finite (consider the Koszul resolution for $R_{\infty} / \m_{R_{\infty}}$ or see \cite[Tag 00O7]{Sta}) and we can apply the Auslander–Buchsbaum formula $$\pd_{R_{\infty}}H^{q_0+\ell_0}(C_{\infty}^{\ast}) = \depth \, R_{\infty}- \depth_{R_{\infty}} \, H^{q_0+\ell_0}(C_{\infty}^{\ast}) \leq 0.$$ 
		Therefore $H^{q_0+\ell_0}(C_{\infty}^{\ast})$ is free over $R_{\infty}$. 
		\item The Künneth spectral sequence (see \cite[Theorem 5.6.4]{Weib94}, we use the cohomological version)
		\begin{align*}
		E_2^{p,q} = \Tor_{-p}^{S_{\infty}}(H^q (C_{\infty}^{\ast}), &\cO) \Rightarrow H^{p+q}(C_{\infty}^{\ast} \otimes_{S_{\infty}} \cO)
		\end{align*}
		collapses because $E_2^{p,q} = 0$ unless $q=q_0+\ell_0$, so we get the desired isomorphism.
		\item The above results imply that $H^{q_0+\ell_0}(C^{\ast}) \cong H^{q_0+\ell_0}(C_{\infty}^{\ast}) \otimes_{S_{\infty}} \cO$ is free over $R_{\infty} \otimes_{S_{\infty}} \cO$.  Then since the module structure on $H^{q_0+\ell_0}(C^{\ast})$ factors through $R_{\infty} \otimes_{S_{\infty}} \cO \twoheadrightarrow R_{\cS}$, the map $R_{\infty} \otimes_{S_{\infty}} \cO \to R_{\cS}$ is an isomorphism and $H^{q_0+\ell_0}(C^{\ast})$ is free over $R_{\cS}$.
	\end{enumerate}
\end{proof}

\subsection{Graded structure}\label{graded}

Let's discuss the graded structures of $\Tor_{\ast}^{S_{\infty}}(R_{\infty},\cO)$ and $\Tor_{\ast}^{S_{\infty}}(H^{q_0+\ell_0}(C_{\infty}^{\ast}),\cO)$.  A priori, these are graded modules, but in fact $\Tor_{\ast}^{S_{\infty}}(R_{\infty},\cO)$ carries additional structures: it's a graded commutative ring.  

\begin{de}
	\begin{enumerate}
		\item A graded commutative ring is a graded ring $A=\oplus_{i \geq 0} A_i$, such that the multiplication satisfies $a \cdot b= (-1)^{mn}b \cdot a$ for $a\in A_m$ and $b\in A_n$.
		\item Let $A=\oplus_{i \geq 0} A_i$ be a graded commutative ring.  A graded $A$-module is an $A$-module $M$ equipped with a graded structure $M= \oplus_{i \geq 0} M_i$ such that the scalar multiplication sends $A_m\times M_n$ to $M_{m+n}$.
	\end{enumerate}
\end{de}

\begin{de}
	\begin{enumerate}
		\item A differential graded ring is a graded commutative ring $A=\oplus_{i \geq 0} A_i$ equipped with a differential $d \colon A \to A$ ({\it i.e.}, a group homomorphism for the additive structure of $A$) satisfying  
		\begin{enumerate}
			\item $d$ sends $A_i$ to $A_{i-1}$;
			\item $d\circ d=0$;
			\item $d(a \cdot b) = (da) \cdot  b+(-1)^m a \cdot  (db)$ for $a\in A_m$ and $b\in A_n$.
		\end{enumerate}
		\item Let $A$ be a differential graded ring with differential $d$.  A differential graded $A$-module is a graded $A$-module $M= \oplus_{i \geq 0} M_i$ equipped with a differential $d_M\colon M \to M$ ({\it i.e.}, a group homomorphism for the additive structure of $M$) satisfying 
		\begin{enumerate}
			\item $d_M$ sends $M_i$ to $M_{i-1}$;
			\item $d_M\circ d_M=0$;
			\item $d_M(a \cdot x) = (da) \cdot  x+(-1)^m a \cdot  (d_M x)$ for $a\in A_m$ and $x\in M_n$.
		\end{enumerate}
	\end{enumerate}
\end{de}

A differential graded ring or module has a natural chain complex structure, and we write $H_{\ast}(-)$ for the homology.  Note that if $M$ is a differential graded $A$-module, then $H_{\ast}(M)$ is naturally a graded $H_{\ast}(A)$-module. 

When $A$ is a ring and $B_1,B_2$ are $A$-algebras, the Tor-algebra $\Tor^A_{\ast}(B_1,B_2)$ can be calculated as $\pi_{\ast}(B_1 \otimes_A c(B_2))$ where $c(B_2)$ is a cofibrant replacement of $B_2$ in ${}_A \backslash \sCR$ (see Section \ref{simplicial commutative ring} and \cite[Section 7.11]{Gil13}).  In fact, $\Tor^A_{\ast}(B_1,B_2)$ is a strictly graded commutative $A$-algebra equipped with divided powers (see \cite[Section 8.5]{Gil13}).  

In our situation, the Koszul resolution of the $S_{\infty}$-algebra $\cO$ is a differential graded ring, and by \cite[Theorem 11.8]{BMR13}, one can calculate the Tor-algebra $\Tor_{\ast}^{S_{\infty}}(R_{\infty},\cO)$ using this differential graded resolution instead of the simplicial resolution.  

\begin{lem}
	$\Tor_{\ast}^{S_{\infty}}(H^{q_0+\ell_0}(C_{\infty}^{\ast}),\cO)$ is naturally a graded module over the graded commutative $\cO$-algebra $\Tor_{\ast}^{S_{\infty}}(R_{\infty},\cO)$, freely generated by $\Tor_0^{S_{\infty}}(H^{q_0+\ell_0}(C_{\infty}^{\ast}),\cO)$.
\end{lem}

\begin{proof}
	Let $E\cong (S_{\infty})^{nr}$ and let $\{ \e_1, \dots, \e_{nr} \}$ be the canonical basis.  Since $(X_1, \dots, X_{nr})$ is a regular sequence in $S_{\infty}$ and $S_{\infty} / (X_1, \dots, X_{nr}) \cong \cO$,  the Koszul complex $K_{\ast}(s)$ associated to the $S_{\infty}$-linear map $s\colon E \to S_{\infty}$ which sends $\e_i$ to $X_i$ is a free resolution of $\cO$.  Recall that 
	\begin{displaymath}
	K_{\ast}(s) \colon 0 \rightarrow \bigwedge^{nr} E \xrightarrow{d_{nr}} \cdots \xrightarrow{d_2} \bigwedge^1 E \xrightarrow{d_1} \bigwedge^0 E \cong S_{\infty} \rightarrow 0,
	\end{displaymath}
	where $d_k(a_1 \wedge \dots \wedge a_k) = \sum_{i=1}^k (-1)^{i-1} s(a_i) a_1 \wedge \cdots \wedge \hat{a}_i \wedge \cdots \wedge a_k.$
	
	Note that $K_{\ast}(s)$ is naturally a differential graded ring with the multiplication defined by 
	\begin{displaymath}
	(a_1 \wedge \cdots \wedge a_i) \cdot (b_1 \wedge \cdots \wedge b_j) = a_1 \wedge \cdots \wedge a_i \wedge b_1 \wedge \cdots \wedge b_j.
	\end{displaymath} 
	Then together with the $R_{\infty}$-module structure on $H^{q_0+\ell_0}(C_{\infty}^{\ast})$, $K_{\ast}(s) \otimes_{S_{\infty}} H^{q_0+\ell_0}(C_{\infty}^{\ast})$ is naturally a differential graded $K_{\ast}(s) \otimes_{S_{\infty}} R_{\infty}$-module.  By the foregoing comment, $\Tor_{\ast}^{S_{\infty}}(H^{q_0+\ell_0}(C_{\infty}^{\ast}),\cO) \cong H_{\ast}(K_{\ast}(s) \otimes_{S_{\infty}} H^{q_0+\ell_0}(C_{\infty}^{\ast}))$ is a graded module over the graded commutative ring $\Tor_{\ast}^{S_{\infty}}(R_{\infty},\cO) \cong H_{\ast}(K_{\ast}(s) \otimes_{S_{\infty}} R_{\infty})$.  Moreover, $\Tor_{\ast}^{S_{\infty}}(H^{q_0+\ell_0}(C_{\infty}^{\ast}),\cO)$ is freely generated by $\Tor_0^{S_{\infty}}(H^{q_0+\ell_0}(C_{\infty}^{\ast}),\cO)$ because $H^{q_0+\ell_0}(C_{\infty}^{\ast})$ is free over $R_{\infty}$.	
\end{proof}

Note that $H^{\ast}(X_{\bG}^U, \widetilde{V}_{\lambda}(\cO))_{\m} \cong H^{\ast}(C^{\ast})$ is equipped with a graded structure (note the switch of indexes $i \mapsto q_0 +\ell_0-i$) via the isomorphism $H^{q_0 +\ell_0-i}(C^{\ast}) \cong \Tor_i^{S_{\infty}}(H^{q_0+\ell_0}(C_{\infty}^{\ast}),\cO)$.  The following corollary is straightforward:

\begin{cor}
	$H^{\ast}(X_{\bG}^U, \widetilde{V}_{\lambda}(\cO))_{\m}$ is a graded $\Tor_{\ast}^{S_{\infty}}(R_{\infty},\cO)$-module, freely generated by $ H^{q_0+\ell_0}(X_{\bG}^U, \widetilde{V}_{\lambda}(\cO))_{\m}$.
\end{cor}

\begin{rem}
	As mentioned in the introduction, an unsatisfactory aspect is that $\Tor_{\ast}^{S_{\infty}}(R_{\infty},\cO)$ depends on various non-canonical choices.  Note the isomorphism 
	$$\Tor_{\ast}^{S_\infty}(R_\infty, \cO) \cong \pi_{\ast} (R_{\infty} \uotimes_{S_{\infty}} \cO)$$ 
	as graded commutative $\cO$-algebras, where $R_{\infty} \uotimes_{S_{\infty}} \cO$ is a simplicial ring which represents the derived tensor product $R_{\infty} \lotimes_{S_{\infty}} \cO$ (see Section \ref{simplicial commutative ring}).  The insight of Galatius and Venkatesh is that one can extend deformations to simplicial rings and reinterpret $R_{\infty} \uotimes_{S_{\infty}} \cO$ as a derived representing ring, thus canonically.  In the following chapters we will discuss the derived deformation functors and derived deformation rings, which are the principal subjects of this paper.
\end{rem}

\section{Simplicial backgrounds}

The derived deformation functors are more or less functors from simplicial commutative $\cO$-algebras to simplicial sets.  In this section we will prepare the necessary foundations to study these functors.  In Section 2.1 we will recall some basic facts on simplicial model categories, and an important objective is to understand the structure of $\sCROk$.  In Section 2.2 we will focus on the pro-representabilty of functors from the Artinian subcategory $\sArt$ of $\sCROk$ to simplicial sets.

\subsection{Simplicial model categories}

\subsubsection{Simplicial sets}

We denote by $\bm{\Delta}$ the cosimplicial indexing category: the objects are totally ordered sets $[n] = \{0,1,\dots,n\}$ for $n\geq 0$, and the morphisms are order-preserving functions between these sets.  Let $d^i \colon [n-1] \to [n]$ $(0 \leq i \leq n)$ and $s^j \colon [n+1]\to [n]$ $(0\leq j \leq n)$ be the morphisms defined by 
\begin{displaymath}
d^i (\{0,1,\dots,n-1\}) = \{0,1,\dots,i-1,i+1,\dots,n\},
\end{displaymath}
and
\begin{displaymath}
s^j (\{0,1,\dots,n+1\}) = \{0,1,\dots,j,j,\dots,n\}.
\end{displaymath}

\begin{de}
	For a category $\cC$, we define $s\cC$ to be the category of functors $\bm{\Delta}^{\op} \to \cC$.
\end{de}

In fact, an object $X \in s\cC$ can be regarded as a sequence of $X_n \in \cC$ for $n\geq 0$ ($X_n$ being the image of $[n]$) together with morphisms $d_i \colon X_n \to X_{n-1}$ $(0 \leq i \leq n)$ and $s_j \colon X_n\to X_{n+1}$ ($0 \leq j \leq n$) satisfying the relations
\begin{displaymath}
\left\{\begin{array}{l}{  d_j d_i = d_i d_{j+1} \qquad \text{ if } i \leq j }\\ { s_js_i =  s_i s_{j-1} \qquad\text{ if } i \leq j-1  } \\ { d_j s_i=s_i d_{j-1}  \qquad\text{ if } i \leq j-2  }  \\ { d_js_{j-1}=d_j s_j=\id  } \\ {d_j s_i=s_{i-1} d_j \qquad\text{ if } i\geq j+1.  } \end{array}\right.
\end{displaymath}

We call $\sSet$ the category of simplicial sets, $s\Gp$ the category of simplicial groups...

\begin{expl}
	\begin{enumerate}
		\item $\Delta^n = \Hom_{\bm{\Delta}}(- , [n])$ ($n \geq 0$) is a simplicial set, we call it the standard $n$-simplex.
		\item We denote by $\partial\Delta^n$ the smallest sub-simplicial set of $\Delta^n$ which contains $d_i(\id_{[n]})$, $0\leq i\leq n$.  We call $\partial\Delta^n$ the boundary of $\Delta^n$.  Explicitly, $\partial\Delta^n_k$ is the set of non-surjective order-preserving functions $\{0,1, \dots, k\} \to \{0,1, \dots, n\}$.
		\item Let $n\geq 1$ and $0\leq m \leq n$.  We denote by $\Lambda^n_m$ the smallest sub-simplicial set of $\Delta^n$ which contains $d_i(\id_{[n]})$ for $0\leq i\leq n$ and $i\neq m$.  We call $\Lambda^n_m$ the $m$-th horn of $\Delta^n$.  Explicitly, $(\Lambda^{n}_m)_k$ is the set of order-preserving functions $\{0,1, \dots, k\} \to \{0,1, \dots, n\}$ such that the image doesn't contain $\{0,1, \dots,m-1,m+1,\dots, n\}$.
	\end{enumerate}
\end{expl}

\begin{de}
	\begin{enumerate}
		\item A morphism of $\sSet$ is a cofibration if it is injective in every simplicial degree.
		\item Let $X$ and $Y$ be simplicial sets.  A morphism $p\colon X\to Y$ is a fibration if for every $n\geq 1$, $0\leq k \leq n$ and solid arrow commutative diagram as follows:
		\begin{displaymath}
		\xymatrix{
			\Lambda^n_k \ar[r]\ar[d]^{i} & X \ar[d]^{p} \\
			\Delta^n \ar[r]\ar@{-->}[ru] & Y,
		}
		\end{displaymath}
		where $i\colon \Lambda^n_k \hookrightarrow \Delta^n$ is the natural inclusion, there is a dotted arrow making the diagram commute.  We say a simplicial set $X$ is fibrant (or Kan), if $X \to \ast$ is a fibration (here $\ast$ refers to $\Delta^0$, which is the terminal object of $\sSet$).
	\end{enumerate}
\end{de}

A morphism $\Lambda^n_k \to X$ can be regarded as an $n$-tuple $(z_0, \dots, \hat{z}_k, \dots, z_n)$ of $z_i \in X_{n-1}$ such that $d_{j-1}z_i = d_i z_j$ for $i <j$.  Thus $p\colon X\to Y$ is a fibration if and only if for every $n\geq 1$ and $n$-tuple $(z_0, \dots, \hat{z}_k, \dots, z_n)$ as above satisfying $p(z_i) = d_i y$ for some $y\in Y_n$, there exists $x\in X_n$ such that $p(x) = y$ and $z_i = d_i x$.  

\begin{lem}\label{simplicial group fibration}
	Every simplicial group is fibrant as a simplicial set, and every morphism of simplicial groups $f\colon G\to H$ which induces surjective $G_n \to H_n$ for every $n\geq 1$ is a fibration as a morphism of simplicial sets. 
\end{lem}

\begin{proof}
	For the first statement, see \cite[Lemma \RNum{1}.3.4]{GJ09}.  For the second statement, it suffices to show that for every $n\geq 1$ and $n$-tuple $(z_0, \dots, \hat{z}_k, \dots, z_n)$ of elements in $G_{n-1}$ and $y\in H_n$ such that $d_{j-1}z_i = d_i z_j$ for $i <j$ and $f(z_i) = d_i y$, there exists $x\in G_n$ such that $f(x)=y$ and $d_i x=z_i $.  Since $G_n\to H_n$ is surjective, there exists a pre-image $x^{\prime}$ of $y$, by considering $(d_i x^{\prime})^{-1}\cdot z_i $, it reduces to show $\ker(f)$ is fibrant, which follows from the first statement.
\end{proof}

Let $\bm{\Delta} X$ be the category of simplices of $X$ (see \cite[Definition 15.1.16]{Hir09}): the objects are natural transformations $\Delta^n \to X$, and the morphisms from $\Delta^n \to X$ to $\Delta^m \to X$ consist of natural transformations $\Delta^n \to \Delta^m$ which respect the natural transformations to $X$.  By Yoneda's lemma, the objects of $\bm{\Delta} X$ can be identified with $\bigsqcup\limits_{n \geq 0} X_n$, and the morphisms from $x \in X_n$ to $y\in X_m$ can be identified with morphisms $[n] \to [m]$ of $\bm{\Delta}$ such that the induced map $X_m\to X_n$ sends $y$ to $x$.

We have the following well-known lemma:

\begin{lem}\label{adjoint}
	Suppose $\cC$ is a category admitting colimits; let $F\colon \bm{\Delta} \to \cC$ be a covariant functor.  Let $F_{\ast}\colon \cC \to \sSet$ be the functor which sends $A\in \cC$ to the simplicial set $X=(X_n)_{n\geq 0}$ given by $X_n=\Hom_{\cC}(F([n]), A)$ at $n$-th simplicial degree, and let $F^{\ast}\colon \sSet \to \cC$ be the functor which sends $X \in \sSet$ to $\varinjlim\limits_{(n,\sigma) \in \bm{\Delta}  X} F(\sigma)$.  Then $F^{\ast}$ is left adjoint to $F_{\ast}$.
\end{lem}

\begin{proof}
	It's clear that $F_{\ast}$ is well-defined, and $F^{\ast}$ is well-defined since every simplicial set morphism $f \colon X \to Y$ induces a functor $\bm{\Delta}  X \to\bm{\Delta}  Y$.  For $X\in \sSet$ and $A\in \cC$, we have
	\begin{align*}
	\Hom_{\cC}(F^{\ast}(X),A) &\cong \varprojlim\limits_{(\Delta^n \to X) \in (\bm{\Delta}  X)^{\op} } \Hom_{\cC}(F([n]), A) \\ &\cong \varprojlim\limits_{(\Delta^n \to X) \in (\bm{\Delta}  X)^{\op} } \Hom_{\sSet}(\Delta^n, F_{\ast}(A))   \\ &\cong \Hom_{\sSet}(\varinjlim\limits_{(\Delta^n \to X) \in \bm{\Delta}  X} \Delta^n, F_{\ast}(A))  \\ &\cong  \Hom_{\sSet}(X, F_{\ast}(A)),
	\end{align*}
	where the last equation follows from \cite[Proposition 15.1.20]{Hir09}.  So $F^{\ast}$ is left adjoint to $F_{\ast}$.
\end{proof}

\begin{expl}
	Consider the functor $\bm{\Delta} \to \Top$, which sends $[n]$ to $|\Delta^n|$, where $|\Delta^n| = \{ (t_0,\dots,t_n)\in \R^{n+1} \mid \sum\limits_{i=0}^n t_i=1, t_i\geq 0  \}$ is the topological standard $n$-simplex, and sends morphisms of $\bm{\Delta}$ to corresponding linear maps.  The associated left adjoint sends $X\in \sSet$ to $|X| =\varinjlim\limits_{(\Delta^n \to X) \in \bm{\Delta} X}|\Delta^n|$, and the associated right adjoint is the usual singular complex functor.  We call $|X|$ the geometric realization of $X$.
\end{expl}

\begin{de}
	A morphism of simplicial sets $f\colon X \to Y$ is a weak equivalence, if the induced map $|f|\colon |X|\to |Y|$ is a topological weak equivalence.
\end{de}

\begin{de}
	Let $X$ be a simplicial set and let $v\colon \ast \to X$ be a vertex of $X$.  We also use $v$ to denote the corresponding point of the geometric realization $ |X|$.  Then for $n\geq 1$ the $n$-th homotopy group of $(X,v)$ is defined by $\pi_n(X, v)=\pi_n(|X|, v)$.  We also define $\pi_0(X) = \pi_0(|X|)$.
\end{de}

For fibrant $X$, the group structures on $\pi_n(X, v)$ for $n\geq 1$ can be defined combinatorially without refering to the geometric realization (see \cite[Section \RNum{1}.7]{GJ09}).  In particular, this is the case when $X\in s\Gp$.  Henceforth, when $X$ is a simplicial group with unit $e$ and $n\geq 1$, we will abbreviate $\pi_n(X,e)$ by $\pi_n(X)$.  Since changing the vertex $v$ induces group isomorphisms of homotopy groups $\pi_n(X,v)$ natural in $X$, a morphism $f\colon X \to Y$ of simplcial groups is a weak equivalence in $\sSet$ if and only if $\pi_n(f)\colon \pi_n(X) \to \pi_n(Y)$ is an isomorphism for all $n$. 

The reason for introducing cofibrations, fibrations and weak equivalences of $\sSet$ is that with these strucutures, the category $\sSet$ becomes a model category.

\subsubsection{Model categories}

\begin{de}
	A category $\cC$ is a model category, if it is equipped with three classes of morphisms: cofibrations, fibrations and weak equivalences (we say a cofibration or fibration is trivial if it is also a weak equivalence), such that the following axioms hold:  
	\begin{enumerate}
		\item[\textbf{CM1}:] $\cC$ is complete and cocomplete.
		\item[\textbf{CM2}:] Given composable morphisms $f,g$ of $\cC$, if any two of $f,g$ and $fg$ are weak equivalences, then so is the third.
		\item[\textbf{CM3}:] If $f$ is a retract of $g$ and $g$ is a cofibration, fibration or weak equivalence, then so is $f$.
		\item[\textbf{CM4}:] If either $i$ is a trivial cofibration and $p$ is a fibration, or $i$ is a cofibration and $p$ is a trivial fibration, then $i$ has the left lifting property with respect to $p$ (tautologically $p$ has the right lifting property with respect to $i$), {\it i.e.}, for every solid arrow commutative diagram
		\begin{displaymath}
		\xymatrix{
			A \ar[r]\ar[d]^i & X \ar[d]^p \\
			B \ar[r]\ar@{-->}[ru] & Y,
		}
		\end{displaymath}
		there exists a dotted arrow making the diagram commutative.
		\item[\textbf{CM5}:] Any morphism $f\colon X\to Y$ can be factored in two ways:
		\begin{enumerate}
			\item $f=pi$, where $p$ is a fibration and $i$ is a trivial cofibration.
			\item $f=q j$, where $q$ is a trivial fibration and $j$ is a  cofibration.
		\end{enumerate}
	\end{enumerate}
\end{de}

\begin{rem}
	\begin{enumerate}
		\item It's customary to write $\hookrightarrow$ for a cofibration, $\twoheadrightarrow$ for a fibration, and $\iso$ for a weak equivalence.
		\item The axiom \textbf{CM1} implies that $\cC$ has an initial object $\emptyset$ and a terminal object $\ast$.  We say an object $A \in \cC$ is cofibrant if $ \emptyset \hookrightarrow A$, and fibrant if $A \twoheadrightarrow \ast$.
		\item We say $B$ is a cofibrant replacement of $A$ if $\emptyset \hookrightarrow B \tilde{\twoheadrightarrow} A$.  We say $B$ is a fibrant replacement of $A$ if $ A \tilde{\hookrightarrow} B \twoheadrightarrow \ast$.
		\item If $\cC$ is a model category, then the opposite category $\cC^{\op}$ also carries a model category structure: a morphism of $\cC^{\op}$ is a cofibration, fibration or weak equivalence if and only if its dual is a fibration, cofibration or weak equivalence of $\cC$ respectively.  So if we can prove some statement under axioms of model category, then the dual statement is also true.
		\item It follows from the axioms \textbf{CM3}, \textbf{CM4} and \textbf{CM5} that a morphism is a cofibration if and only if it has the left lifting property with respect to all trivial fibrations, and a morphism is a trivial cofibration if and only if it has the left lifting property with respect to all fibrations.  Similarly, a morphism is a fibration if and only if it has the right lifting property with respect to all trivial cofibrations, and a morphism is a trivial fibration if and only if it has the right lifting property with respect to all cofibrations.
	\end{enumerate}
\end{rem}

Let's review the theory of cofibrantly generated model categories for the first infinite cardinal, which is sufficient for our purpose.  See \cite[Chapters 10 and 11]{Hir09} for transfinite generalizations.

\begin{de}\label{small}
	Let $\cC$ be a category.
	\begin{enumerate}
		\item Let $\mathfrak{U}$ be a class of morphisms of $\cC$.  We say an object $X\in \cC$ is small relative to $\mathfrak{U}$ if for every (countable) sequence
		\begin{displaymath}
		Y_0 \to Y_1 \to \dots \to Y_i \to \dots
		\end{displaymath}
		where each $Y_i \to Y_{i+1}$ belongs to $\mathfrak{U}$, the natural map $\varinjlim\limits_i \Hom_{\cC}(X, Y_i) \to  \Hom_{\cC}(X, \varinjlim\limits_i Y_i)$ is an isomorphism.
		\item Let $I$ be a set of morphisms of $\cC$.  We say $I$ permits the small object argument if the sources of morphisms of $I$ are small relative to the class of morphisms consisting of pushouts of coproducts of $I$.
	\end{enumerate}
\end{de}

\begin{de}
	A model category $\cC$ is cofibrantly generated, if it satisfies the following two conditions: 
	\begin{enumerate}
		\item There is a set of morphisms $I$, such that $I$ permits the small object argument, and a morphism is a trivial fibration if and only if it has the right lifting property with respect to all elements of $I$.  We call such $I$ a set of generating cofibrations.
		\item There is a set of morphisms $J$, such that $J$ permits the small object argument, and a morphism is a fibration if and only if it has the right lifting property with respect to all elements of $J$.  We call such $J$ a set of generating trivial cofibrations.
	\end{enumerate}
\end{de}

The small object argument of Quillen implies that the factorizations in $\textbf{CM5}$ can be chosen functorial.  We say a morphism $f\colon X_0 \to X$ is an $\N$-composition of morphisms in some class $\mathfrak{U}$ if there exists $X_0 \to X_1 \to \dots \to X_i \to \dots$ such that each $X_i \to X_{i+1}$ belongs to $\mathfrak{U}$ and $f$ coincides with $X_0 \to \varinjlim\limits_i X_i$.

\begin{lem}
	Let $\cC$ be a cofibrantly generated model category with a set of generating cofibrations $I$ and a set of generating trivial cofibrations $J$.
	\begin{enumerate}
		\item There is a functorial factorization of every morphism of $\cC$ into a cofibration followed by a trivial fibration, such that the cofibration is an $\N$-composition of pushouts of coproducts of $I$.
		\item There is a functorial factorization of every morphism of $\cC$ into a trivial cofibration followed by a fibration, such that the trivial cofibration is an $\N$-composition of pushouts of coproducts of $J$.
	\end{enumerate}
\end{lem}

\begin{proof}
	See \cite[Corollary 11.2.6]{Hir09}.
\end{proof}

\begin{cor}\label{cofibration}
	Let notations be as above.  Then a morphism of $\cC$ is a cofibration if and only if it is a retract of an $\N$-composition of pushouts of coproducts of $I$.
\end{cor}

\begin{proof}
	See \cite[Corollary 10.5.23]{Hir09}.
\end{proof}

\begin{expl}\label{chain complex}
	$\Ch_{\geq 0}(R)$, the category of chain complexes of $R$-modules concentrated in non-negative degrees for a commutative ring $R$, is a cofibrantly generated model categroy.  The cofibrations, fibrations and weak equivalences are characterized as follows:
	\begin{enumerate}
		\item $f\colon C_{\ast}\to D_{\ast}$ is a cofibration if $C_n\to D_n$ is injective with projective cokernel for $n\geq 0$.
		\item $f\colon C_{\ast}\to D_{\ast}$ is a fibration if $C_n\to D_n$ is surjective for $n\geq 1$.
		\item $f\colon C_{\ast}\to D_{\ast}$ is a weak equivalence if $H_{\ast}f$ is an isomorphism.
	\end{enumerate}
	Thus every $C_{\ast} \in \Ch_{\geq 0}(R)$ is fibrant, and taking cofibrant replacement means exactly taking projective resolution in the sense of homological algebra.
	
	For $n\geq 0$, let $R[n]$ be the chain complex with $R$ on $n$-th degree and with $0$ elsewhere, and let $R\left\langle n+1\right\rangle $ be the chain complex
	\begin{displaymath}
	\dots \to 0 \to \stackrel{n+1}{R}  =  \stackrel{n}{R} \to 0 \to \dots
	\end{displaymath}
	Then the generating cofibrations may be taken to be $0 \to R[0]$ together with natural inclusions $R[n] \to R\left\langle n+1\right\rangle$, and the generating trivial cofibrations may be taken to be $0 \to R\left\langle n+1\right\rangle$.
\end{expl}

For a model category $\cC$ there are (left or right) homotopy relations for morphisms $f,g \colon X\to Y$ of $\cC$.  For our purpose we will only focus on the case where $X$ is cofibrant and $Y$ is fibrant, and in this case the left and right homotopy relations coincide and define an equivalence relation (see \cite[Section 7.3 and 7.4]{Hir09} for details).

\begin{de}
	Let $\cC$ be a model category and $f\colon A\to B$ be a given morphism of $\cC$.  We define the over and under category ${}_A \backslash\cC /_B$, such that the objects are arrows $A\to X \to B$ with composition $f$, and the morphisms from $A\to X\to B$ to $A\to Y \to B$ are the morphisms $X\to Y$ which respect the morphisms from $A$ and to $B$.
\end{de}

\begin{lem}\label{over and under category}
	The category ${}_A \backslash\cC /_B$ is a model categories, with cofibrations, fibrations and weak equivalences being those of $\cC$.
\end{lem}

\begin{proof}
	It suffices to check the axioms \textbf{CM1} to \textbf{CM5} hold, and they follow directly from the corresponding properties for $\cC$.
\end{proof}

We can regard ${}_A \backslash\cC /_B$ as a subcategory of $\cC$.  Then if two morphisms $f,g\colon X\to Y$ are (left or right) homotopic in ${}_A \backslash\cC /_B$, they are (left or right) homotopic in $\cC$ (see \cite[Proposition 7.6.8]{Hir09}).

\subsubsection{Homotopy categories and derived functors}

For a model category $\cC$, the localization with respect to weak equivalences exists.  More precisely, there is an associated homotopy category $\Ho(\cC)$ with a functor $\gamma\colon \cC \to \Ho(\cC)$, such that $\gamma(f)$ is an isomorphism if and only if $f$ is a weak equivalence, and if $F\colon \cC \to \cD$ is a functor which sends weak equivalences to isomorphisms, then there is a unique functor $F_{\ast}\colon \Ho(\cC) \to \cD$ such that $F_{\ast} \circ \gamma = F$.  We remark that $\Ho(\cC)$ has same objects as $\cC$ and the functor $\gamma\colon \cC \to \Ho(\cC)$ is identity on objects.  The morphisms of $\Ho(\cC)$ satisfy 
\begin{displaymath}
\Hom_{\Ho(\cC)}(A, X) \cong \Hom_{\cC}(B, Y)/(\text{homotopy}),
\end{displaymath} 
where $B$ is any cofibrant replacement of $A$ and $Y$ is any fibrant replacement of $X$.  See \cite[Section 8.3]{Hir09} for details. 

\begin{lem}
	Let $\cC$ be a model category and $\cA$ be any category.  We fix simultaneously a cofibrant replacement $X^{\prime}$ for every $X \in \cC$.  Suppose $F\colon \cC \to \cA$ is a functor which sends trivial cofibrations between cofibrant objects to isomorphisms.  Then there is a well-defined functor
	\begin{displaymath}
	\mathbf{L}F \colon \Ho(\cC) \to \cA
	\end{displaymath}
	which sends $X$ to $F(X^{\prime})$.  We say that $\mathbf{L}F$ is the total left derived functor for $F$.
\end{lem}

\begin{proof}
	See \cite[Lemma 7.7.1]{Hir09} and \cite[Lemma \RNum{2}.7.3]{GJ09}.
\end{proof}

Note that the total left derived functor depends on the system of cofibrant replacements up to natural isomorphism.

We may dually define the total right derived functor $\mathbf{R}$ for a functor which sends trivial fibrations between fibrant objects to isomorphisms.

\begin{de}
	Let $\cC,\cD$ be two model categories and let $F\colon \cC \rightleftarrows \cD \colon G$ be a pair of adjoint functors.  We say $(F,G)$ is a Quillen pair if one of the following equivalent conditions holds:
	\begin{enumerate}
		\item $F$ preserves cofibrations and trivial cofibrations.
		\item $G$ preserves fibrations and trivial fibrations.
	\end{enumerate}
	In this case we say $F$ is a left Quillen functor and $G$ is a right Quillen functor.
\end{de}

\begin{thm}\label{Quillen}
	Let $\cC,\cD$ be two model categories and let $F\colon \cC \rightleftarrows \cD \colon G$ be a pair of adjoint functors.  Suppose $(F,G)$ is a Quillen pair.  Then $\mathbf{L}F \colon \Ho(\cC) \to \Ho(\cD)$ and $\mathbf{R}G \colon \Ho(\cD) \to \Ho(\cC)$ exist, and $\mathbf{R}G$ is right adjoint to $\mathbf{L}F$.  If furthermore for cofibrant $A \in \cC$ and fibrant $X\in \cD$, the map $A \to GX$ is a weak equivalence if and only if the adjoint map $FA \to X$ is a weak equivalence, then $\mathbf{L}F$ and $\mathbf{R}G$ induce an adjoint equivalence of categories $\Ho(\cC) \cong \Ho(\cD)$.
\end{thm}

\begin{proof}
	See \cite[Theorem 8.5.18 and Theorem 8.5.23]{Hir09}.
\end{proof}

\begin{expl}
	Let $\cC$ be a model category and let $I$ be a small category.
	
	Suppose that there exists a model category structure on $\cC^I$ such that a morphism $A \to B$ is a fibration or weak equivalence if and only if every $A(i)\to B(i)$ ($i\in I$) is a fibration or weak equivalence in $\cC$ (this holds when $\cC$ is cofibrantly generated, see \cite[Theorem 11.6.1]{Hir09}); we call it the projective model structure and denote it by $\cC^I_{\proj}$.  Then the constant functor $\Delta\colon \cC \to \cC^I_{\proj}$ preserves fibrations and weak equivalences, so the left adjoint functor $\varinjlim \colon \cC^I_{\proj} \to \cC$ is left Quillen, and the total left derived functor $\mathbf{L}\varinjlim$ exists.  For convenience, we denote the colimit of some cofibrant replacement by $\hocolim$ (it is defined up to weak equivalence) and call it the homotopy colimit.
	
	Dually, supppose that there exists a model category structure on $\cC^I$ such that a morphism $A \to B$ is a cofibration or weak equivalence if and only if every $A(i)\to B(i)$ ($i\in I$) is a cofibration or weak equivalence in $\cC$ (this holds when $\cC$ is combinatorial, see \cite[Proposition A.2.8.2]{Lur09}); we call it the injective model structure and denote it by $\cC^I_{\inj}$.  Then the constant functor $\Delta\colon \cC \to \cC^I_{\inj}$ preserves cofibrations and weak equivalences, so the right adjoint $\varprojlim \colon \cC^I_{\inj} \to \cC$ is right Quillen, and the total right derived functor $\mathbf{R}\varprojlim$ exists.  For convenience, we denote the limit of some fibrant replacement by $\holim$ (it is defined up to weak equivalence) and call it the homotopy limit.
	
	We will primarily work with certain specific types of $I$, where it has a Reedy category structure.  Then there is a Reedy model cateogory stucture on $\cC^I$ which can be described explicitly (see \cite[Theorem 15.3.4]{Hir09}).  Moreover, the homotopy limit (resp. homotopy colimit) can be computed via Reedy fibrant replacements (resp. Reedy cofibrant replacements).  See \cite[Proposition 15.10.10 and 15.10.12]{Hir09} for the cases of homotopy pulllbacks (homotopy pushouts) and homotopy (co)limits indexed by $\N$.  
	
	When $I$ is represented by the diagram $\bullet \rightarrow \bullet \leftarrow \bullet $, we also write $A_1 \times_{A_0}^{h} A_2$ for $\holim \, (A_1 \to A_0 \leftarrow A_2)$.  We say a diagram
	\begin{displaymath}
	\xymatrix{
		A \ar[r]\ar[d]& A_1 \ar[d] \\
		A_2 \ar[r] & A_0 
	}
	\end{displaymath}
	is a homotopy pullback square if the natural map $A \to A_1 \times_{A_0}^{h} A_2$ is a weak equivalence.  
\end{expl}

\begin{lem}\label{holim commute}
	Let $\cC$ and $\cD$ be two model categories.  Let $F\colon \cC \to \cD$ and $G\colon \cD \to \cC$ be two functors and suppose that $(F,G)$ is a Quillen pair.  Let $I$ be a small category.
	\begin{enumerate}
		\item If $\cC^I_{\proj},\cD^I_{\proj}$ exist, then $\mathbf{L}F \circ \mathbf{L}\varinjlim$ is naturally isomorphic to $\mathbf{L}\varinjlim \circ \mathbf{L}F$.
		\item If $\cC^I_{\inj},\cD^I_{\inj}$ exist, then $\mathbf{R}G \circ \mathbf{R}\varprojlim$ is naturally isomorphic to $\mathbf{R}\varprojlim \circ \mathbf{R}G$.
	\end{enumerate}
\end{lem}

\begin{proof}
	We prove the second part, and the proof for the first part is similar.
	
	Let $F^I\colon \cC^I_{\inj} \to \cD^I_{\inj}$ and $G^I\colon \cD^I_{\inj} \to \cC^I_{\inj}$ be the degreewise extensions of $F$ and $G$.  Then it's easy to see that $F^I$ is left adjoint to $G^I$, and $F^I$ preserves cofibrations and weak equivalences, so $(F^I,G^I)$ is a Quillen pair.  Since $\Delta\colon \cC \to \cC^I_{\inj}$ and $F\colon \cC\to \cD$ preserve cofibrant objects, the following diagram commutes up to natural isomorphism:
	\begin{displaymath}
	\xymatrix{
		\Ho(\cC^I_{\inj}) \ar[r]^{\mathbf{L}F^I}& \Ho(\cD^I_{\inj})  \\
		\Ho(\cC) \ar[r]^{\mathbf{L}F}\ar[u]^{\mathbf{L}\Delta} & \Ho(\cD) \ar[u]^{\mathbf{L}\Delta}.
	}
	\end{displaymath}
	Therefore the adjoint diagram
	\begin{displaymath}
	\xymatrix{
		\Ho(\cC^I_{\inj}) \ar[d]^{\mathbf{R}\varprojlim} & \Ho(\cD^I_{\inj}) \ar[l]_{\mathbf{R}G^I}\ar[d]^{\mathbf{R}\varprojlim} \\
		\Ho(\cC)  & \Ho(\cD) \ar[l]_{\mathbf{R}G}
	}
	\end{displaymath}
	commutes up to natural isomorphism.
\end{proof}

\subsubsection{Simplicial model categories}

\begin{de}
	A category $\cC$ is a simplicial category if there is a mapping space functor
	\begin{displaymath}
	\sHom_{\cC}(-, -)\colon \cC^{\op} \times \cC \to \sSet, 
	\end{displaymath}
	with the following properties:
	\begin{enumerate}
		\item $\sHom_{\cC}(A,B)_0 = \Hom_{\cC}(A,B)$.
		\item The functor $\sHom_{\cC}(A, -) \colon \cC \to \sSet$ has a left adjoint
		\begin{displaymath}
		A \otimes - \colon \sSet \to \cC
		\end{displaymath}
		natural in $A$. 
		\item The functor $A \otimes -$ is associative in the sense that there is an isomorphism
		\begin{displaymath}
		A \otimes (K \times L) \cong (A\otimes K) \otimes L
		\end{displaymath}
		natural in $A\in \cC$ and $K,L \in \sSet$.
		\item The functor $\sHom_{\cC}(-, B) \colon \cC^{\op} \to \sSet$ has a left adjoint
		\begin{displaymath}
		\shom_{\cC}(-, B)\colon \sSet \to \cC^{\op}
		\end{displaymath}
		natural in $B$.
	\end{enumerate}
\end{de}

\begin{de}
	A category $\cC$ is a simplicial model category, if it is both a model category and a simplicial category, and satisfies the additional axiom:
	\begin{enumerate}
		\item[\textbf{SM7}:] Suppose $j\colon A\to B$ is a cofibration and $q\colon X\to Y$ is a fibration.  Then 
		\begin{displaymath}
		\sHom_{\cC}(B, X) \xrightarrow{(j^{\ast},q_{\ast})} \sHom_{\cC}(A, X) \times_{\sHom_{\cC}(A, Y)} \sHom_{\cC}(B, Y)
		\end{displaymath}
		is a fibration in $\sSet$, which is trivial if $j$ or $q$ is trivial.
	\end{enumerate}
\end{de}

\begin{rem}\label{simplicial model category}
	\begin{enumerate}
		\item The above definitions imply that $\sHom_{\cC}(A, -) \colon \cC \to \sSet$ is right Quillen with left adjoint $A \otimes -$ when $A$ is cofibrant, $\sHom_{\cC}(-, X) \colon \cC^{\op} \to \sSet$ is right Quillen with left adjoint $\shom_{\cC}(-, X)$ when $X$ is fibrant, and $-\otimes K\colon \cC \to \cC$ is left Quillen with right adjoint $\shom_{\cC}(K, -)$ for $K\in \sSet$.
		\item There is a simplicial homotopy relation for morphisms $X\to Y$ in a simplicial model category $\cC$ (see \cite[Definition 9.5.2]{Hir09}), which coincides with the left and right homotopy relations if the source $X$ is cofibrant and the target $Y$ is fibrant (see \cite[Proposition 9.5.24]{Hir09}).  In particular, if $X\in \cC$ is cofibrant and $Y \in \cC$ is fibrant, then $\Hom_{\Ho(\cC)}(X, Y) \cong \pi_0\sHom_{\cC}(X,Y).$
	\end{enumerate}
\end{rem}

\begin{expl}\label{simplicial set}
	The category of simplicial sets $\sSet$ is a simplicial model category, with the specified classes of cofibrations, Kan fibrations and weak equivalences.  In this case, the tensor product is just the usual product, and $\shom_{\sSet}$ coincides with $\sHom_{\sSet}$ (see \cite[Proposition \RNum{1}.5.1, Theorem \RNum{1}.11.3 and Proposition \RNum{1}.11.5]{GJ09}).  
	
	Moreoover, $\sSet$ is cofibrantly generated.  Note that every simplicial set with finitely many non-degenerate simplices is small relative to all morphisms, we can take the set of generating cofibrations $I = \{ \partial \Delta^n \hookrightarrow \Delta^n \mid n\geq 0 \} $ (see \cite[Theorem \RNum{1}.11.2]{GJ09}), and the set of generating trivial cofibrations $J = \{ \Lambda^n_k \hookrightarrow \Delta^n \mid n\geq 1, 0\leq k \leq n \} $.  
\end{expl}

We explain how to generate cofibrantly generated simplicial model categories from already known ones.

For a complete and cocomplete category $\cC$, the category $s\cC$ has a simplicial category structure: for $A \in s\cC$ and $K \in \sSet$, we define $A \otimes K \in s\cC$ by $(A \otimes K)_n = \bigsqcup\limits_{k \in K_n} A_n$, where $\bigsqcup$ denotes the coproduct in $\cC$, with connecting morphisms naturally induced from those of $A$ and $K$.  Note that the definition is consistent for $\sSet$.

Let $\cC$ and $\cD$ be complete and cocomplete categories.  Suppose there is an adjoint pair of functors 
$$F\colon \cC \rightleftarrows \cD \colon G,$$ 
then the level-wise extended pair $F\colon s\cC \rightleftarrows s\cD \colon G$ is still an adjoint pair between the simplicial categories, and there are natural isomorphisms $F(A\times K) \cong F(A)\otimes K$ for $A\in s\cC$ and $K \in \sSet$ since $F$ preserves coproducts.  

\begin{pro}\label{generate cofibrantly generated simplicial model categories}
	Let notations be as above.  Suppose $s\cC$ is a cofibrantly generated simplicial model category with a set of generating cofibrations $I$ and a set of generating trivial cofibrations $J$.  Let $FI=\{Fi \mid i\in I\}$ and $FJ=\{Fj \mid j\in I\}$.  Suppose
	\begin{enumerate}
		\item[(a)] both $FI$ and $FJ$ permit the small object argument (see \defref{small}), and
		\item[(b)] $G\colon s\cD \to s\cC$ sends $\N$-compositions of pushouts of coproducts of $FJ$ to weak equivalences in $s\cC$.
	\end{enumerate}
	Then there is a cofibrantly generated simplicial model category structure on $s\cD$, such that $FI$ is a set of generating cofibrations and $FJ$ is a set of generating trivial cofibrations.  With this model category structure, $(F, G)$ is a Quillen pair.
\end{pro}

\begin{proof}
	See \cite[Theorem 11.3.2]{Hir09} and \cite[Theorem \RNum{2}.4.4]{GJ09}.
\end{proof}

\begin{rem}
	\begin{enumerate}
		\item The sets $FI$ and $FJ$ already determine the weak equivalences, fibrations and cofibrations of $\cD$.  They can be characterized as follows:
		\begin{enumerate}
			\item $f$ is a weak equivalence if and only if $Gf$ is a weak equivalence in $\cC$.
			\item $f$ is a fibration if and only if $Gf$ is a fibration in $\cC$.
			\item $f$ is a cofibration if and only if it is a retract of an $\N$-composition of pushouts of coproducts of $FI$ (see \corref{cofibration}).
		\end{enumerate}
		\item When $G$ preserves filtered colimits and the sources of $I$ and $J$ are small relative to all morphisms, assumption (a) holds by the proof of \cite[Theorem \RNum{2}.4.1]{GJ09}.  For a condition to ensure assumption (b), see \cite[Lemma \RNum{2}.5.1]{GJ09}.
	\end{enumerate}
\end{rem}

\begin{expl}\label{sMod and sCR}
	Let $R$ be a commutative ring.  We denote by $\sMod_R$ the category of simplicial $R$-modules and denote by $\sCR$ the category of simplicial commutative rings.  Assumptions (a) and (b) of \propref{generate cofibrantly generated simplicial model categories} hold in the following situations:
	\begin{enumerate}
		\item Consider the adjoint pair $F\colon \sSet \rightleftarrows \sMod_R \colon G$, where $F$ is the free module functor and $G$ is the forgetful functor.  We take $I$ and $J$ as in \explref{simplicial set}.  Then $\sMod_R$ is a cofibrantly generated simplicial model category.  In the next section we will show that the model structure of $\sMod_R$ is essentially the same as the model structure of $\Ch_{\geq 0}(R)$ defined in \explref{chain complex}, and a more convenient choice of generating cofibrations and generating trivial cofibrations is by transfering those of $\Ch_{\geq 0}(R)$ in \explref{chain complex} via the Dold-Kan equivalence.
		\item Consider the adjoint pair $F\colon \sMod_{\Z} \rightleftarrows \sCR \colon G$, where $F$ is the symmetric algebra functor and $G$ is the forgetful functor.  We take $I = \{ 0 \to \Z \} \cup \{ \DK(\Z[n] \to \Z\left\langle n+1\right\rangle) \mid n\geq 0 \} $ and $J= \{ \DK(0 \to \Z\left\langle n+1\right\rangle) \mid n\geq 0 \}$ as remarked above.  Then $\sCR$ is a cofibrantly generated simplicial model category.  The weak equivalences and fibrations are those of $\sMod_{\Z}$, and the cofibrations are retracts of $\N$-compositions of pushouts of coproducts of $FI$.
	\end{enumerate}
\end{expl}

\subsubsection{Dold-Kan correspondence}

Let $R$ be a commutative ring.  Our goal here is to recall an equivalence of model categories between $\sMod_R$ and $\Ch_{\geq 0}(R)$.  

When $M \in \sMod_R$, we write $M_n$ for the $R$-module on $n$-th simplicial degree.  Let $N(M)$ be the chain complexes of $R$-modules with $N(M)_n=\bigcap\limits_{i=0}^{n-1}\ker(d_i) \subseteq M_n$ and $n$-th differential map 
\begin{displaymath}
(-1)^n d_n \colon \bigcap\limits_{i=0}^{n-1}\ker(d_i)\subseteq M_n \to \bigcap\limits_{i=0}^{n-2}\ker(d_i)\subseteq M_{n-1}.
\end{displaymath}
Then obviously $M\mapsto N(M)$ is natural in $M$, and we call $N(M)\in\Ch_{\geq 0}(R)$ the normalized complex of $M$.  

The Dold-Kan functor $\DK\colon \Ch_{\geq 0}(R) \to \sMod_R$ is the quasi-inverse of $N$.  Explicitly, for a chain of $R$-modules $C_{\ast}=(C_0 \leftarrow C_1 \leftarrow C_2 \leftarrow \dots)$, we define $\DK(C_{\ast}) \in \sMod_R$ as follows:
\begin{enumerate}
	\item $\DK(C_{\ast})_n = \bigoplus\limits_{[n] \twoheadrightarrow [k]} C_k$.
	\item For $\theta\colon [m] \to [n]$, we define the corresponding $\DK(C_{\ast})_n \to \DK(C_{\ast})_m$ 
	on each component of $\DK(C_{\ast})_n$ indexed by 
	$[n] \stackrel{\sigma}{\twoheadrightarrow} [k]$ as follows: suppose 
	$[m] \stackrel{t}{\twoheadrightarrow} [s] \stackrel{d}{\hookrightarrow} [k]$ 
	is the epi-monic factorization of the composition $[m] \stackrel{\theta}{\rightarrow}  [n] \stackrel{\sigma}{\twoheadrightarrow} [k]$,
	then the map on component $[n] \stackrel{\sigma}{\twoheadrightarrow} [k]$ is 
	\begin{displaymath}
	C_k \stackrel{d^{\ast}}{\rightarrow} C_s \hookrightarrow \bigoplus\limits_{[m] \twoheadrightarrow [r]} C_r.
	\end{displaymath}
\end{enumerate}

\begin{rem}
	Let $M[1]$ be the chain complex with $M$ on degree $1$ and $0$ elsewhere.  Then $\DK(M[1])$ is the nerve of the abelian group $M$ (see \explref{nerve}).
\end{rem}

\begin{thm}
	\begin{enumerate}
		\item (Dold-Kan) The functors $\DK$ and $N$ are quasi-inverse and form an equivalence of categories.  Moreover, two morphisms $f,g \in \Hom_{\sMod_R }(M, N)$ are simplicially homotopic if and only if $N(f)$ and $N(g)$ are chain homotopic.  
		\item The functors $\DK$ and $N$ preserve the model category stuctures of $\Ch_{\geq 0}(R)$ and $\sMod_R$ defined above.
	\end{enumerate}
\end{thm}

\begin{proof}
	See \cite[Theorem 8.4.1]{Weib94} and \cite[Lemma 2.11]{GJ09}.  Note that (1) is valid for any abelian category instead of $\sMod_R$.  
\end{proof}

\begin{rem} 
	Let $\Ch(R)$ be the category of complexes $(C_i)_{i\in\Z}$ of $R$-modules and $\Ch_{\geq 0}(R)$ the subcategory of complexes for which $C_i =0$ for $i<0$.
	The category $\Ch_{\geq 0}(R)$ is naturally enriched over simplicial $R$-modules, and we have 
	$$\sHom_{\Ch_{\geq 0}(R)}(C_{\ast}, D_{\ast}) \cong \sHom_{\sMod_R}(\DK(C_{\ast}), \DK(D_{\ast})).$$
	Given $C_{\ast}, D_{\ast} \in \Ch_{\geq 0}(R)$.  Let $[C_{\ast}, D_{\ast}] \in \Ch(R)$ 
	be the mapping complex, more precisely, $[C_{\ast}, D_{\ast}]_n = \prod_m \Hom_R(C_m, D_{m+n})$ 
	and the differential maps are natural ones.  Let $\tau_{\geq0}$ be the functor which sends a chain complex $X_{\ast}$ 
	to the truncated complex
	\begin{displaymath}
	0 \leftarrow \ker(X_0 \to X_{-1}) \leftarrow X_1 \leftarrow \dots
	\end{displaymath}
	Then there is a weak equivalence
	\begin{displaymath}\label{tronc} 
	\sHom_{\Ch_{\geq 0}(R)}(C_{\ast}, D_{\ast})\simeq \DK (\tau_{\geq 0}[C_{\ast}, D_{\ast}])
	\end{displaymath}
	(see \cite[Remark 11.1]{Lur09}).  And it's clear that $\pi_n\sHom_{\Ch_{\geq 0}(R)}(C_{\ast}, D_{\ast})$ is isomorphic 
	to the chain homotopy classes of maps from $C_{\ast}$ to $D_{\ast +n}$.
\end{rem}

\subsubsection{Simplicial commutative rings}\label{simplicial commutative ring}

In \explref{sMod and sCR} we introduce a model category structure on $\sCR$ such that the fibrations and weak equivalences are those of $\sMod_{\Z}$ (or equivalently $\sSet$).  The description of cofibrations is a bit complicated, but we mention that a cofibration $A\to B$ must be degreewise flat (see \cite[Lemma 7.10.2]{Gil13}).  One can deduce from this fact that the degreewise tensor product $-\otimes_A B \colon {}_A \backslash\sCR \to {}_B \backslash\sCR$ 
is a left Quillen functor, so it makes sense to define its total left derived functor 
$$- \lotimes_A B \colon \Ho({}_A \backslash\sCR) \to  \Ho({}_B \backslash\sCR).$$  
We also use $C \uotimes_A B$ to denote some $c(C) \otimes_A B \in {}_B \backslash\sCR$, where $c(C)$ is a cofibrant replacement of $C$ in ${}_A \backslash\sCR$; it is well defined up to weak equivalence and it represents $C \lotimes_A B$.

In what follows, we will explain the graded commutative ring structure on $\pi_{\ast}(A)$ for $A\in \sCR$.  Here it's natural to consider together the modules over simplicial commutative rings.

\begin{de}\label{modules over simplicial commutative rings}
	Fix $A\in \sCR$.  We define the category $\Mod(A)$ as follows: the objects are simplicial abelian groups $M$ such that each $M_n$ is an $A_n$-module and each morphism $[m] \to [n]$ of $\bm{\Delta}$ induces $M_n\to M_m$ compatible with $A_n\to A_m$, and the morphisms from $M$ to $N$ consist of $A_n$-module morphisms $M_n\to N_n$ ($n\geq 0$) compatible with $\bm{\Delta}$-morphisms $[m] \to [n]$.	
\end{de}

Note if $\underline{A} \in \sCR$ is the constant simplicial ring associated to $A \in \CR$, then $\Mod(\underline{A})$ is naturally isomorphic to $\sMod_A$.  

For $A \in \sCR$ and $M \in \Mod(A)$, the unnormalized chain complex is $C(M) = \bigoplus_{n=0}^{\infty} M_n$ with differential
\begin{displaymath}
\sum_{i=0}^{n} (-1)^i d_i \colon M_n \to M_{n-1}.
\end{displaymath}
It's clear that the above construction is natural in $M$.  Moreover, the inclusion of abelian group complexes $N(M) \to C(M)$ (by the way one can check the boundary and cycle in $N(M)_n$ are $A_n$-modules) is a homotopy equivalence and induces $H_{\ast}(N(M)) \iso H_{\ast} (C(M))$ (see \cite[Lemma 5.1.2]{Gil13}).

In the following we define multiplications of $C(A)$ on $C(M)$, making $C(M)$ a differential graded module over $C(A)$ (see Section \ref{graded}).

For $m,n \geq 0$, the set of surjective morphisms $[m+n] \to [m]$ of $\bm{\Delta}$ is in one-to-one correspondence with the set $\{\sigma=(\sigma_i)_{i=1}^m \mid 1\leq \sigma_1 < \sigma_2 < \dots <\sigma_m \leq m+n  \} $, where $\sigma=(\sigma_i)_{i=1}^m$ corresponds to the morphism $[m+n] \to [m]$ sending $\sigma_i, \sigma_i+1, \dots, \sigma_{i+1}-1$ to $i$ (we put $\sigma_0 = 0$ and $\sigma_{m+1} = m+n+1$ for convenience).  Let $P_{m,n}$ be the set of permutations $(\sigma, \tau)$ of $\{1,2,\dots,m+n\}$ where $\sigma=(\sigma_i)_{i=1}^m$ satisfies  $1\leq \sigma_1 < \sigma_2 < \dots <\sigma_m \leq m+n$ and $\tau=(\tau_i)_{i=1}^n$ satisfies  $1\leq \tau_1 < \tau_2 < \dots <\tau_n \leq m+n$.  Then $(\sigma, \tau) \in P_{m,n}$ determines surjective morphisms $\sigma\colon[m+n] \to [m]$ and $\tau\colon[m+n] \to [n]$.  Let $\text{sign}(\sigma, \tau)$ be the sign of the permutation $(\sigma, \tau)$.  Then for $(\sigma, \tau) \in P_{m,n}$, we have $(\tau, \sigma) \in P_{n,m}$ and $\text{sign}(\sigma, \tau) = (-1)^{mn}\text{sign}(\tau,\sigma)$.

The multiplication of $C(A)$ on $C(M)$ is defined by 
\begin{displaymath}
a\cdot x = \sum\limits_{(\sigma, \tau) \in P_{m,n}} \text{sign}(\sigma, \tau) A(\sigma)(a) M(\tau)(x),
\end{displaymath}
for $a\in A_m$ and $x\in M_n$, where $A(\sigma) \colon A_m \to A_{m+n}$ corresponds to $\sigma\colon[m+n] \to [m]$ and $M(\tau) \colon M_n \to M_{m+n}$ corresponds to $\tau \colon[m+n] \to [n]$.  Then one has the following lemma:

\begin{lem}
	Let $A\in \sCR$ and let $M\in \Mod(A)$. 
	\begin{enumerate}
		\item $C(A)$ is a strictly graded commutative ({\it i.e.}, $a\cdot a = 0$ for every $a\in A_i$ for every odd $i$) differential graded ring.  Moreover, with the multiplication induced from $C(A)$, the normalized chain complex $N(A)$ is a sub-differential graded ring of $C(A)$.
		\item $C(M)$ is a differential graded module over $C(A)$.  Moreover, with the multiplication induced from $C(M)$, the normalized chain complex $N(M) \subseteq C(M)$ is a differential graded module over $N(A) \subseteq C(A)$. 
		\item The multiplication is well-defined for homology groups.  In particular, under the isomorphisms $\pi_{\ast}(A) \cong H_{\ast}(N(A)) \cong H_{\ast}(C(A))$ and $\pi_{\ast}(M) \cong H_{\ast}(N(M)) \cong H_{\ast}(C(M))$, $\pi_{\ast}(A)$ is a graded commutative ring and $\pi_{\ast}(M)$ is a graded $\pi_{\ast}(A)$-module.
	\end{enumerate}
\end{lem}

\begin{proof}
	See \cite[Lemma 8.3.2]{Gil13}.
\end{proof}

\subsection{Representability of functors}

\subsubsection{Simplicial Artinian rings}

Recall that $\cO$ is the ring of integers in a $p$-adic number field $K$, and $k$ is the residue field of $\cO$.  We regard $\cO$ and $k$ as constant objects in $\sCR$.  

For $A\in \sCRO$, we have shown that $\bigoplus_i\pi_i A$ is naturally a graded commutative $\cO$-algebra.  Recall that $\sCR$ is cofibrantly generated, so we can fix a functorial factorization $\cO \hookrightarrow c(A) \stackrel{\sim}{\twoheadrightarrow} A$ for $A\in \sCRO$.  Now let's define an Artinian subcategory of $ \sCROk$.

\begin{de}\label{simplicial Artinian algebra}
	The simplicial Artinian $\cO$-algebras over $k$, which we denote by $\sArt$, is the full subcategory of $ \sCROk$ consisting of objects $A \in \sCROk$ such that:
	\begin{enumerate}
		\item $\pi_0 A$ is an Artinian local $\cO$-algebra in the usual sense.
		\item $\pi_{\ast} A = \oplus_{i\geq 0}\pi_i A$ is finitely generated as a module over $\pi_0 A$. 
	\end{enumerate}
\end{de}

Note that $\sArt$ is not a model category, and cofibrations, fibrations and weak equivalences in $\sArt$ are used to indicate those in $\sCROk$.  Nevertheless, $\sArt$ is closed under weak equivalences since the definition only involves homotopy groups.  We also remark that every $A\in \sArt$ is fibrant since $A \to k$ is degreewise surjective.

\begin{expl}
	If $M \in \sMod_k$ and $\dim_k (\pi_{\ast}(M)) < \infty$, then the object $k\oplus M \in \sCROk$ defined by square-zero extension on each simplicial degree is an object of $\sArt$.  In particular, $k \oplus \DK(k[n]) \in \sArt$ for $n\geq 0$ (here $k[n]$ is the chain complex with $k$ on $n$-th degree and $0$ elsewhere).  For simplicity we write $k\oplus k[n]$ for $k \oplus \DK(k[n])$.
\end{expl}

\subsubsection{Formally cohesive functors}

\begin{de}
	A functor $\cF\colon \sArt \to \sSet$ is called formally cohesive if it satisfies the following conditions:
	\begin{enumerate}
		\item $\cF$ is homotopy invariant ({\it i.e.} preserves weak equivalences).
		\item Suppose that 
		\begin{displaymath}
		\xymatrix{
			A \ar[r]\ar[d] & B \ar[d] \\
			C \ar[r] & D}
		\end{displaymath}
		is a homotopy pullback square with at least one of $B \to D$ and $C \to D$ being degreewise surjective ({\it i.e.}, a fibration with surjective $\pi_0$, see \cite[Lemma \RNum{3}.2.11]{GJ09}), then
		\begin{displaymath}
		\xymatrix{
			\cF(A) \ar[r]\ar[d] & \cF(B) \ar[d] \\
			\cF(C) \ar[r] & \cF(D)}
		\end{displaymath}
		is a homotopy pullback square (in this case we say $\cF$ preserves homotopy pullbacks for simplicity).
		\item $\cF(k)$ is contractible.
	\end{enumerate}
\end{de}

\begin{expl}\label{representable functor example}
	If $R \in \sCROk$ is cofibrant, then the functor 
	$$\sHom_{\sCROk}(R, -) \colon \sArt \to \sSet$$ 
	is a restriction of a right Quillen functor and obviously Kan-valued.  In addition, it extends to 
	\begin{displaymath}
	\sHom_{\sCROk}(A, B) \to \sHom_{\sSet}(\sHom_{\sCROk}(R, A), \sHom_{\sCROk}(R, B))
	\end{displaymath}
	(this is called the simplicial enrichment), which is given by the adjoint 
	$$\sHom_{\sCROk}(A, B) \times  \sHom_{\sCROk}(R, A) \to \sHom_{\sCROk}(R, B)$$
	defined just below \cite[Lemma \RNum{2}.2.2]{GJ09}.  Moreover, the functor is formally cohesive:
	\begin{enumerate}
		\item Since a right Quillen functor preserves weak equivalences between fibrant objects (\cite{Hir09}, Proposition 8.5.7) and every object of $\sArt$ is fibrant, $\sHom_{\sCROk}(R, -)$ is homotopy invariant.
		\item Note that $B \times_D^h C \in \sArt$ (see \cite[Lemma 2.3]{GV18}).  Write $\cF=\sHom_{\sCROk}(R, -)$ for simplicity.  By \lemref{holim commute} we have $\mathbf{R}\cF(B \times_D^h C) \cong \mathbf{R}\cF(B) \times_{\mathbf{R}\cF(D)}^h \mathbf{R}\cF(C)$ in the homotopy category, then use the fact that $\cF$ is homotopy invariant, we get the chain of weak equivalences $\cF(A) \simeq \cF(B \times_D^h C) \simeq \cF(B) \times_{\cF(D)}^h \cF(C).$
		\item $\sHom_{\sCROk}(R, k)$ is obviously contractible.
	\end{enumerate}
\end{expl}

We can construct formally cohesive functors from known ones:

\begin{lem}\label{formally cohesive functors from known ones}
	\begin{enumerate}
		\item Let $X$ be a simplicial set and let $\cF$ be a Kan-valued, homotopy invariant functor.  Then the functor $A \mapsto \sHom_{\sSet}(X, \cF(A))$ is formally cohesive (resp. preserves homotopy pullbacks) if $\cF$ is formally cohesive (resp. preserves homotopy pullbacks).
		\item Let $\cC$ be a small category and let $(\cF_c)_{c\in \cC}$ be a $\cC$-system of homotopy invariant functors from $\sArt$ to $\sSet$.  Define $\cF = \holim_{c\in \cC} \, \cF_c$ to be the objectwise homotopy limit, then $\cF$ is formally cohesive (resp. preserves homotopy pullbacks) if every $\cF_c$ $(c\in \cC)$ is formally cohesive (resp. preserve homotopy pullbacks).
		\item Let $I$ be a small filtered category and let $(\cF_i)_{i\in I}$ be a filtered system of homotopy invariant functors.  Define $\cF(A) = \hocolim_I \, \cF_i(A)$.  Then $\cF$ is formally cohesive (resp. preserves homotopy pullbacks) if all $\cF_i$ $(i\in I)$ are formally cohesive (resp. preserve homotopy pullbacks).
	\end{enumerate}
\end{lem}

\begin{proof}
	First note $\sHom_{\sSet}(X, \cF(-))$ and $\holim_{c\in \cC} \, \cF_c$ are homotopy invariant under our assumptions, then since both $\sHom_{\sSet}(X, -)$ and the homotopy limit functor are right Quillen, (1) and (2) are consequences of \lemref{holim commute} (see also \cite[Lemma 4.29 and Lemma 4.30]{GV18}).  Part (3) follows from \lemref{filtered colimit} below.
\end{proof}

\begin{lem}\label{filtered colimit}
	Let $I$ be a small filtered category.
	\begin{enumerate}
		\item The functor $\varinjlim_I\colon \sSet^I_{\proj} \to \sSet$ preserves fibrations and trivial fibrations.
		\item The functor $\varinjlim_I\colon \sSet^I_{\proj} \to \sSet$ preserves weak equivalences.
		\item The functor $\varinjlim_I\colon \sSet^I \to \sSet$ commutes with homotopy pullbacks.
	\end{enumerate}
\end{lem}

\begin{proof}
	\begin{enumerate}
		\item Fibrations and trivial fibrations are characterized by right lifting properties with respect to morphisms $\partial \Lambda^n_k \hookrightarrow \Delta^n$ and $\partial \Delta^n \hookrightarrow \Delta^n$ respectively, and all objects involved are small in the sense of Quillen, so the result follows.
		\item By part (1) and \cite[Proposition 8.5.7]{Hir09}, the functor $\varinjlim_I$ preserves weak equivalences between fibrant objects.  The result follows because Kan's $\Ex^{\infty}$ functor (see \cite[\RNum{3}.4]{GJ09}) gives fibrant replacements and preserves filtered colimits.
		\item Let $(B_i \rightarrow D_i \leftarrow C_i)_{i\in I}$ be a system of diagrams.  Let $B_i^{\prime} \rightarrow D_i^{\prime} \leftarrow C_i^{\prime}$ be a fibrant replacement of $B_i \rightarrow D_i \leftarrow C_i$, then by lifting properties $(B_i^{\prime} \rightarrow D_i^{\prime} \leftarrow C_i^{\prime})_{i\in I}$ forms a direct system.  From parts (1) and (2), we see $\varinjlim_IB_i^{\prime} \rightarrow \varinjlim_ID_i^{\prime} \leftarrow \varinjlim_I C_i^{\prime}$ is fibrant and is weakly equivalent to $\varinjlim_IB_i \rightarrow \varinjlim_I D_i \leftarrow \varinjlim_I C_i$, so 
		\begin{displaymath}
		\varinjlim_IB_i\times_{\varinjlim_ID_i}^h \varinjlim_I C_i \simeq  \varinjlim_IB_i^{\prime}\times_{\varinjlim_ID_i^{\prime}} \varinjlim_IC_i^{\prime} \simeq \varinjlim_I B_i^{\prime}\times_{D_i^{\prime}}C_i^{\prime},
		\end{displaymath}
		where the second weak equivalence is because filtered colimits commute with finite limits.
	\end{enumerate}
\end{proof}

\subsubsection{Pro-representable functors}

\begin{de}
	Let $\cF$ and $\cG$ be two functors from $\sArt$ to $\sSet$.  
	\begin{enumerate}
		\item A natural transformation $T\colon \cF \to \cG$ is a weak equivalence if it induces weak equivalences $\cF(A) \iso \cG(A)$ for all $ A\in \sArt$.
		\item $\cF$ and $\cG$ are weakly equivalent if there exists a finite zig-zag of weak equivalences between $\cF$ and $\cG$.
	\end{enumerate}
\end{de} 

\begin{de}
	A functor $\cF \colon \sArt \to \sSet$ is pro-representable, if there is a projective system $R = (R_n)_{n\in \N}$ with each $R_n \in \sArt$ cofibrant, such that $\cF$ is weakly equivalent to $\varinjlim\limits_n \sHom_{\sCROk}(R_n, -)$.  
\end{de}

In this case we say $R = (R_n)$ is a representing (pro-)ring for $\cF$ (we will often omit "pro" for convenience).  For a pro-ring $R = (R_n)$ we shall write 
$$\sHom_{\sCROk}(R, -) = \varinjlim\limits_n\sHom_{\sCROk}(R_n, -)$$
for simplicity.

\begin{rem}
	\begin{enumerate}
		\item The pro-representability defined above is called the sequential pro-representability in \cite{GV18}, but we will only encounter this case.
		\item By \lemref{filtered colimit}, one can replace the colimit by the homotopy colimit.  As pointed out in \cite[Section 2.6]{GV18}, the homotopy colimit is easier to map out of, while the usual colimit preserves fibrations. 
		\item The representing ring is not uniquely determined up to natural isomorphism.  However, since filtered colimits of $\sSet$ commute with $\pi_0$, it's easy to see that the representing ring is uniquely determined up to natural isomorphism as a pro-object in $\Ho(\sCROk)$.  So if $R$ pro-represents $\cF$ then $\pi_{\ast} R$ is well-defined.
	\end{enumerate}
\end{rem}

We expect that a natural transformation of pro-representable functors induces a morphism between the corresponding pro-rings, at least modulo homotopy.  For this we require the representing pro-ring $R$ to be nice in the sense of \cite[Definition 2.23]{GV18}.  When $R = (R_n)$ is degreewise cofibrant, then the niceness condition means exactly that the pro-ring $R$ is Reedy fibrant in the standard Reedy model category $(\sCROk)^{\N}$, so one can always make such a choice by taking fibrant replacements in the Reedy model category.

\begin{lem}\label{functor  map to ring map}
	Let $\cF$ and $\cG$ be two Kan-valued functors from $\sArt$ to $\sSet$.  We use $T\colon \cF \dashrightarrow \cG$ to denote a zigzag of natural tansformations 
	\begin{displaymath}
	\cF \stackrel{\sim}{\leftarrow} \cF_1 \to \cF_2 \stackrel{\sim}{\leftarrow} \cF_3 \to \cF_4 \stackrel{\sim}{\leftarrow} \dots \to \cG
	\end{displaymath}
	where all left arrows are weak equivalences.  Suppose $R = (R_n)$ (resp. $S = (S_n)$) is a representing pro-ring for $\cF$ (resp. $\cG$) and $R$ is fibrant in the Reedy model category ({\it i.e.}, nice), then there is a morphism $S \xrightarrow{\alpha} R$ of pro-simplicial rings such that for $A\in \sArt$, the diagram 
	\begin{displaymath}
	\xymatrix{
		\cF(A) \ar@{-->}[r]^{T}  & \cG(A)  \\
		\sHom_{\sCROk}(R , A) \ar[r]^{\alpha^{\ast}}\ar@{-->}[u]^{\simeq} & \sHom_{\sCROk}(S , A) \ar@{-->}[u]^{\simeq}
	}
	\end{displaymath}
	is commutative after taking homotopy groups $\pi_i$ $(i\geq 0)$ (note the dotted arrows become true arrows after taking homotopy groups, since	weak equivalences become isomorphisms). 
\end{lem}

\begin{proof}
	First of all we can replace the zigzag $T$ by $\cF \stackrel{\sim}{\leftarrow} \cF^{\ast} \to \cG$, where $\cF^{\ast}$ is the homotopy limit of the diagram 
	\begin{displaymath}
	\cF \stackrel{\sim}{\leftarrow} \cF_1 \to \cF_2 \stackrel{\sim}{\leftarrow} \cF_3 \to \cF_4 \stackrel{\sim}{\leftarrow} \dots \to \cG 
	\end{displaymath}
	(see discussions around \cite[(7.3)]{GV18}).  Then as \cite[Lemma 2.25]{GV18} there exists horizontal arrows in the second and third lines which make the diagram 
	\begin{displaymath}
	\xymatrix{
		\cF^{\ast}(A) \ar[r] & \cG(A)  \\
		\hocolim_n \, \sHom_{\sCROk}(R_n , A) \ar[r]\ar[u]\ar[d] & \hocolim_n \, \sHom_{\sCROk}(S_n , A) \ar[u]\ar[d] \\
		\sHom_{\sCROk}(R , A) \ar[r] &\sHom_{\sCROk}(S , A)
	}
	\end{displaymath}
	commute modulo simplicial homotopy.  Note the niceness of $R$ implies that 
	$$\lim_n \sHom_{\sCROk}(S , R_n) \to \holim_n \, \sHom_{\sCROk}(S , R_n)$$
	is a weak equivalence, and the arrow in the third line exists by the enriched Yoneda's lemma.
\end{proof}

By \lemref{filtered colimit} and \explref{representable functor example}, any pro-representable functor is formally cohesive.  Conversely, Lurie's criterion asserts that a formally cohesive functor is pro-representable if additionally its tangent complex is not far from the tangent complexes of simplicial commutative rings.  We will introduce tangent complexes and Lurie's criterion below.

\subsubsection{(Co)tangent complexes of simplicial commutative rings}\label{tangent complexes of simplicial commutative rings}

Let's recall Quillen's cotangent and tangent complexes of simplicial commutative rings.  

Let $\Omega_{R/\cO}$ be the module of differentials with the canonical $R$-derivation $d\colon R \to \Omega_{R/ \cO}$ for an $\cO$-algebra $R$.  Let $\Der_{\cO}(R, -)$ be the covariant functor which sends an $R$-module $M$ to the $R$-module 
\begin{displaymath}
\Der_{\cO}(R, M) = \{ D\colon R \to M \mid D \text{ is } \cO \text{-linear and } D(xy) = xD(y)+yD(x),\text{ } \forall x,y\in R \}.
\end{displaymath}
It's well-known that $\Hom_R(\Omega_{R/\cO}, -)$ is naturally isomorphic to $\Der_{\cO}(R, -)$ via $\phi \mapsto \phi \circ d$.  

For any $k$-module $M$ and any $R \in {}_{\cO}\backslash \CR/_{k}$, we have natural isomorphisms
\begin{displaymath}
\Hom_k (\Omega_{R/\cO} \otimes_R k, M) \cong \Der_{\cO}(R, M)\cong \Hom_{{}_{\cO}\backslash \CR/_{k}}(R, k \oplus M).
\end{displaymath}
where $k\oplus M$ is the $k$-algebra with square-zero ideal $M$. So the functor $R \mapsto \Omega_{R/\cO} \otimes_R k$ is left adjoint to the functor $M \mapsto k \oplus M$.  

The above adjunction has level-wise extensions to simplicial categories (see \cite{GJ09} 
Lemma \RNum{2}.2.9 and Example \RNum{2}.2.10).  For $R \in {}_{\cO}\backslash \sCR$, we can form degreewisely $\Omega_{R/\cO}\otimes_R k \in \sMod_k$, and we have 
\begin{displaymath}
\sHom_{\sMod_k }(\Omega_{R/\cO} \otimes_R k, M) \cong \sHom_{\sCROk}(R, k \oplus M).
\end{displaymath}
The functor $M \mapsto k \oplus M$ from $\sMod_k$ to $\sCROk$ 
preserves fibrations and weak equivalences (we may see this via the Dold-Kan correspondence), 
so the left adjoint functor $R \mapsto \Omega_{R/\cO} \otimes_R k$ is left Quillen and it admits a total left derived functor.  

\begin{de}
	For $R \in {}_{\cO}\backslash \sCR$, we define the cotangent complex of $R$ to be 
	$$L_{R/\cO} = \Omega_{c(R)/\cO} \otimes_{c(R)} R \in \Mod(R) $$
	(here $\otimes$ is the degreewise tensor product, and see \defref{modules over simplicial commutative rings} for $\Mod(R)$).
\end{de}

Then the total left derived functor of $R \mapsto \Omega_{R/\cO} \otimes_R k$ is $R \mapsto L_{R/\cO} \otimes_R k$.  

By construction, $L_{R/\cO} \otimes_R k$ is cofibrant as it's the image of the cofibrant object $c(R)$ under a total left derived functor, and it is fibrant in $\sMod_k$ (all objects are fibrant there).  It follows that $L_{R/\cO} \otimes_R k$ is determined up to homotopy equivalence (by the Whitehead theorem \cite[Theorem 7.5.10]{Hir09}).  Using the Dold-Kan equivalence, we can form the normalized complex (determined up to homotopy equivalence)
$$N( L_{R/\cO} \otimes_R k )\in \Ch_{\geq 0}(k).$$
We will often abuse the language and also use $L_{R/\cO}  \otimes_R k$ to denote its image under $N$.

Recall that for $M,N\in\Ch(k)$, the internal Hom $[M,N]\in\Ch(k)$ is defined as 
$$[M,N]_n=\prod_m \Hom_k(M_m, N_{m+n}).$$
When $R \in \sCROk$ and $C_{\ast} \in \Ch_{\geq 0}(k)$, we have (by \remref{tronc}):
\begin{align*}
\sHom_{\sCROk}( c(R), k \oplus \DK(C_{\ast})) 
&\cong \sHom_{\sMod_k}(L_{R/\cO}  \otimes_R k,\DK(C_{\ast})) 
\\ &\simeq  \DK(\tau_{\geq 0}[L_{R/\cO}  \otimes_R k, C_{\ast}]).
\end{align*}

\begin{de}
	The tangent complex $\t R$ is the internal hom complex $[L_{R/\cO} \otimes_R k, k] \in \Ch_{\leq 0}(k).$
\end{de}

Note that $\t R$ is well-defined up to chain homotopy equivalence since it is the case for $L_{R/\cO} \otimes_R k$.  Also note $H_{-i}(\t R) = 0$ for $i <0$.  When convenient, we may identify $\Ch_{\leq 0}(k)=\Ch^{\geq 0}(k)$ via $C^i=C_{-i}$.  

\begin{rem}
	For a field $k$, the functor $\Hom_k (- ,k)$ on $k$-vector spaces is exact and there are no significant differences between $\t R$ and $L_{R/\cO}  \otimes_R k$.  On the other hand, in studying the adjoint Selmer groups, \cite{TU21} considers derived deformations over $\rho_B \colon \Gamma_S \to G(B)$ for some Artinian $\cO$-algebra $B$, where $L_{R/\cO}  \otimes_R B$ appears to be the more appropriate object.
\end{rem}

\subsubsection{(Co)tangent complexes of formally cohesive functors and Lurie's criterion}

The tangent complexes of formally cohesive functors is constructed in \cite[Section 4]{GV18}.  The key result is the following:

\begin{pro}
	Let $\cF \colon \sArt \to \sSet$ be a formally cohesive functor.  Then there exists $L_{\cF} \in \Ch(k)$ such that $\cF(k\oplus \DK(C_{\ast}))$ is weakly equivalent to $\DK(\tau_{\geq0} [L_{\cF},C_{\ast}])$ for every $C_{\ast} \in \Ch_{\geq 0}(k)$ with $H_{\ast}(C_{\ast})$ finite.
\end{pro}

\begin{proof}
	See \cite[Lemma 4.25]{GV18}.
\end{proof}

\begin{de}
	Let $\cF \colon \sArt \to \sSet$ be a formally cohesive functor.  We define $\t \cF = [L_{\cF}, k]$ to be the tangent complex of $\cF$.  
\end{de}

\begin{rem}
	It's easy to see that $L_{\cF}$ and $\t\cF$ are well-defined up to quasi-isomorphisms.  Comparing with above discussions for simplicial commutative rings, we call $L_{\cF}$ the cotangent complex of $\cF$.
\end{rem}

\begin{rem}
	In \cite[Section 4]{GV18}, the authors showed the existence of tangent complexes for general formally cohesive functors.  On the other hand, we can explicitly calculate the tangent complexes of the derived deformation functors we are interested in.
\end{rem}

It's convenient to regard $\t \cF$ as a cochain complex via $C^i=C_{-i}$, and we denote $\t^i\cF = H_{-i}\t\cF$.  Then for $i, n\geq 0$, we have $\pi_i \cF(k \oplus k[n]) \cong H_i ([L_{\cF}, k[n]]) \cong H_{i-n} ([L_{\cF}, k]) \cong \t^{n-i} \cF$.  

If $R \in \sCROk$ is cofibrant and $\cF_R = \sHom_{\sCROk}(R , -)$, then the cotangent complexes $L_{\cF_R}$ and $L_{R/\cO} \otimes_R k$ are quasi-isomorphic, since 
$$\DK(\tau_{\geq0} [L_{\cF_R},k[n]]) \simeq \sHom_{\sCROk}(R, k\oplus k[n]) \simeq \DK(\tau_{\geq0} [L_{R/\cO} \otimes_R k,k[n]]).$$

Now we see any pro-representable functor $\cF$ is formally cohesive and satisfies $\t^i \cF = 0$ ($\forall i <0$).  The following theorem gives the converse:

\begin{thm}[Lurie's criterion]\label{Lurie}
	Let $\cF$ be a formally cohesive functor.  If $\dim_k \t^i \cF$ is finite for every $i \in \Z$ and $\t^i\cF = 0$ for every $i <0$, then $\cF$ is (sequentially) pro-representable.
\end{thm}	

\begin{proof}
	See \cite[Corollary 6.2.14]{Lur04} and \cite[Theorem 4.33]{GV18}.
\end{proof}

The following lemma illustrates the conservativity of the tangent complex functor:

\begin{lem}\label{conservativity of tangent complex}
	Suppose $\cF,\cG \colon \sArt \to \sSet$ are formally cohesive functors.  Then a natural transformation $\cF \to \cG$ is a weak equivalence if and only if it induces isomorphisms $\t^i\cF \to \t^i\cG$ for all $i$.
\end{lem}

\begin{proof}
	One direction is clear and we prove the other.  If the natural transformation induces isomorphisms $\t^i \cF \to \t^i\cG$ for all $i$, then $\cF(k\oplus k[n]) \to \cG(k\oplus k[n])$ is a weak equivalence for every $n\geq 0$.  Hence by simplicial artinian induction \cite[Lemma 2.8]{GV18} and the formal cohesiveness of $\cF$ and $\cG$, the map $\cF(A) \to \cG(A)$ is a weak equivalence for every $A \in \sArt$.
\end{proof}

The following lemma says that tangent complexes commute with homotopy limits:

\begin{lem}\label{tangent complex calculation}
	Let $\cC$ be a small category and let $(\cF_c)_{c\in \cC}$ be a $\cC$-system of formally cohesive functors from $\sArt$ to $\sSet$.  Define $\cF = \holim_{c\in \cC} \, \cF_c$ to be the objectwise homotopy limit, then $\t \cF = \holim_{c\in \cC}  \, \t \cF_c$.  In particular, for the objectwise homotopy pullback diagram 
	\begin{displaymath}
	\xymatrix{
		\cF \ar[r]^{f_1}\ar[d]^{f_2} & \cF_1 \ar[d]^{p_1} \\
		\cF_2 \ar[r]^{p_2} & \cF_0
	}
	\end{displaymath}
	with $\cF_i$ $(i=0,1,2)$ formally cohesive, we have the long exact sequence
	\begin{displaymath}
	\t^n \cF \xrightarrow{((f_1)_{\ast}, (f_2)_{\ast} )} \t^n \cF_1 \oplus \t^n \cF_2 \xrightarrow{ (p_1)_{\ast} - (p_2)_{\ast} } \t^n \cF_0 \to \t^{n+1} \cF \to \dots.
	\end{displaymath}
\end{lem} 

\begin{proof}
	The functor $\cF$ is formally cohesive by \lemref{formally cohesive functors from known ones}.  The equation $\t \cF = \holim_{c\in \cC}  \,  \t \cF_c$ follows immediately from $\cF(k\oplus \DK(C_{\ast})) \simeq \DK(\tau_{\geq0} (\t \cF \otimes C_{\ast}))$ ($C_{\ast} \in \Ch_{\geq 0}(k)$ with $H_{\ast}(C_{\ast})$ finite). 
\end{proof}

\section{Derived deformation functors}

In this section, we will define the derived deformation functors with prescribed local deformation conditions, and study the homotopy of the pro-representing rings.  The main result is \thmref{main theorem}, where we show that \cite[Theorem 14.1]{GV18} holds in our more general setting.  

In Section 3.1, we will introduce the derived universal deformation functor with an emphasis on the center-modified version following \cite[Section 5.4]{GV18}.  In Section 3.2, we will define the derived local deformation problems using the classical framed local deformation rings; this can be thought of as the reverse procedure of \remref{framed derived functor}, where we define the derived framed deformation functor from the unframed one.  In Section 3.3 we will impose local conditions to the derived global deformation functor, and in Section 3.4 we will verify the calculations of \cite[Section 11 and Section 14]{GV18} in our more general setting and then prove \thmref{main theorem}.

\subsection{Derived universal deformation functor}

\subsubsection{Reformulation of $\Def_S$}\label{reformulation}

Let $\bar\rho \colon \Gamma_S \to G(k)$ be a fixed residual representation.  Recall we defined $\Def_S \colon \CNL_{\cO} \to \Set$ by associating $A\in \CNL_{\cO}$ to the set of $\ker(G(A) \to G(k))$-conjugacy classes of continuous liftings $\rho\colon \Gamma_S \to G(A)$ which make the following diagram commute:
\begin{displaymath}
\xymatrix{
	\Gamma_S \ar[r]^{ \rho }\ar[dr]^{\bar\rho} & G(A) \ar[d]\\
	& G(k).
}
\end{displaymath}
It's convenient to work with Artinian local $\cO$-algebras $\Art_{\cO}$ instead of $\CNL_{\cO}$ to avoid the issue of continuity, so we often regard $\Gamma_S$ as the projective limit of finite groups $\Gamma_i$ and restrict $\Def_S$ to $\Art_{\cO}$.

In the following we shall explain the simplicial interpretation of $\Def_S \colon \Art_{\cO} \to \Set$.

Let $\Gpd$ be the category of small groupoids (recall a groupoid is a category such that all homomorphisms between two objects are isomorphisms).  Note that a group $G$ can be regarded as a one point groupoid $\bullet$ with $\End(\bullet)=G$.  One reason for introducing groupoids is that $\Gpd$ is a model category (see \cite[Theorem 6.7]{Str00}), while $\Gp$ is not.  Let's recall a morphism $f\colon G\to H$ of $\Gpd$ is 
\begin{enumerate}
	\item a weak equivalence if it is an equivalence of categories;
	\item a cofibration if it is injective on objects;
	\item a fibration if for all $a\in G$, $b\in H$ and $h\colon f(a) \to b$ there exists $g\colon a\to a^{\prime}$ such that $f(a^{\prime}) = b$ and $f(g)=h$.
\end{enumerate}
Moreover, the empty groupoid is the initial object and the unit groupoid consisting in a unique object with a unique isomorphism is the final object, every object of $\Gpd$ is both cofibrant and fibrant, and the homotopy category $\Ho(\Gpd)$ is the quotient category of $\Gpd$ modulo natural isomorphisms.  By regarding a group $G$ as a one point groupoid, the functor $\Gp \to \Ho(\Gpd)$ so obtained has the effect of modulo conjugations, so, for any finite group $\Gamma_i$, we have 
\begin{displaymath}
\Hom_{\Gp}(\Gamma_i, G(A))/G^{\ad}(A) \cong \Hom_{\Ho(\Gpd)}(\Gamma_i, G(A)).
\end{displaymath}

Let $\Cat$ be the category of small categories.  Let's recall the nerve construction for $\Cat$ and $\Gpd$; it's an application of \lemref{adjoint}:

\begin{expl}\label{nerve}
	\begin{enumerate}
		\item Let $\bm{\Delta} \to \Cat$ be the functor defined by regarding $[n]$ as a posetal category: its objects are $0,1,\ldots n$ and $\Hom_{[n]}(k,\ell)$ has at most one element, and is non-empty if and only if $k\leq \ell$.  We write $P\colon \sSet \to \Cat$ and $B\colon \Cat \to \sSet$ for the associated left adjoint and right adjoint respectively.  The functor $B$ is called the nerve functor.  The simplicial set $B\cC=(X_n)$ is defined by sets $X_n\subset \mathrm{Ob}(\cC)^{[n]}$ of $(n+1)$-tuples $(C_0,\ldots,C_n)$ of objects of $\cC$ with morphisms $C_k\to C_\ell$ when $k\leq \ell$, which are compatible when $n$ varies; it is a fibrant simplicial set if and only if $\cC\in \Gpd$ (see \cite[Lemma \RNum{1}.3.5]{GJ09}).  In a word, for $B\cC$ to be fibrant, it must have the extension property with respect to inclusions of horns in $\Delta^n$ ($\forall n\geq 1$).  For $n=2$, it amounts to saying that all homomorphisms in $\cC$ are invertible; for $n>2$, the extension condition is automatic (details in the reference above).  For $\cC\in \Cat$, we have $PB\cC\cong \cC$, so $\Hom_{\Cat}(\cC, \cD) \cong \Hom_{\sSet}(B\cC, B\cD)$ ($\forall \cC,\cD\in \Cat$).  Note that $B(\cC\times [1]) \cong B\cC \times \Delta[1]$ (product is taken degreewise); in consequence, when $\cC\in \Cat$ and $\cD \in \Gpd$, two functors $f,g\colon \cC\to\cD$ are naturally isomorphic if and only if $Bf$ and $Bg$ are homotopic.
		\item As a corollary of (1), we have $\Hom_{\Gpd}(GPX, H) \cong \Hom_{\sSet}(X, BH)$ for $X\in \sSet$ and $H \in \Gpd$, where $GPX$ is the free groupoid associated to $PX$.  We remark that $GPX$ and $\pi_1|X|$ (the fundamental groupoid of the geometric realization) are isomorphic in $\Ho(\Gpd)$ (see \cite[Theorem \RNum{3}.1.1]{GJ09}). 
	\end{enumerate}
\end{expl}

\begin{lem} 
	The nerve functor $B\colon \Gpd \to \sSet$ is fully faithful and Kan-valued.  Moreover, it is right Quillen.
\end{lem}

\begin{proof}
	For the first statement, we know by the above example that $\Hom_{\Cat}(\cC, \cD) \cong \Hom_{\sSet}(B\cC, B\cD)$ ($\forall \cC,\cD\in \Cat$) and $B\cC$ is fibrant for a groupoid $\cC$.
	
	For the second statement, note that $B$ obviously preserves weak equivalences; moreover, by definition, $Bf\colon BG\to BH$ is a fibration if and only if it has the right lifting property with respect to inclusions of horns in $\Delta^n$, $\forall n\geq 1$ (see \cite[page 10]{GJ09}). For $n=1$ this means exactly that $f$ is a fibration, while for $n\geq2$ it's automatic (see the proof of \cite[Lemma \RNum{1}.3.5]{GJ09}).
\end{proof}

For convenience, for $\Gamma_S = \varprojlim \Gamma_i$, we understand $B\Gamma_S$ as the pro-simplicial set $(B\Gamma_i)$ (here each $\Gamma_i$ is regarded as the one object groupoid $\bullet$ such that $\End(\bullet)=\Gamma_i$).  For $A\in \Art_{\cO}$, by applying the above lemma and then passing to homotopy categories, we get  
\begin{align*}
\Hom_{\Gp}(\Gamma_i, G(A))/G^{\ad}(A) &\cong \Hom_{\Ho(\Gpd)}(\Gamma_i, G(A)) \\
&\cong \Hom_{\Ho(\sSet)}(B\Gamma_i, BG(A))\\
&\cong \pi_0\sHom_{\sSet}(B\Gamma_i, BG(A)).
\end{align*}
Passing to the limit, $\pi_0\sHom_{\sSet}(B\Gamma_S, BG(A))$ is isomorphic to the set of $G^{\ad}(A)$-conjugacy classes of continuous maps from $\Gamma_S$ to $G(A)$.

We shall consider the deformations of $\bar\rho$, so it's natural to work with the overcategory $\sSet/_{BG(k)}$.  It is also a simplicial model category: the cofibrations, fibrations, weak equivalences and tensor products are those of $\sSet$ (see \cite[Lemma \RNum{2}.2.4]{GJ09} for the only non-trivial part of the statement).  Note that $\bar\rho \colon \Gamma_S \to G(k)$ induces a map $B\Gamma_S \to BG(k)$, which makes $B\Gamma_S$ a pro-object of $\sSet/_{BG(k)}$.  Similar to preceding discussions, we have 
\begin{equation*}
\Def_S(A) \cong \Hom_{\Ho(\sSet/_{BG(k)})}(B\Gamma_S, BG(A)) \cong \pi_0\sHom_{\sSet/_{BG(k)}}(B\Gamma_S, BG(A)) 
\end{equation*}
for $A\in \Art_{\cO}$.  Note that $\sHom_{\sSet/_{BG(k)}}(B\Gamma_S, BG(A))$ is the fiber over $\bar\rho$ of the fibration map $$\sHom_{\sSet}(B\Gamma_S, BG(A)) \to \sHom_{\sSet}(B\Gamma_S, BG(k)),$$ so it is actually the homotopy fiber (see \cite[Theorem 13.1.13 and Proposition 13.4.6]{Hir09}).

\begin{rem}\label{framed derived functor}
	The same argument gives a simplicial interpretation of the framed universal deformation functor $\Def_S^{\square}$.  Let $\Gpd_{\ast}$ be the category of based groupoids ({\it i.e.}, the under category ${}_{\ast} \backslash \Gpd$).  Now one has  
	\begin{displaymath}
	\Hom_{\Gp}(\Gamma_i, G(A)) \cong \Hom_{\Ho(\Gpd_{\ast})}(\Gamma_i, G(A)).
	\end{displaymath}
	
	We regard $B\Gamma_S$ as a pro-object of the over and under category ${}_{\ast} \backslash \sSet /_{BG(k)}$ under $\bar\rho \colon \Gamma_S \to G(k)$ (note ${}_{\ast} \backslash \sSet /_{BG(k)}$ is also a simplicial model category: the cofibrations, fibrations, weak equivalences are those of $\sSet$, and the tensor product of $X\in {}_{\ast} \backslash \sSet /_{BG(k)}$ and $K\in \sSet$ is the pushout of $\ast \leftarrow \ast \otimes K \rightarrow X\otimes K$).  Proceeding as the unframed case, one gets
	\begin{equation*}
	\Def_S^{\square}(A) \cong \Hom_{\Ho({}_{\ast} \backslash \sSet /_{BG(k)})}(B\Gamma_S, BG(A)) \cong \pi_0\sHom_{{}_{\ast} \backslash \sSet /_{BG(k)}}(B\Gamma_S, BG(A)) 
	\end{equation*}
	for $ A\in \Art_{\cO}$.  
	
	By the description of the tensor product in ${}_{\ast} \backslash \sSet /_{BG(k)}$, one sees that $\sHom_{{}_{\ast} \backslash \sSet /_{BG(k)}}(B\Gamma_S, BG(A))$ is isomorphic to the fiber over the base point of the fibration map
	$$\sHom_{\sSet/_{BG(k)}}(B\Gamma_S, BG(A)) \to \sHom_{\sSet/_{BG(k)}}(\ast, BG(A)).$$
	In other words, one has the homotopy pullback square
	\begin{displaymath}
	\xymatrix{
		\sHom_{{}_{\ast} \backslash \sSet /_{BG(k)}}(B\Gamma_S, BG(A)) \ar[r]\ar[d] & \ast \ar[d] \\
		\sHom_{\sSet/_{BG(k)}}(B\Gamma_S, BG(A)) \ar[r] & \sHom_{\sSet/_{BG(k)}}(\ast, BG(A)).
	}
	\end{displaymath}
\end{rem}

\subsubsection{Derived universal deformation functor}

Let's extend the functor $\sHom_{\sSet/_{BG(k)}}(B\Gamma_S, BG(-))$ to the category $\sArt$ (see \defref{simplicial Artinian algebra}).

Define $\cO_{N_{\bullet}G} \in \Alg_{\cO}^{\bm{\Delta}}$ ({\it i.e.}, a functor $\bm{\Delta}\to \Alg_{\cO}$, also called a cosimplicial object in $\Alg_{\cO}$) as follows: in codegree $p$ we have $\cO_{N_p G} = \cO_{G}^{\otimes p}$, and the coface and codegeneracy maps are induced from the comultiplication and the coidentity of the Hopf algebra $\cO_{G}$ respectively.  Then for $A\in \Alg_{\cO}$, the nerve $BG(A)$ is exactly $\Hom_{\Alg_{\cO}}(\cO_{N_{\bullet}G}, A)$, with face and degeneracy maps induced by the coface and codegeneracy maps in $\cO_{N_{\bullet}G}$.  When $A\in \sCRO$, the na\"\i ve analogy is the diagonal of the bisimplicial set $([p],[q]) \mapsto \Hom_{\Alg_{\cO}}(\cO_{N_p G}, A_{q})$ (recall that the diagonal of a bisimplicial set is a simplicial set model for its geometric realization).  However, we need to make some modifications using cofibrant replacements to ensure the homotopy invariance.  Recall that $\sCR$ is cofibrantly generated, so there is a functorial factorization $\cO \hookrightarrow c(A) \stackrel{\sim}{\twoheadrightarrow} A$ for $A\in \sCRO$.

\begin{de}
	\begin{enumerate}
		\item For $A\in \sCRO$, we define $\Bi(A)$ to be the bisimplicial set $$([p],[q])\mapsto  \Hom_{\sCRO}(c(\cO_{N_p G}),A^{\Delta[q]}),$$ with face and degeneracy maps induced by the coface and codegeneracy maps in $\cO_{N_{\bullet}G}$ and the face and degeneracy maps in $A^{\Delta[\bullet]}$.
		\item The diagonal $\diag \Bi(A)$ is the simplicial set induced from the diagonal embedding 
		$$\bm{\Delta}^{\op} \to \bm{\Delta}^{\op} \times \bm{\Delta}^{\op} \xrightarrow{\Bi(A)} \Set.$$ 
	\end{enumerate}
\end{de}

When $A$ is an $\cO$-algebra regarded as a constant object in $\sCRO$, we have $$\Bi(A)_{p,q} = \Hom_{\sCRO}(c(\cO_{N_p G}),A^{\Delta[q]}) \cong \Hom_{\Alg_{\cO}}(\cO_{N_p G},A),$$ where the latter isomorphism is because the constant embedding functor is right adjoint to $\pi_0\colon \sCRO \to \Alg_{\cO}$.  Hence $\Bi(A)$ is just a disjoint union of copies of $BG(A)$ in index $q$.  In particular, for $A\in \sArt$ there is a natural map $\Bi(A)_{\bullet,q} \to BG(k)$ for each $q\geq 0$, so we may regard $\Bi(A) \in (\sSet/_{BG(k)})^{\bm{\Delta}^{\op}}$ via the association $[q]\mapsto \Bi(A)_{\bullet, q}$.  Recall that any morphism $X\to Y$ in $\sSet$ admits a functorial factorization $$X \stackrel{\sim}{\hookrightarrow} \widetilde{X} \twoheadrightarrow Y$$ into a trivial cofibration and a fibration.

\begin{de}\label{BG}
	For $A\in \sArt$, the simplicial set $\BG(A)$ is defined by the functorial trivial cofibration-fibration factorization $\diag \Bi(A) \stackrel{\sim}{\hookrightarrow} \BG(A) \twoheadrightarrow BG(k)$.
\end{de}

It's clear that $\BG\colon \sArt \to \sSet/_{BG(k)}$ defines a functor.  If $A\in \Art_{\cO}$ is a constant simplicial ring, then $\diag \Bi(A) = BG(A) \twoheadrightarrow BG(k)$ is a fibration, and hence $BG(A)$ is a strong deformation retract of $\BG(A)$ in $\sSet/_{BG(k)}$ (see \cite[Definition 7.6.10]{Hir09}).  In particular, these two are indistinguishable in our applications.

The following lemma explains the reason for taking cofibrant replacements of $\cO_{N_p G}$:

\begin{lem}
	If $A \to B$ is a weak equivalence, then so is $\BG(A) \to \BG(B)$.
\end{lem}

\begin{proof}
	If $A \to B$ is a weak equivalence, then 
	$$\sHom_{\sCRO}(c(\cO_{N_p G}),A) \to \sHom_{\sCRO}(c(\cO_{N_p G}),B)$$ 
	is a weak equivalence for each $p\geq 0$, so are $\diag \Bi(A) \to \diag \Bi(B)$ (see \cite[Theorem 15.11.11]{Hir09}) and $\BG(A) \to \BG(B)$.
\end{proof}

\begin{de}
	\begin{enumerate}
		\item The derived universal deformation functor $s\Def_S \colon \sArt \to \sSet$ is defined by 
		$$s\Def_S(A)=\sHom_{\sSet/_{BG(k)}}(B\Gamma_S, \BG(A)). $$
		\item The derived universal framed deformation functor $s\Def_S^{\square} \colon \sArt \to \sSet$ is defined by 
		$$s\Def_S^{\square}(A)=\hofib_{\ast} (s\Def_S (A) \to \sHom_{\sSet/_{BG(k)}}(\ast, \BG(A))). $$
	\end{enumerate}
\end{de}

Note $s\Def_S(A)$ can be defined alternatively as 
$$\hofib_{\bar\rho} (\sHom_{\sSet}(B\Gamma_S, \BG(A)) \to \sHom_{\sSet}(B\Gamma_S, BG(k))).$$

The following proposition summarizes the properties of the derived functors:

\begin{pro}
	The functors $s\Def_S$ and $s\Def_S^{\square}$ are formally cohesive.  
\end{pro}

\begin{proof}
	We first verify three conditions in the above definition for $s\Def_S$:
	\begin{enumerate}
		\item If $A \to B$ is a weak equivalence, then $\BG(A)\to \BG(B)$ is a weak equivalence between fibrant objects in $\sSet/_{BG(k)}$, so 
		$$\sHom_{\sSet/_{BG(k)}}(B\Gamma_S, \BG(A)) \to \sHom_{\sSet/_{BG(k)}}(B\Gamma_S, \BG(B))$$ 
		is also a weak equivalence.
		\item By \cite[Lemma 4.31]{GV18}, to prove 
		\begin{displaymath}
		\xymatrix{
			\BG(A) \ar[r]\ar[d] & \BG(B) \ar[d] \\
			\BG(C) \ar[r] & \BG(D)}
		\end{displaymath}
		is a homotopy pullback square (one can regard this diagram in either $\sSet/_{BG(k)}$ or $\sSet$), it suffices to check:
		\begin{enumerate}
			\item the functor $\Omega\BG\colon \sArt \to \sSet$ preserves homotopy pullbacks, and 
			\item $\pi_1\BG(C) \to \pi_1\BG(D)$ is surjective whenever $C \to D$ is degreewise surjective.
		\end{enumerate}
		Part (a) follows from \cite[Lemma 5.2]{GV18}, and part (b) follows from \cite[Corollary 5.3]{GV18}.  
		
		Then since $\BG$ is homotopy invariant and take fibrant values in $\sSet/_{BG(k)}$, we can apply \lemref{formally cohesive functors from known ones} to deduce that $s\Def_S = \sHom_{\sSet/_{BG(k)}}(B\Gamma_S, \BG(-))$ preserves homotopy pullback squares.
		\item It's clear that $s\Def_S(k)$ is contractible.
	\end{enumerate}
	
	The same argument applies for $A \to \sHom_{\sSet/_{BG(k)}}(\ast, \BG(A))$.  So $s\Def_S^{\square}$ is formally cohesive as it is the homotopy pullback of formally cohesive functors.
\end{proof}

Now it's clear that $s\Def_S$ and $s\Def_S^{\square}$ are indeed generalizations of $\Def_S$ and $\Def_S^{\square}$:

\begin{pro}
	When $A$ is homotopy discrete ({\it i.e.}, $A$ is weakly equivalent to $\pi_0 A$), we have $\pi_0s\Def_S (A) \cong \Def_S(\pi_0 A)$ and $\pi_0s\Def_S^{\square} (A) \cong \Def_S^{\square} (\pi_0 A)$.
\end{pro}

\begin{proof}
	By the formal cohesiveness, we may suppose $A$ is a constant simplicial ring.  Then since $BG(A)$ is a strong deformation retract of $\BG(A)$ in $\sSet/_{BG(k)}$, the proposition follows from the discussions in Section \ref{reformulation}.
\end{proof}

It's natural to ask if the functors $s\Def_S$ and $s\Def_S^{\square}$ are pro-representable, and for this one has to calculate their tangent complexes.  From now on, we will use calligraphic letters for the pro-representing rings of derived deformation functors to distinguish them from the classical representing rings.

\begin{lem}\label{tangent complex of universal derived deformation functor}
	\begin{enumerate}
		\item We have $\t^i s\Def_S = H^{i+1}(\Gamma_S, \g_k)$ for all $i \in \Z$. 
		\item We have $\t^i s\Def_S^{\square} = \left\{\begin{array}{l}{  0 \text{ if } i<0;  } \\ { Z^1(\Gamma_S, \g_k) \text{ if } i=0; }  \\ {  H^{i+1}(\Gamma_S, \g_k) \text{ if } i>0.  }\end{array}\right.$
	\end{enumerate}
\end{lem}

\begin{proof}
	\begin{enumerate}
		\item See \cite[Lemma 5.10]{GV18}.  Here we give a slightly different approach.  
		
		Without loss of generality, we temporarily forget the pro-issue on $X = B\Gamma_S$.  Then by \cite[Proposition 18.9.2]{Hir09}, $X$ is weakly equivalent to $\hocolim_{(\bm{\Delta}X)^\op} \, \ast$ (i.e., the homotopy colimit of the single-point simplicial set indexed by $(\bm{\Delta}X)^\op$), and hence 
		$$\sHom_{\sSet}(X, \BG(k \oplus k[n])) \simeq \holim_{\bm{\Delta}X}  \, \BG(k \oplus k[n]).$$
		Since homotopy limits commute with homotopy pullbacks, we deduce 
		$$s\Def_S(k \oplus k[n]) \simeq \holim_{\bm{\Delta}X}  \, \sHom_{\sSet/_{BG(k)}}(\ast, \BG(k \oplus k[n])).$$
		So $\t (s\Def_S)$ is the homotopy limit indexed by $\bm{\Delta}X$ of $\t (\sHom_{\sSet/_{BG(k)}}(\ast, \BG(-)))$.  
		The homotopy groups of $\hofib_{\ast}(\BG(k \oplus k[j]) \to BG(k))$ are trivial except at degree $j+1$, where it is $\g_k$ (see \cite[Lemma 5.5]{GV18}), so $\t (\sHom_{\sSet/_{BG(k)}}(\ast, \BG(-)))$ is concentrated on degree $-1$, where it is $\g_k$.  The $\bm{\Delta}X$-diagram of complexes on $X$ forms a cohomological coefficient system in the sense of \cite[Page 28]{GM13}, or local system in the sense of \cite[Definition 4.34]{GV18}, and the $\pi_1 (X,\ast)$-action on $\g_k$ is exacly the adjoint action.  
		
		By shifting (co)degrees $i \mapsto i+1$, it suffices to calculate $\holim_{\bm{\Delta}X} \, \g_k$ where $\g_k$ is the cochain complex concentrated on degree $0$.  By \cite[Lemma 18.9.1]{Hir09}, $\holim \, \g_k$ is naturally isomorphic to $\holim_{\bm{\Delta}} \, Z$ where $Z$ is the cosimplicial object in $\Ch^{\geq 0}(k)$ whose codegree $[n]$ term is $\prod_{\sigma \in X_n} \g_k$.  The coface maps of $Z$ can be described as follows: 
		
		The $k[\Gamma_S]$-module $\g_k$ defines a functor $D$ from the one-object groupoid $\bullet$ with $\End(\bullet)=\Gamma_S$ to $\Ch^{\geq 0}(k)$, such that $D(\bullet) = \g_k$, and $D(\Gamma_S)$ acts on $\g_k$ by the adjoint action.  Then $Z^n$ is $\prod\limits_{i_0 \to \dots \to i_n} D(i_n)$ (all $i_k$'s are equal to the object $\bullet$ here, but keeping the difference helps to clarify the process).  Let $d_k$ be the $k$-th face map from $\Gamma_S^{n+1}$ to $\Gamma_S^n$, in other words, $d_k$ maps $(i_0 \to \dots \to i_{n+1})$ to $(j_0 \to \dots \to j_{n})$ by "covering up" $i_k$.  Then the corresponding $D(j_n) \to D(i_{n+1})$ is the identity map if $k\neq n+1$, and is $D(i_n \to i_{n+1})$ if $k=n+1$.  
		
		By \cite[Proposition 19.10]{Dug08}, $\holim_{\bm{\Delta}} \, Z$ is quasi-isomorphic to the total complex of the alternating double complex defined by $Z$.  Since each $Z^n$ is concentrated on degree $0$, the total complex is simply 
		$$\dots \to \prod\limits_{\Gamma_S^n} \g_k \to \prod\limits_{\Gamma_S^{n+1}} \g_k \to \dots$$
		and the alternating sum $\prod\limits_{\Gamma_S^n} \g_k \to \prod\limits_{\Gamma_S^{n+1}} \g_k$ is exactly the one which computes the group cohomology.  We deduce $\holim_{\bm{\Delta}X} \, \g_k \simeq C^{\bullet}(\Gamma_S, \g_k)$, and hence ($+1$ arises from the degree-shifting) $\t^i s\Def_S = H^{i+1}(\Gamma_S, \g_k)$ for all $i \in \Z$. 
		\item From \lemref{tangent complex calculation} and 
		$$s\Def_S^{\square}(A)=\hofib_{\ast} (s\Def_S (A) \to \sHom_{\sSet/_{BG(k)}}(\ast, \BG(A))),$$ 
		we get the long exact sequence 
		\begin{align*}
		\t^i s\Def_S^{\square}  \to \t^i s\Def_S \to \t^i \sHom_{\sSet/_{BG(k)}}(\ast, \BG(-)) \stackrel{[1]}{\to} \dots.
		\end{align*}
		In the proof of (1), we know $\t^i s\Def_S = H^{i+1}(\Gamma_S, \g_k)$ ($\forall i\in \Z$) and 
		\begin{displaymath}
		\t^i \sHom_{\sSet/_{BG(k)}}(\ast, \BG(-)) = \left\{\begin{array}{l}{  \g_k \text{ if } i=-1;  } \\ { 0 \text{ if } i\neq -1 }.  \end{array}\right.
		\end{displaymath}
		So the conclusion follows from the above long exact sequence; note by \lemref{tangent complex calculation} all maps there are natural ones.
	\end{enumerate}
\end{proof}

By Lurie's criterion \ref{Lurie}, the functor $\t^i s\Def_S^{\square}$ is always pro-representable, while the functor $s\Def_S$ can't be pro-representable unless $H^0(\Gamma_S, \g_k) = 0$.  If $G$ has a nontrivial center $Z$, we need a variant $s\Def_{S,Z}$ of the functor $s\Def_S$, in order to allow pro-representability.

\subsubsection{Modifying the center}\label{modifying the center}

We follow \cite[Section 5.4]{GV18} for this modification.  Define $PG=G/Z$, then the short exact sequence $1\to Z(A)\to G(A)\to PG(A)\to 1$  yields a fibration sequence $BG(A)\to BPG(A)\to B^2Z(A)$.  Indeed, given a simplicial group $H$ and a simplicial sets $X$ with a left $H$-action, we can form the bar construction $N_\ast(\ast,H,X)$ at each simplicial degree (see \cite[Example 3.2.4]{Gil13}), which gives the bisimplicial set $([p],[q]) \mapsto  H_p^q\times X_p=\colon N_q (\ast, H_p, X_p)$.  Consider the action $Z(A)\times G(A)\to G(A)$, and the corresponding simplicial action $N_pZ(A)\times N_pG(A)\to N_pG(A)$ (note that $N_\ast Z(A)$ is a simplicial group because $Z(A)$ is abelian).  We identify for each $p\geq 0$, $$BG(A)_p=N_p(\ast,\ast,N_pG(A)),$$ $$BPG(A)_p=N_p(\ast,N_pZ(A),N_pG(A)),$$ and we put $$B^2Z(A)_p=N_p(\ast,N_pZ(A),\ast)$$ (with diagonal face and degeneracy maps).  The desired fibration is given by the canonical morphisms of simplicial sets which in degree $p$ are: 
$$N_p(\ast,\ast,N_pG(A))\to N_p(\ast,N_pZ(A),N_pG(A))\to N_p(\ast,N_pZ(A),\ast).$$  

The functor $s\Def_{S,Z} \colon \sArt \to \sSet$ is defined by the homotopy pullback square (here the base maps are those induced from $BG(k) \to BPG(k) \to B^2Z(k)$) 
\begin{displaymath}
\xymatrix{
	s\Def_{S,Z}(A) \ar[r]\ar[d] & \sHom_{\sSet/_{B^2Z(k)}}(\ast,B^2Z(A)) \ar[d] \\
	\sHom_{\sSet/_{BPG(k)}}(B\Gamma_S,\BPG(A)) \ar[r] & \sHom_{\sSet/_{B^2Z(k)}}(B\Gamma_S,B^2Z(A)).}
\end{displaymath}
Then $s\Def_{S,Z}$ is formally cohesive becasue it is the homotopy pullback of formally cohesive functors.  Observe that $s\Def_{S,Z}$ and $s\Def_S$ coincide when $Z$ is trivial.  

\begin{rem}\label{modifying the center functorial}
	Note the construction $s\Def_{S,Z}$ is functorial both in $\Gamma_S$ and $G$.  
\end{rem}

Consider the diagram 
\begin{displaymath}
\xymatrix{
	s\Def_S(A) \ar[r]\ar[d] &  \ast \ar[d] \\
	s\Def_{S,Z}(A) \ar[r]\ar[d] & \sHom_{\sSet/_{B^2Z(k)}}(\ast,B^2Z(A)) \ar[d] \\
	\sHom_{\sSet/_{BPG(k)}}(B\Gamma_S,\BPG(A)) \ar[r] & \sHom_{\sSet/_{B^2Z(k)}}(B\Gamma_S,B^2Z(A)).}
\end{displaymath}
By above discussions, the lower square and the combined square are homotopy fiber squares, so is the upper square (see \cite[Proposition 13.3.15]{Hir09}).  Now we can calculate the tangent complex of $s\Def_{S,Z}$.

\begin{lem}
	We have 
	$$\t^i s\Def_{S,Z} = \left\{\begin{array}{l}{  H^0(\Gamma_S, \g_k)/\z_k \text{ if } i=-1;  } \\ {  H^i(\Gamma_S, \g_k) \text{ otherwise} }.  \end{array}\right.$$
\end{lem}

\begin{proof}
	Using the above homotopy fiber square, the proof is similar to \lemref{tangent complex of universal derived deformation functor}.
\end{proof}

Since we've made the assumption $H^0(\Gamma_S, \g_k) = \z_k$, the functor $s\Def_{S,Z}$ is pro-representatble.

\begin{lem}\label{center fiber sequence}
	$s\Def_{S,Z}$ fits into the fiber sequence 
	$$s\Def_S^{\square}(A) \to s\Def_{S,Z} (A) \to \sHom_{\sSet/_{BPG(k)}}(\ast, \BPG(A)). $$	
\end{lem}

\begin{proof}
	Consider the diagram 
	\begin{displaymath}
	\xymatrix{
		s\Def_S^{\square}(A) \ar[r]\ar[d] & s\Def_S(A) \ar[r]\ar[d] & s\Def_{S,Z}(A) \ar[d] \\
		\ast \ar[r] & \sHom_{\sSet/_{BG(k)}}(\ast, \BG(A)) \ar[r]\ar[d] & \sHom_{\sSet/_{BPG(k)}}(\ast, \BPG(A)) \ar[d] \\
		& \ast \ar[r] &\sHom_{\sSet/_{B^2Z(k)}}(\ast, B^2Z(A))	.}
	\end{displaymath}
	It suffices to apply \cite[Proposition 13.3.15]{Hir09} twice.  Since the right composed square and the lower square are homotopy fiber squares, so is the upper right square.  Then since the left square is also a homotopy fiber square, we deduce that the upper composed square is a homotopy fiber squares.  
\end{proof}

Similar results of this section hold for the derived universal (framed) deformation functors for $\Gamma_v \to G(k)$.  In this case we just replace the subscript $S$ by $v$ in our notations.  Note even after modifying the center, the functor $s\Def_{v,Z}$ is generally not pro-representable, as generally $H^0 (\Gamma_v, \g_k) \neq \z_k$.

\subsection{Derived local deformation problem}

Let $v$ be a finite place of $F$.  Following \cite[Definition 9.1]{GV18}, a derived local deformation problem at $v$ means a functor $\sArt \to \sSet$ equipped with a natural transformation to $s\Def_{v,Z} \colon \sArt \to \sSet$ (note the center-modification here).  Let $\cD_v$ be a local deformation problem and let $R_v$ be the framed deformation ring for $\cD_v$ (so $R_v$ is a quotient of $R_v^{\square}$).  It's natural to try to associate a derived local deformation problem to $\cD_v$.  

Note the conjugation action of $\ker(G(A) \to G(k))$ on a lifting $\Gamma_v \to G(A)$ together with the functorial cofibrant replacement $c$ induce a cosimplicial object $[p] \mapsto c(R_v \otimes \cO_{N_p G}) \in \sCROk$.  To take into account the continuity, we regard $R_v$ as a pro-Artinian object in the following.

\begin{de}\label{derived local problem}
	Associated to $\cD_v$, we define 
	\begin{enumerate}
		\item $s\cD_v^{\square} \colon \sArt \to \sSet$ to be the functor $A \mapsto \sHom_{\sCROk}(c(R_v), A)$;
		\item $s\cD_v \colon \sArt \to \sSet$ to be the functor which sends $A \in \sArt$ to the fixed fibrant replacement in $\sSet /_{BG(k)}$ of the diagonal of $([p],[q])\mapsto  \Hom_{\sCROk}(c(R_v \otimes \cO_{N_p G} ),A^{\Delta[q]})$.
	\end{enumerate}
\end{de}

The definition of $s\cD_v$ is inspired by the simplicial bar construction (one may compare with \cite[Lemma 5.7]{GV18}).  The natural $\bm{\Delta}$-equivariant map $c(\cO_{N_{\bullet} G}) \to c(R_v \otimes \cO_{N_{\bullet} G}) \to c(R_v)$ induces 
$$s\cD_v^{\square}(A) \to s\cD_v(A) \to \sHom_{\sSet/_{BG(k)}}(\ast, \BG(A))$$
for $A \in \sArt$, which is a fibration sequence by \cite[Lemma 4.6.6]{Lan15}.  Using the long exact sequence for homotopy groups, one sees that $s\cD_v$ preserves homotopy pullbacks, since this is the case for $s\cD_v^{\square}$ and $\sHom_{\sSet/_{BG(k)}}(\ast, \BG(-))$.  Then we deduce 

\begin{lem}
	$s\cD_v$ is formally cohesive.  
\end{lem}

Now we construct the natural transformation $s\cD_v \to s\Def_v$.  

\begin{pro}
	There is a natural transformation $s\cD_v \to s\Def_v$ making the diagram
	\begin{displaymath}
	\xymatrix{
		s\cD_v^{\square}(A) \ar[r]\ar[d] & s\cD_v(A) \ar[r]\ar[d] &\sHom_{\sSet/_{BG(k)}}(\ast, \BG(A)) \ar@{=}[d] \\
		s\Def_v^{\square}(A) \ar[r] &  s\Def_v (A) \ar[r] &\sHom_{\sSet/_{BG(k)}}(\ast, \BG(A)).}
	\end{displaymath}
	commutative up to weak equivalence.  Here the first vertical arrow is induced from $\cR_v^{\square} \to \pi_0 \cR_v^{\square} \to R_v$.
\end{pro}

\begin{rem}
	When the representing ring $\cR_v^{\square}$ for $s\Def_v^{\square}$ is homotopy discrete, the map $s\cD_v(A) \to s\Def_v (A)$ is the natural one induced from the quotient map $R_v^{\square} \to R_v$.  Note the homotopy discreteness of $\cR_v^{\square}$ is equivalent to the conjecture below \cite[(1.5)]{GV18} which says $R_v^{\square}$ is a complete intersection ring of expected dimension.  Here we don't need $\cR_v^{\square}$ to be homotopy discrete, which illustrates in a certain sense the comment of $\loccit$ that one of the advantages of the derived deformation ring is to circumvent the conjecture mentioned above.
\end{rem}

\begin{proof}
	Fix $A \in \sArt$.  We write $Z = s\Def_v(A)$ and write $\bX$ for the bisimplicial set $([p],[q])\mapsto  \Hom_{\sCROk}(c(R_v \otimes \cO_{N_p G} ),A^{\Delta[q]})$.  Note $\bX$ can be viewed as a simplicial object in $\sSet$ through $[p]\mapsto \bX_p = \sHom_{\sCROk}(c(R_v \otimes \cO_{N_p G} ),A)$.  By \cite[Theorem 15.11.6]{Hir09}, $\diag \bX$ is naturally isomorphic to the realization $|\bX|$, or in other words the coend $\bX \otimes_{\bm{\Delta}^{\op}} \Delta$ where $\Delta$ is the cosimplicial standard simplex.  
	
	So it suffices to construct a $\sSet$-morphism $\bX \otimes_{\bm{\Delta}^{\op}} \Delta \to s\Def_v(A)$, or equivalently a system of $\sSet$-morphisms $\Delta^n \to \sHom_{\sSet}(\bX_n, Z)$ which is $\bm{\Delta}$-compatible in $[n]$.  Given $[n] \in \bm{\Delta}$, we construct $\Delta^n_k \to \Hom_{\sSet}(\bX_n \otimes \Delta^k, Z)$ by induction on $k$: for $k=0$ a map $[0] \to [n]$ gives naturally $\bX_n \to \bX_0 \to Z$ where the second arrow is induced from $s\Def_v^{\square}(A) \to s\Def_v(A)$; for $k>0$, each of the $(k+1)$ maps $\bX_k \to \bX_0 \to Z$ factors through $s\Def_v^{\square}(A) \to s\Def_v(A)$, so we can choose a morphism $\bX_k \otimes \Delta^k \to Z$ such that for $[l] \to [k]$ with $l<k$ it is compatible with $\bX_k \otimes \Delta^l \to \bX_l \otimes \Delta^l \to Z$ via the embedding $\Delta^l \to \Delta^k$, and $\bX_n \otimes \Delta^k \to Z$ associated to $[k] \to [n]$ is the composition $\bX_n \otimes \Delta^k \to \bX_k \otimes \Delta^k \to Z$.  Thus we get a $\sSet$-morphism $\Delta^n \to \sHom_{\sSet}(\bX_n, Z)$, and this construction is clearly $\bm{\Delta}$-compatible in $[n]$.  It's direct to check that the map $s\cD_v(A) \to s\Def_v(A)$ make the above diagram commutative up to weak equivalence.
\end{proof}

We will always take the center-modification into account.  For this it suffices to replace $G$ by $PG=G/Z$ in \defref{derived local problem}, and henceforth we will instead write $s\cD_v$ for the fibrant replacement of the diagonal of $([p],[q])\mapsto  \Hom_{\sCROk}(c(R_v \otimes \cO_{N_p (PG)} ),A^{\Delta[q]})$ to simplify our notations.  Analogous to the above proposition and using \lemref{center fiber sequence}, we have the following:

\begin{cor}
	There is a natural diagram 
	\begin{displaymath}
	\xymatrix{
		s\cD_v^{\square}(A) \ar[r]\ar[d] & s\cD_v(A) \ar[r]\ar[d] &\sHom_{\sSet/_{BPG(k)}}(\ast, \BPG(A)) \ar@{=}[d] \\
		s\Def_v^{\square}(A) \ar[r] &  s\Def_{v,Z} (A) \ar[r] &\sHom_{\sSet/_{BPG(k)}}(\ast, \BPG(A)).}
	\end{displaymath}
	which is commutative up to weak equivalence and whose rows are fiber sequences.
\end{cor}

\begin{rem}
	In some cases we can define the derived local deformation problem more arithmetically.  For the unramified condition, see the example on \cite[Page 91]{GV18}.  For the (nearly) ordinary condition, one can also define the derived local deformation functor directly by replacing the role of $G$ by its Borel $B$, and under the regularity and dual regularity conditions (see \cite[Propostion 6.2 and Propostion 6.3]{Til96}), this definition coincides with the one using the framed ring (see discussions after \cite[Definition 2.13]{CT20}).
\end{rem}

\begin{lem}
	When $A \in \Alg_{\cO}$ is regarded as a constant simplicial ring, $\pi_0 s\cD_v (A)$ is isomorphic to $\cD_v(A)$. 
\end{lem}

\begin{proof}
	For $A \in \Alg_{\cO}$, we have (the canonical base point is omitted for brevity)
	$$\pi_1 \sHom_{\sSet/_{BPG(k)}}(\ast, BPG(A)) \cong \ker(PG(A) \to PG(k)),$$ 
	and it acts by conjugation on $\pi_0 s\cD_v^{\square}(A) \cong \cD_v(A)$. 
	
	Moreover, We have the sequence of maps 
	$$\pi_1 \sHom_{\sSet/_{BPG(k)}}(\ast, BPG(A)) \to \pi_0 s\cD_v^{\square}(A) \to \pi_0  s\cD_v(A)$$
	such that $\pi_0 s\cD_v^{\square}(A) \to \pi_0  s\cD_v(A)$ is surjective and two elements of $\pi_0 s\cD_v^{\square}(A)$ have the same image if and only if they are in the same orbit for the $\pi_1 \sHom_{\sSet/_{BPG(k)}}(\ast, BPG(A))$-action.  The conclusion follows easily.
\end{proof}

Recall $R_v$ is said to be formally smooth if it's a power series ring over $\cO$.  

\begin{lem}\label{local tangent complex}
	Suppose $R_v$ is formally smooth, then we have 
	$$\t^i s\cD_v = \left\{\begin{array}{l}{  H^0(\Gamma_v, \g_k)/\z_k \text{ if } i=-1;  } \\ {  L_v \text{ if } i= 0;  }   \\ {  0 \text{ if } i> 0.  }\end{array}\right.$$
\end{lem}

\begin{proof}
	By \lemref{formally cohesive functors from known ones}, we have the long exact sequence 
	\begin{align*}
	0 \to &\t^{-1}s\cD_v^{\square} \to \t^{-1}s\cD_v \to \t^{-1}\sHom_{\sSet/_{BPG(k)}}(\ast, \BPG(-)) \\
	\to &\t^0s\cD_v^{\square} \to \t^0s\cD_v \to \t^0 \sHom_{\sSet/_{BPG(k)}}(\ast, \BPG(-)) \\
	\to &\t^1 s\cD_v^{\square} \to \t^1 s\cD_v \to \t^1 \sHom_{\sSet/_{BPG(k)}}(\ast, \BPG(-)) \\
	\to &\dots.
	\end{align*}
	Since $s\cD_v^{\square} = \sHom_{\sCROk}(c(R_v), -)$ and $R_v$ is formally smooth, $\t^i s\cD_v^{\square}$ in concentrated on degree $0$, where it is $\widetilde{L}_v = \Hom_{\CNL_{\cO}}(R_v, k \oplus k[0])$.  On the other hand 
	$$\t^i \sHom_{\sSet/_{BPG(k)}}(\ast, \BPG(A))$$ 
	is $\g_k /\z_k$ concentrated on degree $-1$.  Hence $\t^i s\cD_v$ fits into the exact sequence 
	\begin{align*}
	0 \to &0 \to \t^{-1}s\cD_v \to \g_k/\z_k \\
	\to &\widetilde{L}_v \to \t^0s\cD_v \to 0 \\
	\to &0 \to \t^1 s\cD_v \to 0 \\
	\to &\dots
	\end{align*}
	where all maps are natural, and the conclusion follows.
\end{proof}

\subsubsection{Some local deformation problems}\label{local deformation problems}

We discuss some local deformation problems for $\bar\rho \colon \Gamma_v \to G(k)$ for specific groups used in this thesis and \cite{TU21}.  

\subsubsection{Minimal deformations}

Let $v \in S \backslash S_p$.  We would like to formulate a deformation condition which controls the ramifications and is formally smooth (or liftable in \cite{CHT08}).  

For general linear groups, the minimal conditions are defined in \cite[Section 2.4.4]{CHT08}.  It's noted by \cite{Boo19} that the key feature to define a lifting $\rho$ to be minimal is to require $\rho(\tau)$ to have "the same unipotent structure" as $\bar\rho(\tau)$ (for $\tau\in I_v$).  In \loccit the author reinterpreted the definition of \cite{CHT08} using unipotent orbits, and then defined analogously the minimal conditions for symplectic and orthogonal similitude groups.

Let's illustrate some ideas for $G = \GL_N$.  We say $\bar\rho \colon \Gamma_v \to \GL_N(k)$ is minimal if $\bar\rho(I_v)$ contains a regular unipotent element.  Let $J_N$ be the standard Jordan block of size $N$ (note $J_N$ is regular nilpotent) and $t_v\colon I_v \to \Z_p$ be the character defined by $\frac{\tau(\varpi_v^{1/p^n})}{\varpi_v^{1/p^n}}  = \zeta_{p^n}^{t_v(\tau)}$ (for $n\geq 1$ and $\tau\in I_v$).  Without loss of generality, we can suppose $\bar\rho(\tau)  = \exp(t_v(\tau)J_N)$, and we say a lifting $\rho \colon \Gamma_v \to G(A)$ of $\bar\rho$ is minimal if there exists $g_v \in \ker(\GL_N(A) \to \GL_N(k))$ such that $g_v \rho(\tau) g_v^{-1} = \exp(t_v(\tau)J_N)$.  

We write $\cD_v^{\min}$ for the framed minimal deformation functor at $v$, then by \cite[Lemma 1]{TU21} the representing ring is a power series ring in $N^2$ variables over $\cO$, in other words, $\cD_v^{\min}$ is formally smooth and for $L_v \subseteq H^1(\Gamma_v,\g_k)$ associated to $\cD_v^{\min}$ we have $\dim_k L_v - \dim_k H^0(\Gamma_v,\g_k) = 0$ (see also \cite[Corollary 2.4.21]{CHT08}).  Note the unframed deformation ring doesn't exist unless $p\notdivides q_v^i-1$ for all $1 \leq i \leq N-1$.  

For symplectic and orthogonal similitude groups, the minimal deformation condition is defined in \cite[Chapter 4]{Boo17} using the classification of nilpotent orbits by the Bala-Carter data (see \cite[Definition 4.4.2.1]{Boo17}).  By \cite[Proposition 4.4.2.3]{Boo17}, $\cD_v^{\min}$ is formally smooth and $\dim_k L_v - \dim_k H^0(\Gamma_v,\g_k) = 0$ for these groups.

\subsubsection{Ordinary deformations}

In the ordinary case $G$ is allowed to be arbitrary.  Let $B=TN\subseteq G$ be a Borel subgroup scheme ($T$ is a maximal split torus and $N$ is the unipotent radical of $B$, and all these groups are defined over $\cO$).  Let $\Phi$ be the root system associated to $(G,T)$ and $\Phi^+$ the subset of positive roots associated to $(G,B,T)$.

Let $v \in S_p$.  A representation $\bar\rho \colon \Gamma_v \to G(k)$ is call ordinary if there exists $\bar{g}_v\in G(k)$ such that $\bar{\rho}$ takes values in $\bar{g}_v^{-1}  B(k) \bar{g}_v$.  We require the following regularity and dual regularity conditions:

$(Reg_v)$ for any $\alpha\in\Phi^+$, $\alpha\circ\overline{\chi}_v\neq 1$,  and

$(Reg^\ast_v)$ for any $\alpha\in\Phi^+$, $\alpha\circ\overline{\chi}_v\neq \omega$. 

The framed nearly ordinary deformation functor $\cD_v^{\nord}$ is defined such that $\rho \in \cD_v^{\nord}(A)$ if and only if there exists $g_v\in G(A)$ which lifts $\bar{g}_v$ such that $\rho$ takes values in $g_v^{-1}\cdot B(A)\cdot g_v$. 
Note that this implies that the homomorphism $\chi_{\rho,v}\colon\Gamma_v\to T(A)$ given by $g_v\cdot \rho\cdot g_v^{-1}$ lifts $\overline{\chi}_v$.  A lifting $\rho \in \cD_v^{\nord}(A)$ is called ordinary of weight $\mu$ if after conjugation by $g_v$, the cocharacter $\rho\vert_{I_v}\colon I_v\to T(A)=B(A)/N(A)$ is given (via the Artin reciprocity map $\rec_v$) by $\mu\circ \rec_v^{-1}\colon I_v\to \cO_v^\times\to T(A)$, and we write $\cD_v^{\ord, \mu}$ for the framed ordinary deformation functor of weight $\mu$.  We also define $\cD_v^{\ord}$ to be the framed ordinary deformation functor without fixing the weight $\mu$.  By \cite[Lemma 2]{TU21}, the functors $\cD_v^{\nord}$, $\cD_v^{\ord, \mu}$ and $\cD_v^{\ord}$ are all formally smooth, and one has $\dim_k L_v - \dim_k H^0(\Gamma_v,\g_k) = [F_v : \Q_p] (\dim G - \dim B)$.

\subsubsection{Fontaine-Laffaille deformations}

Let $v \in S_p$ be unramified.  For $G = \GL_N$, we write $\cD_v^{\FL}$ for the framed Fontaine-Laffaille deformation functor ({\it i.e.}, $\rho \in \cD_v^{\FL}(A)$ if there exists a $\phi$-filtered $A$-module $M$ free of rank $N$ over $A$, such that $\rho$ is isomorphic to $V_{\crys}(M)$).  By \cite[Corollary 2.4.3]{CHT08}, $\cD_v^{\FL}$ is formally smooth and one has $\dim_k L_v - \dim_k H^0(\Gamma_v,\g_k) = [F_v : \Q_p] (\dim G - \dim B)$.

For a symplectic or orthogonal similitude group, the Fontaine-Laffaille condition with fixed similitude lifting is defined in \cite[Definition 3.2.1.2]{Boo17}, and when the Fontaine-Laffaille weights are multiplicity-free, \cite[Definition 3.2.1.3]{Boo17} proved that $\cD_v^{\FL}$ is formally smooth with $\dim_k L_v - \dim_k H^0(\Gamma_v,\g_k^{\prime}) = [F_v : \Q_p] (\dim G - \dim B)$.

\subsection{Derived deformation functor with local conditions}

Let $\cS =(S, \{\cD_v\}_{v \in S})$ be a global deformation problem (see \defref{global deformation problem}) and let $\cD_{\cS}$ be the deformation functor of type $\cS$.  

\begin{de}
	The derived deformation functor of type $\cS$ is defined to be the homotopy limit 
	$$s\cD_{\cS} =s\Def_{S,Z} \times^h_{\prod_{v\in S} s\Def_{v,Z}}  \prod_{v\in S} s\cD_v.$$
\end{de}

Since each functor on the right hand side is formally cohesive, so is $s\cD_{\cS}$.

\begin{lem}\label{pi0 sD_S}
	When $A\in \Alg_{\cO}$ is regarded as a constant simplicial ring, we have $\pi_0 s\cD_{\cS}(A) \cong \cD_{\cS}(A)$.
\end{lem}

\begin{proof}
	We fix compatible base points.  Firstly, from the fiber sequence  
	$$s\Def_S^{\square}(A) \to s\Def_{S,Z} (A) \to \sHom_{\sSet/_{BPG(k)}}(\ast, BPG(A))$$
	and $H^0(\Gamma_S,\g_k) = \z_k$, we see $\pi_1 s\Def_{S,Z}(A)$ is trivial.  On the other hand, from the diagram 
	\begin{displaymath}
	\xymatrix{
		s\cD_v^{\square}(A) \ar[r]\ar[d] & s\cD_v(A) \ar[r]\ar[d] &\sHom_{\sSet/_{BPG(k)}}(\ast, BPG(A)) \ar@{=}[d] \\
		s\Def_v^{\square}(A) \ar[r] &  s\Def_{v,Z} (A) \ar[r] &\sHom_{\sSet/_{BPG(k)}}(\ast, BPG(A)), }
	\end{displaymath}
	we deduce that $\pi_1 s\Def_{v,Z}(A)$ doesn't contribute to $\pi_0 s\cD_{\cS}(A)$.  Note every functor defining $s\cD_{\cS}$ has the desired $\pi_0$, so $\pi_0 s\cD_{\cS}(A)$ is the fiber of 
	$$\Def_S(A) \oplus \bigoplus_{v \in S} \cD_v(A) \to \bigoplus_{v \in S} \Def_v(A),$$
	and the conclusion follows.
\end{proof}

From now on we suppose every representing ring $R_v$ for $\cD_v$ is formally smooth.  By \lemref{formally cohesive functors from known ones} and \lemref{local tangent complex}, the tangent complex of $s\cD_{\cS}$ fits into the exact sequence 
\begin{align*}
0 \to &\t^{-1} s\cD_{\cS} \to H^0(\Gamma_S, \g_k)/\z_k \oplus \bigoplus_{v \in S} H^0(\Gamma_v,\g_k)/\z_k \to \bigoplus_{v\in S} H^0(\Gamma_v,\g_k)/\z_k \\
\to &\t^0 s\cD_{\cS} \to H^1(\Gamma_S, \g_k) \oplus \bigoplus_{v \in S} L_v \to \bigoplus_{v\in S} H^1(\Gamma_v,\g_k) \\
\to &\t^1 s\cD_{\cS} \to H^2(\Gamma_S, \g_k) \to \bigoplus_{v\in S} H^2(\Gamma_v,\g_k)  \\
\to &\t^2 s\cD_{\cS} \to 0.
\end{align*}
Hence $\t^{-1} s\cD_{\cS} = 0$ and $s\cD_{\cS}$ is pro-representable, say by $\cR_{\cS}$.

\begin{lem}\label{tangent complex of cR_S}
	$\t^i \cR_{\cS} \cong H_{\cS}^{i+1} (\Gamma_S, \g_k)$ for $i\geq 0$.
\end{lem}

\begin{proof}
	This follows directly by comparing the above exact sequence with the exact sequence 
	\begin{align*}
	0 \to &H^0_{\cS}(\Gamma_S, \g_k) \to H^0(\Gamma_S, \g_k) \to 0 \\
	\to &H^1_{\cS}(\Gamma_S, \g_k) \to H^1(\Gamma_S, \g_k) \to  \bigoplus_{v \in S} H^1(\Gamma_v, \g_k)/L_v \\
	\to &H^2_{\cS}(\Gamma_S, \g_k) \to H^2(\Gamma_S, \g_k) \to \bigoplus_{v \in S} H^2(\Gamma_v, \g_k) \\
	\to &H^3_{\cS}(\Gamma_S, \g_k) \to 0.
	\end{align*}
\end{proof}

By \remref{relative dimension remark}, $\t^i \cR_{\cS}$ is concentrated on degrees $0,1$ when $\bar\rho$ has an enormous image and $\zeta_p \notin F$.

\begin{rem}
	Without the assumption that every $R_v$ is formally smooth, the functor $s\cD_{\cS}$ is still pro-representable, but $\t^i \cR_{\cS} \cong H_{\cS}^{i+1} (\Gamma_S, \g_k)$ no longer holds for $i\geq 1$.  We expect a modified version of $\pi_{\ast}\cR_{\cS} \cong \Tor_{\ast}^{S_{\infty}}(R_{\infty}, \cO)$ holds in this case.  
\end{rem}

\begin{rem}
	Let $\Sigma$ be a non-empty subset of $S$.  It's natural to define the derived $\Sigma$-framed deformation of type $\cS$ (see \cite[Page 112]{ACC+18}) as 
	\begin{displaymath}
	s\cD_{\cS}^{\Sigma} =s\Def_{S,Z} \times^h_{\prod_{v\in S} s\Def_{v,Z}}  (\prod_{v\in S\backslash \Sigma} s\cD_v \times \prod_{v\in \Sigma} s\cD_v^{\square}).
	\end{displaymath} 
	Indeed, here it is not even necessary to modify the center.  But in order to have $\pi_0 s\cD_{\cS}^{\Sigma}(A) \cong \cD_{\cS}^{\Sigma}(A)$ for constant ring $A$, we need to suppose $H^0(\Gamma_v, \g_k) = \z_k$ for $v\in \Sigma$ (this is true for minimal conditions).  Then the functor $s\cD_{\cS}^{\Sigma}$ is pro-representable (say by $\cR_{\cS}^{\Sigma}$) and the natural transformation (up to weak equivalence) $s\cD_{\cS}^{\Sigma} \to \prod_{v\in \Sigma} s\cD_v^{\square}$ induces $A^{\Sigma} \to \cR_{\cS}^{\Sigma}$ (up to weak equivalence) where $A^{\Sigma} = \widehat{\otimes}_{v\in \Sigma} R_v$ is regarded as a pro-Artinian ring.  
	
	Under the assumption that every $R_v$ ($v\in S$) is formally smooth, it's not difficult to prove that the relative tangent complex $\t(\cR_{\cS}^{\Sigma}, A^{\Sigma})$ (see \cite[Definition 4.1]{GV18}) satisfies 
	\begin{displaymath}
	\t^i (\cR_{\cS}^{\Sigma}, A^{\Sigma}) \cong H_{\cS, \Sigma}^{i+1} (\Gamma_S, \g_k)
	\end{displaymath}
	for $i\geq 0$ (see \cite[(6.2.22)]{ACC+18} for $H_{\cS, \Sigma}^{\ast} (\Gamma_S, \g_k)$).  
\end{rem}

\subsubsection{Relative derived deformations}

Let $B \in \Art_{\cO}$ and let $\rho_B \colon \Gamma_S \to G(B)$ be a fixed lifting of type $\cS$.  \cite{TU21} considered the derived deformation functor of type $\cS$ over $\rho_B$ (denoted by $s\cD_{\cS, B}$).  Essentially it's the functor $s\cD_{\cS}$ restricted to $\sSet_{/ \BG(B)}$.  Note $\rho_B$ induces a map $\cR_{\cS} \to \pi_0 \cR_{\cS} \to B$, and with this specified map, $\cR_{\cS}$, as a pro-object in $\sArtB$, represents $s\cD_{\cS, B}$.

We calculate $\pi_i s\cD_{\cS, B} (B \oplus M[n])$ instead of the tangent complex, where $M$ is a finite module over $B$ and $M[n]$ means the Dolk-Kan of the chain complex $M$ concentrated on degree $n$.  In fact the procedures of proving \lemref{tangent complex of universal derived deformation functor} and \lemref{tangent complex of cR_S}  can be genralized directly, and one finds 
$$\pi_i s\cD_{\cS, B}(B\oplus M[n]) \cong H^{n-i+1}_{\cS} (\Gamma_S, \g_B \otimes_B M) \quad (\forall i,n \geq 0),$$
where $\g_B=\Lie(G/\cO) \otimes_{\cO} B$.  Moreover, by the discussions in Section \ref{tangent complexes of simplicial commutative rings}, the complex $C^{\ast+1}_{\cS} (\Gamma_S, \g_B \otimes_B M)$ is quasi-isomorphic to $[L_{\cR_{\cS}/\cO}  \otimes_{\cR_{\cS}} B, M]$ (here $\cR_{\cS}$ is regarded as a pro-object to take into account the continuity).

\subsection{Taylor-Wiles descent}

Now we are able to generalize \cite[Theorem 14.1]{GV18}.  We follow the approach of \cite{GV18}, but make minor modifications to fit our more general situation.  

We keep the settings in Section \ref{Calegari-Geraghty}.  Recall that $\zeta_p \notin F$ and $\bar\rho$ is supposed to have an enormous image.  Write $\cQ=(Q_m)_{m \geq 1}$ for a system of disjoint allowable Taylor-Wiles data (see \defref{allowable Taylor-Wiles}) such that each $Q_m$ is of level $m$ and cardinal $r\geq \dim_k H^1_{\cS}(\Gamma_S, \g_k^{\ast})$, and write $\Gamma_m = \Gamma_{S\cup Q_m}$, $\cD_m = \cD_{\cS_{Q_m}}$ and $R_m = R_{\cS_{Q_m}}$.  Let 
$$s\cD_m = s\Def_{S\cup Q_m,Z} \times^h_{\prod_{v\in S} s\Def_{v,Z}}  \prod_{v\in S} s\cD_v.$$

Note we don't put the derived unconditional deformation condition for $v\in Q_m$ for it's not formally smooth, but as \lemref{pi0 sD_S}, it's easy to see that $\pi_0 s\cD_m(A) \cong  \cD_m(A)$ for $A \in \Art_{\cO}$.  Moreover, $\t^{-1} s\cD_m$ is obviously trivial so $s\cD_m$ is pro-representable, say by $\cR_m$.

Let's fix $m\geq 1$.  By the definition of allowable Taylor-Wiles data, we have $H^2_{\cS_{Q_m}}(\Gamma_m, \g_k) = 0$.  Hence we have the exact sequence (see \remref{relative dimension remark}, note $L_v = H^1(\Gamma_v, \g_k)$ for $v\in Q_m$) 
\begin{align}\label{eqn 2}
0 \to &H^1_{\cS_{Q_m}}(\Gamma_m, \g_k) \to H^1(\Gamma_m, \g_k) \xrightarrow{A_m} \bigoplus_{v \in S} H^1(\Gamma_v, \g_k)/L_v \notag \\
\to &0 \to H^2(\Gamma_m, \g_k) \xrightarrow{B_m} \bigoplus_{v \in S \cup Q_m} H^2(\Gamma_v, \g_k) \to 0.
\end{align}
In particular, $B_m$ is an isomorphism.

We use $s\Def_v^{\ur}$ to denote the derived local deformation functor for the unramified condition.  For a Taylor-Wiles prime $v$, recall that $\bar\rho \vert_{\Gamma_v} \colon \Gamma_v \to G(k)$ is conjugated to some $\bar\rho_v^T \colon \Gamma_v \to T(k)$.  We write $s\Def_v^T$ (resp. $s\Def_v^{T,\ur}$) for the derived universal deformation functor for $\bar\rho_v^T \colon \Gamma_v \to T(k)$ (resp. $\bar\rho_v^T \vert_{\Gamma_v/ I_v}\colon \Gamma_v/ I_v \to T(k)$).

\begin{lem}
	Let $v$ be a Taylor-Wiles prime.  In the natural commutative diagram
	\begin{displaymath}
	\xymatrix{
		s\Def_S \ar[r]\ar[d] & s\Def_v^{\ur} \ar[d] & s\Def_v^{T,\ur} \ar_{\sim}[l] \ar[d] \\
		s\Def_{S\cup \{v\}} \ar[r] & s\Def_v  & s\Def_v^{T}, \ar_{\sim}[l] 
	}
	\end{displaymath}
	the first square is a homotopy pullback square, and the arrows with $\sim$ are objectwise weak equivalences.
\end{lem}

\begin{proof}
	See \cite[Section 8.2]{GV18} for the first statement, and \cite[Section 8.3]{GV18} for the second.
\end{proof}

We thus obtain a homotopy pullback square up to weak equivalences
\begin{displaymath}
\xymatrix{
	s\Def_S \ar[r]\ar[d] & s\Def_v^{T,\ur}  \ar[d] \\
	s\Def_{S\cup \{v\}} \ar[r]  & s\Def_v^{T}.
}
\end{displaymath}
In order that the functors involved are pro-representable, we need to modify their centers as in Section \ref{modifying the center}.  We use $s\Def_{v,T}^T$ (resp. $s\Def_{v,T}^{T,\ur}$) to denote the functor obtained from $s\Def_v^{T}$ (resp. $s\Def_v^{T,\ur}$) by modifying the center (the cumbersome notations just say that the center of $T$ is $T$ itself).  By \remref{modifying the center functorial} we have the commutative diagram 
\begin{displaymath}
\xymatrix{
	s\Def_{S,Z} \ar[r]\ar[d] & s\Def_{v,T}^{T,\ur}  \ar[d] \\
	s\Def_{S\cup \{v\}, Z} \ar[r]  & s\Def_{v,T}^{T}.
}
\end{displaymath}

\begin{lem}
	The above diagram is a homotopy pullback square.
\end{lem}

\begin{proof}
	The diagram is a homotopy pullback square if and only if the sequence 
	\begin{align*}
	0 \to &\t^0 s\Def_{S,Z}  \to \t^0 s\Def_{v,T}^{T,\ur} \oplus \t^0 s\Def_{S\cup \{v\}, Z} \to \t^0 s\Def_{v,T}^{T} \\
	\to &\t^1 s\Def_{S,Z}  \to \t^1 s\Def_{v,T}^{T,\ur} \oplus \t^1 s\Def_{S\cup \{v\}, Z} \to \t^1 s\Def_{v,T}^{T} \\
	\to &\dots
	\end{align*}
	is exact.  This follows from the homotopy pullback square before modifying the center and the fact that modifying the center doesn't change $\t^i$ for $i \geq 0$.
\end{proof}

By repeating the procedure of adding Taylor-Wiles primes, we can replace $v$ by a Taylor-Wiles datum $Q_m$.  Moreover, by applying 
$$- \times^h_{\prod_{v\in S} s\Def_{v,Z}}  \prod_{v\in S} s\cD_v $$ 
to the first vertical arrow, we can replace $s\Def_{S,Z} \to s\Def_{S\cup Q_m,Z}$ by $s\cD_{\cS} \to s\cD_m$.  The following corollary is clear:

\begin{cor}\label{homotopy pullback for T-W datum}
	Let $Q_m$ be a Taylor-Wiles datum.  Then we have the homotopy pull back square
	\begin{displaymath}
	\xymatrix{
		s\cD_{\cS} \ar[r]\ar[d] & \prod_{v \in Q_m} s\Def_{v,T}^{T,\ur}  \ar[d] \\
		s\cD_m \ar[r]  & \prod_{v \in Q_m} s\Def_{v,T}^{T},
	}
	\end{displaymath}
	and consequently we have an objectwise weak equivalence 
	$$s\cD_{\cS} \iso s\cD_m \times^h_{\prod_{v \in Q_m} s\Def_{v,T}^{T}} \prod_{v \in Q_m} s\Def_{v,T}^{T,\ur}.$$ 
\end{cor}

Now we pass to the level of rings.  In Section \ref{simplicial commutative ring} we defined the "derived" tensor product $\uotimes$ for simplicial commutative rings; this can be extended for pro-objects in $\sArt$ indexed by natural numbers (we have to take the Postnikov truncations for $\sArt$ is not closed under the tensor product, and then we can suppose the resulting pro-ring is nice for convenience, see discussions around \cite[Definition 3.3]{GV18}), with the property that $\cR_1 \uotimes_{\cR_3} \cR_2$ is a pro-objects of $\sArt$ representing the homotopy pullback of 
$$\sHom_{\sCROk}(\cR_1, -) \to \sHom_{\sCROk}(\cR_3, -) \leftarrow \sHom_{\sCROk}(\cR_2, -).$$ 

We say a map $\cR \to \cS$ between pro-$\sArt$ objects is a weak equivalence if it induces a weak equivalence on represented functors after applying level-wise cofibrant replacements (see \cite[Definition 7.4]{GV18}), and we say a pro-object $\cR$ of $\sArt$ is homotopy discrete if the map $\cR \to \pi_0 \cR$ is a weak equivalence.  

Let $\cS_m$ (resp. $\cS_m^{\ur}$) be a pro-object of $\sArt$ which represents $\prod_{v \in Q_m} s\Def_{v,T}^T$ (resp. $\prod_{v \in Q_m} s\Def_{v,T}^{T,\ur}$).  By \lemref{functor  map to ring map} applying to the weak equivalence 
$$s\cD_{\cS} \iso s\cD_m \times^h_{\prod_{v \in Q_m} s\Def_{v,T}^{T}} \prod_{v \in Q_m} s\Def_{v,T}^{T,\ur},$$ 
there is a weak equivalence of representing rings $\cR_{\cS} \to \cR_m \uotimes_{\cS_m} \cS_m^{\ur}$ (note \lemref{functor  map to ring map} allows us to reverse the arrrow).  This map between pro-$\sArt$ objects is an isomorphism in the pro-homotopy category by \cite[Lemma 3.14]{GV18}, so the isomorphism $\pi_{\ast} \cR_{\cS} \to \pi_{\ast}(\cR_m \uotimes_{\cS_m} \cS_m^{\ur})$ of pro-graded $\cO$-algebras is well-defined.

\begin{lem}
	The pro-objects $\cS_m^{\ur}$ and $\cS_m$ are homotopy discrete.
\end{lem}

\begin{proof}
	Note \cite[Lemma 7.5]{GV18} asserts that a pro-object $\cR$ of $\sArt$ such that $b_i = \dim_k \t^i \cR$ is zero except for $i=0,1$ is homotopy discrete if and only if the complete local ring associated to $\pi_0 \cR$ is isomorphic to a quotient of $\cO[[X_1, \dots, X_{b_0}]]$ by a regular sequence of length $b_1$.  
	
	By \lemref{center fiber sequence}, $\cS_m$ and $\cS_m^{\ur}$ represent the derived framed deformation functors $\prod_{v \in Q_m} s\Def_v^{T,\square}$ and $\prod_{v \in Q_m} s\Def_v^{T,\ur,\square}$.  Hence it suffices to show the classical (framed) universal deformation ring $\Sigma_v$ (resp. $\Sigma_v^{\ur}$) for $\bar\rho_v^T \colon \Gamma_v \to T(k)$ (resp. $\bar\rho_v^T \vert_{\Gamma_v/ I_v}\colon \Gamma_v/ I_v \to T(k)$) where $v$ is a Taylor-Wiles prime is a complete intersection ring of expected dimension.  
	
	\begin{enumerate}
		\item For $\Sigma_v^{\ur}$, it's easy to see that $b_i = \dim_k \t^i s\Def_{v,T}^{T,\ur}$ vanishes for $i\neq 0$, and $b_0=  \dim_k H^1_{\ur}(\Gamma_v, \t_k) = n$.  So it suffices to show $\Sigma_v^{\ur} \cong \cO[[X_1,\dots, X_n]]$.  But $\Sigma_v^{\ur}$ is the classical universal deformation ring for $\Def_v^{T,\ur}$, which is represented by $\cO[[X^{\ast}(T) \otimes \widehat{\Z} ]] \cong \cO[[X_1,\dots, X_n]]$ (see \cite[Proposition 4.2]{Til96}).
		\item For $\Sigma_v$, we have 
		$$b_i = \dim_k \t^i s\Def_{v,T}^T =  \left\{\begin{array}{l}{  \dim_k H^{i+1}(\Gamma_v, \t_k), \text{ if } i\geq 0; } \\ { 0, \text{ if } i<0.}\end{array}\right. $$  
		So $b_0 = 2n$, $b_1=n$ and $b_i=0$ for $i\neq 0,1$.  It suffices to check that $\Sigma_v$ is isomorphic to $\cO[[X_1, \dots, X_{2n}]]/(Y_1, \dots, Y_{n})$ for a regular sequence $(Y_i)$.  By \cite[Proposition 4.2]{Til96}, the classical representing ring for $\Def_v^T$ is isomorphic to $\cO[[X^{\ast}(T) \otimes F_v^{\ast,(p)}]]$ (here $(p)$ means the pro-$p$ completion).  Recall $\Delta_v$ is the Sylow $p$-subgroup of $(k_v^{\ast})^n$.  We have $X^{\ast}(T) \otimes F_v^{\ast,(p)} \cong \Delta_v\times \widehat{\Z}^n$ and hence $\Sigma_v \cong \cO[[X^{\ast}(T) \otimes F_v^{\ast,(p)}]] \cong \cO[\Delta_v][[X_1, \dots, X_n]]$ as expected.
	\end{enumerate}
\end{proof}

Let $\Sigma_m^{\ur} = \cO[[X_1, \dots, X_{nr}]]$ and $\Sigma_m = \cO[\Delta_{Q_m}][[X_1, \dots, X_{nr}]]$ (here $\Delta_{Q_m} = \prod_{v \in Q_m} \Delta_v$).  For convenience we also use $\Sigma_m^{\ur}$ and $\Sigma_m$ to denote the associated pro-Artinian rings.  Then the above lemma just says that $\cS_m$ is weakly equivalent to $\Sigma_m$ and $\cS_m^{\ur}$ is weakly equivalent to $\Sigma_m^{\ur}$.  

Note that $I_v \to \Gamma_v \to \Gal(\bar{k}_v/k_v)$ for $v\in Q_m$ induces $\cO[ \Delta_{Q_m} ] \to \Sigma_m \to \Sigma_m^{\ur}$.

\begin{lem}\label{a homotopy pullback square}
	The commutative diagram 
	\begin{displaymath}
	\xymatrix{
		\cO[ \Delta_{Q_m} ] \ar[r]\ar[d] & \Sigma_m  \ar[d] \\
		\cO \ar[r]  & \Sigma_m^{\ur}
	}
	\end{displaymath}
	induces a homotopy pullback square of represented functors after cofibrant replacements.  
\end{lem}

\begin{proof}
	It suffices to note that $\Sigma_m$ is obtained from $\cO[ \Delta_{Q_m} ]$ by adding $nr$ free variables, and $\Sigma_m^{\ur}$ is obtained from $\cO$ by adding $nr$ free variables.
\end{proof}

Recall in Section \ref{Calegari-Geraghty} we've defined $S_m = S_{\infty}/J_m$ which is a quotient of $\cO[\Delta_{Q_m}]$.  Also we've introduced $\bar{S}_m = S_m /p^m$, $\bar{R}_m = R_m \otimes_{\cO[[\Delta_{Q_m}]]} \bar{S}_m$ and a constant $c(m)$ such that $\bar{R}_m \to \End_{\cO}(H^{\ast}(C_m^{\ast}))$ factors through $\bar{R}_m / \m_{\bar{R}_m}^{c(m)}$.  Without loss of generality, we may suppose $\bar{R}_m / \m_{\bar{R}_m}^{c(m)} \leftarrow \bar{S}_m \rightarrow \cO/p^m$ forms a compatible projective system for $m\in \N^{\ast}$.  We remark that the cohomology of locally symmetric space is not involved explicitly here, but finally we will need $R_{\infty} \cong \varprojlim_m \bar{R}_m / \m_{\bar{R}_m}^{c(m)}$, which is true only if the numerical coincidence holds (see the proof of \corref{Calegari-Geraghty corollary} (3)).  

For each $m\geq 1$ we have (still apply \lemref{functor  map to ring map} to reverse the weak equivalences) 
\begin{displaymath}
f_m \colon \cR_{\cS} \iso \cR_m \uotimes_{\Sigma_m} \Sigma_m^{\ur} \iso \cR_m \uotimes_{\cO[ \Delta_{Q_m} ]} \cO \to \bar{R}_m / \m_{\bar{R}_m}^{c(m)}  \uotimes_{\bar{S}_m} \cO/p^m =\colon	 \cC_m.
\end{displaymath}

We have $\Tor_{\ast}^{S_\infty}(R_\infty, \cO) \cong \pi_{\ast} (R_{\infty} \uotimes_{S_{\infty}} \cO) \cong \varprojlim_m \pi_{\ast} (\bar{R}_m / \m_{\bar{R}_m}^{c(m)} \uotimes_{\bar{S}_m} \cO/p^m)$ as graded commutative $\cO$-algebras.  Here the first isomorphism follows from Section \ref{graded} and the second isomorphism follows from \cite[Lemma 7.6]{GV18}.

For each $n>m$, there is a natural map 
$$e_{n,m} \colon \cC_n =\bar{R}_n / \m_{\bar{R}_n}^{c(n)}  \uotimes_{\bar{S}_n} \cO/p^n \to \cC_m= \bar{R}_m / \m_{\bar{R}_m}^{c(m)}  \uotimes_{\bar{S}_m} \cO/p^m,$$
but a prior, the maps $f_m \colon \cR_{\cS} \to \cC_m$ $(m\geq 1)$ don't form a compatible system under $e_{n,m}$, and we have to do another patching so that $e_{n,m} \circ f_n$ ($n>m$) are compatible modulo homotopy.  The key observation is that $\t \cR_{\cS}$ is finite dimensional, so each homotopy class of maps $\cR_{\cS} \to \cC_m$ as pro-$\sArt$ objects is indeed finite (see \cite[Page 100]{GV18}).  

Consider the projective system of homotopy classes of maps $\cR_{\cS} \to \cC_m$ ($m\geq 1$) induced from $e_{n,m}$, then we can choose a subsequence of $(f_m)$ such that $e_{n,m} \circ f_n$ is homotopic to $f_m$ for every $f_n,f_m$ ($n>m$) in that subsequence.  Without loss of generality, we may simply suppose $(f_m)_{m\geq 1}$ is such a sequence, and then $(f_m)_{m\geq 1}$ induces $\hocolim_m \, \sHom_{\sCROk}(\cC_m, - ) \to \sHom_{\sCROk}(\cR_{\cS}, - )$.

Now we prove $\pi_{\ast}\cR_{\cS} \cong \Tor_{\ast}^{S_{\infty}}(R_{\infty},\cO)$.  Let's recall the setting: 
\begin{enumerate}
	\item $\bG$ is a connected reductive algebraic group defined over a number field $F$ and $G = {}^L \bG$;
	\item $p$ is an odd prime number which is very good for $G$ and satisfies $\zeta_p \notin F$;
	\item $\bar\rho \colon \Gamma_S \to G(k)$ is an absolutely irreducible Galois representation associated to some cuspidal automorphic representation occuring in $H^{\ast}(X_{\bG}^U, \widetilde{V}_{\lambda}(\cO))_{\m}$ which fits our assumption $(\Resm)$;
	\item we assume the conjectures $(\Galm)$ and $(\Vanm)$. 
\end{enumerate}

\begin{thm}\label{main theorem}
	With the above notations, there is an isomorphism of graded commutative $\cO$-algebras $\pi_{\ast}\cR_{\cS} \cong \Tor_{\ast}^{S_{\infty}}(R_{\infty},\cO)$ (where $\pi_{\ast} \cR_{\cS}$ is defined as the projective limit).  Moreover, $H^{\ast}(X_\bG^U, \widetilde{V}_{\lambda}(\cO))_{\m}$ is a graded $\pi_{\ast}\cR_{\cS}$-module freely generated by $H^{q_0+\ell_0}(X_\bG^U, \widetilde{V}_{\lambda}(\cO))_{\m}$.  
\end{thm}

\begin{rem}
	We have supposed special types of local deformation problems in $(\Vanm)$, but essentially what we require are:
	\begin{enumerate}
		\item the numerical coincidence $\dim_k H^1_{\cS}(\Gamma_S, \g_k) - \dim_k H^1_{\cS^{\perp}}(\Gamma_S, \g_k^{\ast}) = -\ell_0$ holds;
		\item the local deformation problems have formally smooth framed representing rings.
	\end{enumerate}
\end{rem}

\begin{proof}
	We will prove the first assertion and the second is an immediate consequence.  
	
	By above discussions, it suffices to prove 
	$$\hocolim_m \, \sHom_{\sCROk}(\cC_m, - ) \to \sHom_{\sCROk}(\cR_{\cS}, - )$$ 
	is a weak equivalence of natural transformations, and by \lemref{conservativity of tangent complex} it suffices to show 
	\begin{displaymath}
	\t^i (\hocolim_m \, \sHom_{\sCROk}(\cC_m, - )) \to \t^i \sHom_{\sCROk}(\cR_{\cS}, - )
	\end{displaymath}
	is an isomorphism for all $i\geq 0$.
	
	For $m\geq 1$, $\t \cC_m$ fits into the exact triangle $\t \cC_m \to \t (\bar{R}_m / \m_{\bar{R}_m}^{c(m)}) \oplus \t (\cO/p^m) \to \t \bar{S}_m$, and by taking colimits over $m$ we get the following exact sequence: 
	\begin{displaymath}
	\t^i (\hocolim_m  \, \sHom_{\sCROk}(\cC_m, - )) \to \t^i R_{\infty} \to \t^i S_{\infty} \stackrel{[1]}{\to} \dots,
	\end{displaymath}
	so the Euler characteristic for $\t (\hocolim_m \, \sHom_{\sCROk}(\cC_m, - ))$ is $\dim R_{\infty} - \dim S_{\infty}$.  On the other hand, by \lemref{tangent complex of cR_S}, the Euler characteristic for $\t (\sHom_{\sCROk}(\cR_{\cS}, - ))$ is $\dim_k H^1_{\cS}(\Gamma_S, \g_k) - \dim_k H^1_{\cS^{\perp}}(\Gamma_S, \g_k^{\ast})$, which is equal to $\dim R_{\infty} - \dim S_{\infty}$ by \lemref{relative dimension 1}.  We also find that both tangent complexes are concentrated on degrees $0$ and $1$.  Thus it suffices to show $\t^i \cC_m \to \t^i \cR_{\cS}$ is an isomorphism for $i=0$ and a surjection for $i=1$, or equivalently by \lemref{a homotopy pullback square}, $\t^i (\bar{R}_m / \m_{\bar{R}_m}^{c(m)} \uotimes_{\bar{\Sigma}_m} \bar{\Sigma}_m^{\ur}) \to \t^i (\cR_m \uotimes_{\Sigma_m} \Sigma_m^{\ur})$ is an isomorphism for $i=0$ and a surjection for $i=1$, where $\bar{\Sigma}_m= \Sigma_m \otimes_{\cO[\Delta_{Q_m}]} \bar{S}_m$ and $\bar{\Sigma}_m^{\ur} = \Sigma_m^{\ur} / p^m$.  
	
	Consider the following commutative diagram with exact rows: 
	\begin{displaymath}
	\begin{array}{l}
	{
		\xymatrix{
			0 \ar[r] & \t^0 (\bar{R}_m / \m_{\bar{R}_m}^{c(m)} \uotimes_{\bar{\Sigma}_m} \bar{\Sigma}_m^{\ur}) \ar[r]\ar[d]^{j_1} &\t^0 (\bar{R}_m / \m_{\bar{R}_m}^{c(m)}) \oplus \t^0 \bar{\Sigma}_m^{\ur} \ar[r]\ar[d]^{f} & \t^0 \bar{\Sigma}_m \ar[d]^{g} \\
			0 \ar[r] & \t^0 (\cR_m \uotimes_{\Sigma_m} \Sigma_m^{\ur}) \ar[r] &\t^0 \cR_m \oplus \t^0 \Sigma_m^{\ur} \ar[r] & \t^0 \Sigma_m 
		}
	}\\
	{
		\xymatrix{
			\ar[r] & \t^1 (\bar{R}_m / \m_{\bar{R}_m}^{c(m)} \uotimes_{\bar{\Sigma}_m} \bar{\Sigma}_m^{\ur}) \ar[r]\ar[d]^{j_2} &\t^1 (\bar{R}_m / \m_{\bar{R}_m}^{c(m)}) \oplus \t^1 \bar{\Sigma}_m^{\ur} \ar[d]
			\\
			\ar[r] & \t^1(\cR_m \uotimes_{\Sigma_m} \Sigma_m^{\ur}) \ar[r]^{\gamma} & \t^1 \cR_m \oplus \t^1 \Sigma_m^{\ur} \ar[r]^{\qquad h} & \t^1 \Sigma_m,
		}
	}\\
	\end{array}
	\end{displaymath}
	The maps $f$, $g$ are clearly isomorphisms.  By a diagram chasing, it suffices to show $h$ is an isomorphism.  Note $\t^1 \Sigma_m^{\ur} = 0$ and $\t^1 \Sigma_m \cong \prod_{v \in Q_m} H^2(\Gamma_v, \t_k) \cong \prod_{v \in Q_m}H^2(\Gamma_v, \g_k)$, so it remains to prove $\t^1 \cR_m \cong \prod_{v \in Q_m}H^2(\Gamma_v, \g_k)$.  As \lemref{tangent complex of cR_S}, we have the exact sequence 
	\begin{align*}
	0\to &\t^0 \cR_m \to H^1(\Gamma_m, \g_k) \to \bigoplus_{v \in S} H^1(\Gamma_v,\g_k) / L_v \\
	\to &\t^1 \cR_m \to H^2(\Gamma_m, \g_k) \to \bigoplus_{v\in S} H^2(\Gamma_v,\g_k) \to 0.
	\end{align*}
	By comparing it with the exact sequence \eqref{eqn 2}, we conclude $\t^1 \cR_m \cong \prod_{v \in Q_m}H^2(\Gamma_v, \g_k)$.	
\end{proof}                                                                   

\begin{rem}
	The formal smoothnesses for local deformation rings play an essential role (especially in \lemref{tangent complex of cR_S}) in the above calculations.  A natural question is to genralize the result without the formally smooth assumptions (for example firstly for local complete intersection rings).  However, we do not yet have a clear answer to this question.
\end{rem}

\section{Examples}\label{examples}

\subsection{General linear groups}

We keep the notations of the previous section.  Suppose $F$ is a number field with $r_1$ real places and $r_2$ complex places, and consider the locally symmetric spaces associated to $\Res_{\Q}^F \GL_N$.  The maximal compact subgroup of $\GL_N(\R)$ is $\OO(N)$ and the maximal compact subgroup of $\GL_N(\C)$ is $\U(N)$, so we have 
\begin{displaymath}
\left\{\begin{array}{l}{  2q_0+\ell_0 =  (N^2-\frac{N(N-1)}{2}) r_1 +(2N^2 - N^2 ) r_2 -1 = \frac{N^2+N}{2} r_1+ N^2 r_2 -1;  } \\ { \ell_0 =  (N- [\frac{N}{2}])r_1 +(2N-N) r_2 -1 = (N-[\frac{N}{2}]) r_1 + Nr_2-1 ,}\end{array}\right.
\end{displaymath}
and consequently $q_0 =[ \frac{N^2}{4}] r_1 + \frac{N^2-N}{2} r_2$.

We suppose 
\begin{enumerate}
	\item $\pi_v$ is minimal for $v\in S\backslash S_p$;
	\item either $\pi_v$ is regular ordinary for every $v\in S_p$, or $p$ is unramified in $F$ and $\lambda_{\tau,1} - \lambda_{\tau,n} <p-n$ for all $\tau$.
\end{enumerate}

In \cite{HLTT16} the authors proved that there exists a Galois representation $\rho_{\pi}\colon  \Gamma_S \to \GL_N(\cO)$ associated to $\pi$ such that $\bar\rho = \rho_{\pi} \pmod \varpi$ satisfies $(\Resm)$.  In the ordinary case, we suppose $\bar\rho \vert_{\Gamma_v}$ is regular and dual regular (see Section \ref{local deformation problems}, and these are called distinguishability and strong distinguishability assumptions in \cite[Page 3-4]{TU21}).  Let $\cD_v^{\min}$, $\cD_v^{\ord}$ and $\cD_v^{\FL}$ be the minimal, ordinary and Fontaine-Laffaille local deformation functors respectively (see Section \ref{local deformation problems}), then we have  

\begin{pro}
	The functors $\cD_v^{\min}$, $\cD_v^{\ord}$ and $\cD_v^{\FL}$ are liftable, and the framed representing rings for these functors are formally smooth.  Moreover, for $\cD_v^{\ord}$ and $\cD_v^{\FL}$, we have $\dim_k L_v - \dim_kH^0(\Gamma_v,\g_k) = [F_v:\Q_p]\frac{N(N-1)}{2}$; for $\cD_v^{\min}$, we have $\dim_k L_v - \dim_k H^0(\Gamma_v,\g_k) =0$.  
\end{pro}

Let $\cS$ to be the global deformation problem for $\bar\rho \colon \Gamma_S \to \LG(k)$ which is simultaneously either ordinary or Fontaine-Laffaille for $v\in S_p$, and minimal for $v \in S\backslash S_p$.  Let's check the condition $\dim_k H^1_{\cS}(\Gamma_S, \g_k) - \dim_k H^1_{\cS^{\perp}}(\Gamma_S, \g_k^{\ast}) = -\ell_0$ (see \lemref{relative dimension}):

\begin{lem}
	Let the notations be as above.  Then 
	$$ - 1 + \sum\limits_{v \divides \infty} \dim_k H^0(\Gamma_v,\g_k) - \sum\limits_{v\in S} (\dim_k L_v - \dim_kH^0(\Gamma_v,\g_k) ) = \ell_0$$
	holds if and only if the action of complex conjugation on $\g_k$ is odd for every real place of $F$.
\end{lem}

\begin{proof}
	By the above proposition we have $\sum\limits_{v\in S} (\dim_k L_v - \dim_kH^0(\Gamma_v,\g_k) ) = \frac{n^2-n}{2}r_1 + (n^2-n) r_2 $.  So the condition  is equivalent to $\sum\limits_{v \divides \infty} \dim_k H^0(\Gamma_v,\g_k) = [\frac{n^2+1}{2}] r_1 + n^2 r_2$.  But for each $v$ real, $H^0(\Gamma_v,\g_k)$ is at least $[\frac{n^2+1}{2}]$, so we must have the equality, which is exactly the oddness condition.
\end{proof}

Now it remains to check $(\Galm)$ and $(\Vanm)$ for \thmref{main theorem}.  To the author's knowledge. the hypothesis $(\Vanm)$ is still far from reach except for $\GL_3$ over $\Q$ and $\GL_2$ over $F$ where $F$ satisfies $r_1+r_2\leq 2$ (or $q_0 \leq 2$), where it can be solved using the congruence subgroup problem (see \cite{PR10}).  For the hypothesis $(\Galm)$, in \cite[Theorem 2.3.7]{ACC+18}, the authors construct a map $\Gamma_S \to \GL_N(\T/I)$, where $I$ is a nilpotent ideal, with desired characteristic polynomials for $v\notin S$.  In subsequent sections 3,4,5 of \cite{ACC+18}, the local-global compatibilities are established for minimal, Fontaine-Laffaille and ordinary places, given some additional restrictions listed there.  The nilpotent ideal $I$ is eliminated in \cite[Theorem 6.1.4]{CGH+20} under the assumption that $p$ splits completely in $F$, however, the local-global compatibility hasn't been established yet.

\subsection{Orthogonal similitude groups}

Consider the locally symmetric spaces associated to the orthogonal similitude groups $\GSO_{a,b}$ over $\Q$.  Recall that 
$$\GO_{a,b}(R) = \{ g\in \GL_{a+b}(R) \mid  g^t \begin{pmatrix}
\id_a & 0  \\
0 & -\id_b
\end{pmatrix} g = \lambda \begin{pmatrix}
\id_a & 0 \\
0 & -\id_b
\end{pmatrix} \text{ for some } \lambda\in R^{\ast}  \}, $$
and $\GSO_{a,b}$ is the connected component of the identity in $\GO_{a,b}$ (so $\GSO_{a,b} = \GO_{a,b}$ if $a+b$ is odd).    

We still have \thmref{main theorem} once all necessary hypotheses are verified.  But when $a+b$ is small, it seems more convenient to approach \thmref{main theorem} via the special (local) isomorphisms listed in \cite[\RNum{10}.6.4]{Hel01} and \cite{MY90}, for the auxiliary group under the isomorphism is often better understood.  

It's easy to see that $\GSO_{a,b}$ is abelian when $a+b\leq 2$, and the center $Z(\GSO_{a,b})$ consists of scalar matrices when $a+b>2$.  In the second case, the invariants $q_0$ and $\ell_0$ satisfy 
\begin{displaymath}
\left\{\begin{array}{l}{  q_0 =  [\frac{ab}{2}];  } \\ { \ell_0 = [\frac{a+b}{2}] - [\frac{a}{2}] - [\frac{b}{2}].}\end{array}\right.
\end{displaymath}

\subsubsection{Derived deformation rings under Langlands transfers}

Let's discuss how the derived deformation rings behave under Langlands transfers in general.  

Let $G$ and $H$ be a connected reductive linear algebraic group over $\Q$.  As in the introduction, we fix a finite set of finite places $S \supseteq S_p$ of $F$, an open compact group $U = U_S \times U^S = (\prod_{v\in S} U_v) \times (\prod_{v\notin S} U_v)$ with $U_v \subseteq \underline{H}(\cO_v)$ and each $U_v$ ($v\notin S \backslash S_p$) hyperspecial maximal.  Suppose $\pi_H$ is a cuspidal automorphic representation occuring in $H^{\ast}(X_H^U, \widetilde{V}_{\lambda}(\cO))_{\m}$ where $\m$ is a non-Eisenstein maximal ideal and we make the assumption $(\Resm)$ for the residual representation $\bar\rho_H \colon \Gamma_S \to \LH(k)$.  Suppose the Langlands transfer $r \colon \LH \to \LG$ is established, then there exists an automorphic $\pi_G$ in the global $L$-packet defined by $\pi_H$ and $r$, and the residual representations $\bar\rho_G \colon \Gamma_S \to \LG(k)$ satisfy $\bar\rho_G = r\circ \bar\rho_H$.  

Let $\cD_G$ (resp. $s\cD_G$) be the deformtion functor (derived deformtion functor) with suitable local conditions for $\bar\rho_G$, and we define similarly $\cD_H$ (resp. $s\cD_H$) with compatible local conditions.  Then there is a natural map $\cD_H \to \cD_G$ (resp. $s\cD_H \to s\cD_G$) induced by $r$, and hence a morphism $R_G \to R_H$ (resp. $\cR_G \to \cR_H$ up to weak equivalence) between the deformation rings (resp. derived deformation rings).  

In the following we will take $H = \GSO_{a,b}$ with $a+b=4$ or $6$.  Note then $\widehat{H} = \GSpin_{a+b}$.  Recall that $\GSpin_4$ can be identified with  
\begin{displaymath}
\{ (A,B)\in \GL_2 \times \GL_2  \mid \det(A) = \det(B) \},
\end{displaymath}
and $\GSpin_6$ can be identified with the subgroup of $\GL_1\times \GL_4$ defined by the exact sequence 
\begin{displaymath}
1 \to \GSpin_6 \to \GL_1\times \GL_4 \to \GL_1 \to 1
\end{displaymath}
with $\GL_1\times \GL_4 \to \GL_1$ is given by $(\lambda, g) \mapsto \lambda^{-2}\det(g)$.  

For $H = \GSO_{3,1}$ and $G = \Res_{\Q}^F \GL_2$ where $F$ is a quadratic imaginary field, the transfer is induced by the natural inclusion $\GSpin_4 \hookrightarrow \GL_2 \times \GL_2$.  Let $\Gamma_{F,S}$ be the Galois group of the maximal $S$-ramified extension of $F$ and let $\Gal(F/\Q) = \{1,c\}$.  Then $\LG = (\GL_2 \times \GL_2) \rtimes \{1,c\} $ and $\LH = \widehat{H} \rtimes \{1,c\} $, and the complex conjugation $c$ acts by exchanging the components in $\GL_2 \times \GL_2$.  Note $\widehat{H}$ can be identified with the subgroup
\begin{displaymath}
\{  
\begin{pmatrix}
a_1 & 0 & b_1 & 0 \\
0 & a_2 & 0 & b_2 \\
c_1 & 0 & d_1 & 0 \\
0 & c_2 & 0 & d_2
\end{pmatrix} \mid a_1d_1-b_1c_1 = a_2d_2-b_2c_2 \}
\end{displaymath}
of $\GSp_4$, and the action of $c$ is extended to the conjugation action by 
$\begin{pmatrix}
0 & 1 & 0 & 0 \\
1 & 0 & 0 & 0 \\
0 & 0 & 0 & 1 \\
0 & 0 & 1 & 0
\end{pmatrix} \in \GSp_4$.

\begin{lem}
	For $H = \GSO_{3,1}$ and $G = \Res_{\Q}^F \GL_2$, the map $\Def_H \to \Def_G$ between unconditional deformation functors is an isomorphism.
\end{lem}

\begin{proof}	
	We also use $c$ to denote the complex conjugation of $\Gamma_S$.  
	
	Let's first consider the functor $\Def_G$.  Let $A\in \Art_{\cO}$ and suppose $\rho_G \colon \Gamma_S \to \LG(A)$ is a lifting of $\bar\rho_G$.  For $\sigma \in \Gamma_{F,S}$, we write $\rho_G(\sigma) = ((M_{\sigma}, N_{\sigma}), 1)$.  Note $\rho_G(c) = ((X,X^{-1}),c)$ for some $X\in \GL_2(A)$, and without loss of generality up to conjugation we may suppose $X$ is the identity matrix.  Then it's easy to see $N_{\sigma} = M_{c\sigma c}$, so the deformation of $\rho_G$ is uniquely determined by the deformation for $\Gamma_{F,S} \to \GL_2$.  
	
	For $\Def_H$, we can only conjugate $\rho_G(c) = ((X,X^{-1}),c)$ to either $((\begin{pmatrix}
	1 & 0 \\
	0 & 1
	\end{pmatrix}, \begin{pmatrix}
	1 & 0 \\
	0 & 1
	\end{pmatrix}), c)$ or $((\begin{pmatrix}
	1 & 0 \\
	0 & -1
	\end{pmatrix}, \begin{pmatrix}
	1 & 0 \\
	0 & -1
	\end{pmatrix}), c)$.  But still the deformation of $\rho_G$ is uniquely determined by the deformation for $\Gamma_{F,S} \to \GL_2$. 	
\end{proof}

\begin{lem}
	For $H = \GSO_{3,1}$ and $G = \Res_{\Q}^F \GL_2$, the map $s\Def_H \to s\Def_G$ between unconditional derived deformation functors is a weak equivalence.
\end{lem}

\begin{proof}
	It suffices to check $s\Def_H \to s\Def_G$ induces a weak equivalence on tangent complexes.  Write $\h_k$ and $\g_k$ for the Lie algebras of $\LH$ and $\LG$ respectively (note $\h_k$ is a direct summand of $\g_k$), then it suffices to show $H^i(\Gamma_S, \h_k) \hookrightarrow H^i(\Gamma_S, \g_k)$ is an isomorphism for $i=1,2$.  
	
	For $i=1$ the isomorphism follows from the above lemma.  Note also the isomorphism for $i=0$ and the Euler characteristics $\chi(\Gamma_S, \h_k) = \chi(\Gamma_S, \g_k)$ (the subspace fixed by $c$ in $\g_k$ lies in $\h_k$), so $H^2(\Gamma_S, \h_k) \hookrightarrow H^2(\Gamma_S, \g_k)$ is also an isomorphism.  
\end{proof}

Similarly, the above lemmas hold for $H = \GSO_{2,2}$ and $G = \GL_2 \times \GL_2$ as well.

In the case $H = \GSO_{3,3}$ and $G = \GL_4$, these groups are split so we can identify $\LH$ with $\GSpin_6$ and identify $\LG$ with $\GL_4$.  The transfer $r \colon \GSpin_6 \to \GL_4$ is given by the second projection 
$$\GSpin_6 \subseteq \GL_1\times \GL_4 \to \GL_4.$$

\begin{lem}
	For $H = \GSO_{3,3}$ and $G = \GL_4$, the map $\Def_H \to \Def_G$ between unconditional deformation functors is an isomorphism.
\end{lem}

\begin{proof}
	Let $A \in \Art_{\cO}$ and let $\rho_H \colon \Gamma_S \to \GSpin_6(A)$ be a lifting of $\bar\rho_H \colon \Gamma_S \to \GSpin_6(k)$.  Suppose $\rho_H(\sigma) = (\lambda_{\sigma}, M_{\sigma})$.  If $M_{\sigma}$ is given, then there is a unique choice for such $\lambda_{\sigma}$ since $\lambda_{\sigma}^2$ and $\lambda_{\sigma} \pmod{\m_A}$ are determined.
\end{proof}

\begin{lem}
	For $H = \GSO_{3,3}$ and $G = \GL_4$, the map $s\Def_H \to s\Def_G$ between unconditional derived deformation functors is a weak equivalence.
\end{lem}

\begin{proof}
	It suffices to check $s\Def_H \to s\Def_G$ induces a weak equivalence on tangent complexes.  Write $\h_k$ and $\g_k$ for the Lie algebras of $\LH$ and $\LG$ respectively, then it's easy to see $\h_k \cong \g_k$, so $H^i(\Gamma_S, \h_k) \to  H^i(\Gamma_S, \g_k)$ is an isomorphism for $i=1,2$, and the conclusion follows. 
\end{proof}

The local conditions for $\rho_H \colon \Gamma_S \to \LH(A)$ should be essentially defined by the corresponding local conditions for $\rho_G \colon \Gamma_S \to \LG(A)$.  So in the cases
\begin{enumerate}
	\item $H = \GSO_{3,1}$ and $G = \Res_{\Q}^F \GL_2$, or
	\item $H = \GSO_{2,2}$ and $G = \GL_2 \times \GL_2$, or 
	\item $H = \GSO_{3,3}$ and $G = \GL_4$,
\end{enumerate}
we have the following: 

\begin{cor}
	The map $s\cD_H \to s\cD_G$ is a weak equivalence, and so is $\cR_G \to \cR_H$.  In particular, the map $\pi_{\ast} \cR_G \to \pi_{\ast} \cR_H$ is an isomorphism of graded commutative $\cO$-algebras. 
\end{cor}

If we could relate $H^{\ast}(X_H^U, \widetilde{V}_{\lambda}(\cO))_{\m}$ and $H^{\ast}(X_G^V, \widetilde{V}_{\lambda}(\cO))_{\m}$, then we are able to deduce \thmref{main theorem} for $H$ if it is known for $G$.  In the following we study the case $\GSO_{3,1}$.

\subsubsection{The case $\GSO_{3,1}$}

Write $H = \GSO_{3,1}$ and $G = \Res_{\Q}^F \GL_2$ where $F$ is an imaginary quadratic field.  By the preceding calculations, we know the $q_0$ and $\ell_0$ for both groups coincide.  We define $\phi \colon G \to H$ as follows:

Let $W = \{ x\in M_2(F) \mid x=x^{ct} \}$, then $\det \colon W \to \Q$ is a quadratic form of signature $(1,3)$, so $\GO_{3,1}$ can be identified with the group of orthogonal similitudes of $W$.  Let $A$ be the kernel of the norm map $N\colon \Res_{\Z}^{\cO_F} \G_m \to \G_m$.  Note that $W$ comes with a structure over $\cO_F$, we have the following commutative diagram of algebraic group schemes over $\Z$ with exact rows over algebraically closed fields:
\begin{displaymath}
\xymatrix{
	0 \ar[r] &A \ar[r]\ar@{=}[d] &\Res_{\Z}^{\cO_F} \G_m \ar[r]^{N}\ar[d] &\G_m \ar[r]\ar[d] &0 \\
	0 \ar[r] &A \ar[r] &\Res_{\Z}^{\cO_F} \GL_2 \ar[r]^{\phi} &\GSO_{3,1} \ar[r] &0. }
\end{displaymath}
Here $\phi$ is induced by associating $g\in \Res_{\Z}^{\cO_F} \GL_2$ to the endomorphism $x \mapsto gxg^{ct}$ on $W$, and the vertical maps are natural inclusions.  Note the similitude character of $\phi(g)$ is $\det(g)\det(g)^c$.

Let $U = U_S \times U^S = (\prod_{v\in S} U_l) \times (\prod_{v\notin S} U_l)$ be an open compact subgroup of $H_f$ such that each $U_v$ ($v\notin S$) is hyperspecial maximal.  We define $V_l$ to be the inverse image of $U_l$ under $\GL_2(\Z_l \otimes_{\Z} \cO_F) \to \GSO_{3,1}(\Z_l)$ and define $V = \prod_l V_l$.  Note the Langlands transfer $r\colon \LH \to \LG$ induces a map between the spherical Hecke algebras $\cH(G^S, V^S) \to \cH(H^S, U^S)$ via the Satake isomorphisms
\begin{displaymath}
C_c^{\infty}(H(\Q_l) \dslash U_l) \iso \C[\widehat{T}_H]^{W(\widehat{H}, \widehat{T}_H)(\C)}
\end{displaymath}  
and 
\begin{displaymath}
C_c^{\infty}(G(\Q_l) \dslash V_l) \iso \C[\widehat{T}_G]^{W(\widehat{G}, \widehat{T}_G)(\C)}.
\end{displaymath}  

Let $\lambda$ be a dominant weight for $\GSO_{3,1}$ and let $V_{\lambda}$ be the irreducible algebraic representation of $\GSO_{3,1}$ of highest weight $\lambda$.  By regarding $V_{\lambda}$ as an irreducible algebraic representation of $\Res_{\Z}^{\cO_F} \GL_2$ via $\phi\colon \Res_{\Z}^{\cO_F} \GL_2 \to \GSO_{3,1}$, we get a natural map $H^{\ast}(X_H^U, \widetilde{V}_{\lambda}(\cO)) \to H^{\ast}(X_G^V, \widetilde{V}_{\lambda}(\cO))$.  We make the assumption that we can choose $F$ such that $V \to U$ is surjective (note $V_l \to U_l$ is surjective for $l$ unramified in $F$).  The following proposition should be known by \cite{HST93} and \cite{Mok14}, nevertheless we will give a proof.

\begin{pro}\label{GSO31}
	The natural map $H^{\ast}(X_G^V, \widetilde{V}_{\lambda}(\cO)) \to H^{\ast}(X_H^U, \widetilde{V}_{\lambda}(\cO))$ is an isomorphism, and we have the commutative diagram of Hecke actions 
	\[
	\begin{tikzcd}
	\cH(H^S, U^S) & \cH(G^S, V^S) \arrow[l] \\
	H^{\ast}(X_H^U, \widetilde{V}_{\lambda}(\cO)) \arrow[loop above] \arrow[r, "\sim"]  & H^{\ast}(X_G^V, \widetilde{V}_{\lambda}(\cO)). \arrow[loop above]
	\end{tikzcd}
	\]
\end{pro}

\begin{cor}
	$H^{\ast}(X_H^U, \widetilde{V}_{\lambda}(\cO))_{\m}$ is a graded $\pi_{\ast}\cR_H$-module which is freely generated by $H^2(X_H^U, \widetilde{V}_{\lambda}(\cO))_{\m}$ (note $q_0+\ell_0 = 2$ here). 
\end{cor}

\begin{rem}
	Once the isomorphism between locally symmetric spaces and the compatibility with the Langlands transfer are established, it's easy to see that $(\Galm)$ and $(\Vanm)$ for $\pi_G$ implies those for $\pi_H$.  Together with the theory of Calegari-Geraghty we know (here $S_{\infty}^H$ and $R_{\infty}^H$ are limiting rings associated to $H$ constructed by the Taylor-Wiles method, and same for $S_{\infty}^G$ and $R_{\infty}^G$) 
	\begin{enumerate}
		\item $H^{\ast}(X_H^U, \widetilde{V}_{\lambda}(\cO))_{\m} \to H^{\ast}(X_G^V, \widetilde{V}_{\lambda}(\cO))_{\m}$ is an isomorphism;
		\item $H^{\ast}_{\m}(X_G^V, \widetilde{V}_{\lambda}(\cO))_{\m}$ is a graded module freely generated by $H^{q_0+\ell_0}(X_G^V, \widetilde{V}_{\lambda}(\cO))_{\m}$ over $\Tor_{\ast}^{S_{\infty}^G}(R_{\infty}^G,\cO)$;
		\item $H^{\ast}_{\m}(X_H^U, \widetilde{V}_{\lambda}(\cO))_{\m}$ is a graded module freely generated by $H^{q_0+\ell_0}(X_H^U, \widetilde{V}_{\lambda}(\cO))_{\m}$ over $\Tor_{\ast}^{S_{\infty}^H}(R_{\infty}^H,\cO)$.
	\end{enumerate}
	So we should have $\Tor_{\ast}^{S_{\infty}^H}(R_{\infty}^H,\cO) \cong \Tor_{\ast}^{S_{\infty}^G}(R_{\infty}^G,\cO)$.  In general it's seemingly more convinient to compare the derived deformation rings.
\end{rem}

We return to the commutative diagram of algebraic group schemes over $\Z$ with exact rows over algebraically closed fields:
\begin{displaymath}
\xymatrix{
	0 \ar[r] &A \ar[r]\ar@{=}[d] &\Res_{\Z}^{\cO_F} \G_m \ar[r]^{N}\ar[d] &\G_m \ar[r]\ar[d] &0 \\
	0 \ar[r] &A \ar[r] &\Res_{\Z}^{\cO_F} \GL_2 \ar[r]^{\phi} &\GSO_{3,1} \ar[r] &0. }
\end{displaymath}
For a field extension $E/\Q$, we have $H^1(\Gal(\bar{E}/E), (\bar{E}\otimes_{\Q} F)^{\ast})= H^1(\Gal(\bar{E}/E), \GL_2(\bar{E}\otimes_{\Q} F) ) = 0$ by Hilbert's Theorem 90, so we obtain the commutative diagram with exact rows 
\begin{displaymath}
\xymatrix{
	0 \ar[r] &A(E) \ar[r]\ar@{=}[d] &(E\otimes_{\Q} F)^{\ast} \ar[r]^{N}\ar[d] & E^{\ast} \ar[r]\ar[d] &H^1(\Gal(\bar{E}/E), A(\bar{E}) ) \ar[r]\ar@{=}[d] &0 \\
	0 \ar[r] &A(E) \ar[r] &\GL_2(E\otimes_{\Q} F) \ar[r]^{\phi} &\GSO_{3,1}(E) \ar[r] &H^1(\Gal(\bar{E}/E), A(\bar{E}) ) \ar[r] &0.}
\end{displaymath}
Therefore $\GSO_{3,1}(E) = E^{\ast} \phi(\GL_2(E\otimes_{\Q} F))$, and 
$$0 \to A(E) \to \GL_2(E\otimes_{\Q} F)\to \GSO_{3,1}(E)  \to E^{\ast}/N (E\otimes_{\Q} F)^{\ast} \to 0$$
is exact.  The above argument also applies for the adele ring $\A$ since $H^1(\Gal(\bar\Q_l/\Q_l), A(\bar\Z_l)) = 0$ for every unramified $l$, so $\GSO_{3,1}(\A) = \A^{\ast} \phi(\GL_2(\A_F))$ and we have the commutative diagram with exact columns and rows 
\begin{equation}
\xymatrix{
	&&&& 0 \ar[d] \\
	0 \ar[r] &A(\Q) \ar[r]\ar[d] &\GL_2(F) \ar[r]^{\phi}\ar[d] &\GSO_{3,1}(\Q) \ar[r]\ar[d] &\Q^{\ast}/ NF^{\ast} \ar[r]\ar[d] &0 \\
	0 \ar[r] &A(\A) \ar[r] &\GL_2(\A_F) \ar[r]^{\phi} &\GSO_{3,1}(\A) \ar[r] &\A^{\ast}/ N\A_F^{\ast} \ar[r]\ar[d] &0   \\
	&&&& \Gal(F/\Q) \ar[d] \\
	&&&& 0.
}
\end{equation}

\begin{pro}\label{HST}
	There is a bijection between cuspidal automorphic representations $\pi_H$ of $\GSO_{3,1}(\A)$ and pairs $(\pi_G, \chi)$ of a cuspidal automorphic representation $\pi_G$ of $\GL_2(\A_F)$ and a grossencharacter $\chi \colon \Q^{\ast} \backslash \A^{\ast} \to \C^{\ast}$ such that $\chi \circ N$ is the central character of $\pi_G$.
\end{pro}

\begin{proof}
	This follows directly from the above discussion (see also \cite[Proposition 1]{HST93}).  Note that $\pi_H$ corresponds to $(\{f\circ \phi \mid f \in \pi_H \} , \chi_{\pi_H})$, where $\chi_{\pi_H}$ is the central character of $\pi_H$.  For the other direction, the pair $(\pi_G, \chi)$ corresponds to the set of functions $f\colon \GSO(\Q) \backslash \GSO(\A) \to \C$ such that $f \circ \phi \in \pi$ and the central character of $f$ is $\chi$.  
\end{proof}

Now we prove \propref{GSO31}.

\begin{proof}
	Let $\pi_H$ be a cuspidal automorphic representation of $\GSO_{3,1}(\A)$ and let $\pi_G$ be the cuspidal automorphic representation of $\GL_2(\A_F)$ obtained as in the above proposition.  Following \cite[Section 3]{HST93}, the association $\pi_H \mapsto \pi_G$ is compatible with the transfer $r$, so it is also compatible with the Hecke morphism $\cH(G^S, V^S) \to \cH(H^S, U^S)$. 
	
	It remain to show the map $X_G^V \to X_H^U$ induced from $\phi$ is an isomorphism.  From the commutative diagram 
	\begin{displaymath}
	\xymatrix{
		0 \ar[r] &A(\R) \ar[r]\ar[d] &\C^{\ast} \ar[r]^{N}\ar[d] &\R^{\ast} \ar[r]\ar[d] &\R^{\ast}/ N\C^{\ast} \ar[r]\ar[d] &0 \\
		0 \ar[r] &A(\A) \ar[r] &\GL_2(\A_F) \ar[r]^{\phi} & \GSO_{3,1}(\A) \ar[r] &\A^{\ast}/ N\A_F^{\ast}  \ar[r] &0,}
	\end{displaymath}
	we deduce an exact sequence 
	\begin{displaymath}
	0\to A(\A)/A(\R) \to \GL_2(\A_F)/ \C^{\ast} \to \GSO_{3,1}(\A) /\R^{\ast} \to \A^{\ast} / (N\A_F^{\ast} \cdot \R^{\ast}) \to 0,
	\end{displaymath}
	which admits a compatible faithful left action from 
	\begin{displaymath}
	0\to A(\Q) \to \GL_2(F) \to \GSO_{3,1}(\Q) \to \Q^{\ast}/NF^{\ast}\to 0.
	\end{displaymath}
	Note that $\Q^{\ast} \backslash \A^{\ast} / (N\A_F^{\ast} \cdot \R^{\ast})$ is trivial.  Following the proof of the snake lemma, we obtain a sequence of maps 
	\begin{displaymath}
	A(\Q)\backslash A(\A) /A(\R) \hookrightarrow \GL_2(F)\backslash \GL_2(\A_{F}) / \C^{\ast} \twoheadrightarrow \GSO_{3,1}(\Q) \backslash \GSO_{3,1}(\A) /\R^{\ast} 
	\end{displaymath}
	such that the second arrow is surjective and each of its fiber is isomorphic to $A(\Q)\backslash A(\A) /A(\R)$ (note the isomorphism is not canonical in general, but here $A(\A)$ lies in the center of $\GL_2(\A_{F})$ so it's canonical).  Now consider the compatible right action from 
	\begin{displaymath}
	0\to \prod_l A(\Z_l) \to \prod_l V_l \to \prod_l U_l \to 0.
	\end{displaymath}
	It's easy to see that $A(\Q)\backslash A(\A) /(A(\R) \cdot \prod_l A(\Z_l) )$ is trivial, so the induced map $X_G^V \to X_H^U$ is an isomorphism.
\end{proof}

{}

\end{document}